%% file: ms.tex
\documentclass[11pt]{article}

\RequirePackage[OT1]{fontenc}
\RequirePackage{amsthm,amsmath,amsfonts,amssymb}
\RequirePackage[numbers]{natbib}
\RequirePackage[colorlinks,citecolor=blue,urlcolor=blue]{hyperref}
\usepackage[dvipsnames]{xcolor}
\usepackage[title,titletoc,toc]{appendix}
\usepackage{latexsym}
\usepackage{graphicx}
\usepackage{mathrsfs}
\usepackage{enumerate}
\usepackage{mathabx}
\input{04-pre-am-short.tex}

\numberwithin{equation}{section}

\newtheorem{theorem}{Theorem}
\newtheorem{lemma}{Lemma}
\newtheorem{proposition}{Proposition}
\newtheorem{corollary}{Corollary}

\def\argmin{\mathop{\rm arg\, min}}
\def\real{\mathop{{\rm I}\kern-.2em\hbox{\rm R}}\nolimits}

\def\Pen{\hbox{\rm Pen}}
\def\dPen{\hbox{\rm $\dot{\hbox{\rm P}}$en}}

\def\RE{\hbox{\rm RE}}
\def\drho{\dot\rho}

\def\RCIF{\hbox{\rm RCIF}}

\title{Sorted Concave Penalized Regression}
   \author{Long Feng\hspace{.2cm}\\
   	Department of Biostatistics,  Yale University\\
   	and \\
   Cun-Hui Zhang\thanks{
\textit{Research partially supported by NSF Grants IIS-1407939, 
	DMS-1513378, DMS-1721495, IIS-1741390 and National Institute on Drug Abuse Grant R01 DA016750.}}
\\
Department of Statistics and Biostatistics, Rutgers University}
\date{}

\begin{document}
		\maketitle 
\begin{abstract}
The Lasso is biased. Concave penalized least squares estimation (PLSE) takes advantage of signal strength to reduce this bias, leading to sharper error bounds in prediction, coefficient estimation and variable selection. For prediction and estimation, the bias of the Lasso can be also reduced by taking a smaller penalty level than what selection consistency requires, but such smaller penalty level depends on the sparsity of the true coefficient vector. The sorted $\ell_1$ penalized estimation (Slope) was proposed for adaptation to such smaller penalty levels. However, the advantages of concave PLSE and Slope do not subsume each other. We propose sorted concave penalized estimation to combine the advantages of concave and sorted penalizations. We prove that sorted concave penalties adaptively choose the smaller penalty level and at the same time benefits from signal strength, especially when a significant proportion of signals are stronger than the corresponding adaptively selected penalty levels. A local convex approximation, which extends the local linear and quadratic approximations to sorted concave penalties, is developed to facilitate the computation of sorted concave PLSE and proven to possess desired prediction and estimation error bounds. We carry out a unified treatment of penalty functions in a general optimization setting, including the penalty levels and concavity of the above mentioned sorted penalties and mixed penalties motivated by Bayesian considerations. Our analysis of prediction and estimation errors requires the restricted eigenvalue condition on the design, not beyond, and provides selection consistency under a required minimum signal strength condition in addition. Thus, our results also sharpens existing results on concave PLSE by removing the upper sparse eigenvalue component of the sparse Riesz condition. 

\end{abstract}



%\newpage
\section{Introduction}

%The purpose of this paper is to study prediction, coefficient estimation, and variable selection properties of concave penalized least squares estimator (PLSE) in linear regression under the restrictive eigenvalue (RE) condition on the design matrix. 

The purpose of this paper is twofold. First, we provide a unified treatment of 
prediction, coefficient estimation, and variable selection properties 
of concave penalized least squares estimation (PLSE) in high-dimensional linear regression 
under the restrictive eigenvalue (RE) condition on the design matrix. 
Second, we propose sorted concave PLSE to combine the advantages of concave and 
sorted penalties, and to prove its superior theoretical properties and computational feasibility 
under the RE condition. 
Along the way, we study penalty level and concavity of multivariate penalty functions, including 
mixed penalties motivated by Bayesian considerations as well as sorted and separable penalties. 
Local convex approximation (LCA) is proposed and studied as a solution for the computation of 
sorted concave PLSE. 

%As modern technology generates tons of data, high-dimensional data has been studied intensively both in statistics and computer science. In many statistical applications (biomedical, finance, engineering, etc.), the number  of variables $p$ can be much larger than the sample size $n$, but the number of important variables may be still smaller than $n$.  In such cases, we seek an approach that could select important variables, estimate coefficients and make prediction. In linear regression, a widely used approach for these goals is the penalized least-square estimation (PLSE).

Consider the linear model
\bel{eq-1-1}
\by=\bX\bbeta^*+\bep,
\eel
where $\bX=(\bx_1,...,\bx_p)\in\mathbb{R}^{n \times p}$ is a design matrix, 
$\by\in \mathbb{R}^n$ is a response vector, 
$\bep \in \mathbb{R}^n$ is a noise vector, %term follows  $\mathcal{N}\left(0,\sigma^2I_{n\times n}\right)$ 
and $\bbeta^*\in \mathbb{R}^p$ is an unknown coefficient vector. 
For simplicity, we assume throughout the paper that the design matrix is column normalized 
with $\|\bx_j\|_2^2=n$. 

Our study focuses on local and approximate solutions for the minimization of 
penalized loss functions of the form 
\bel{eq-1-2}
\|\by-\bX\bbeta\|_2^2/(2n)+ \Pen(\bbeta)
\eel
with a penalty function $\Pen(\cdot)$ satisfying certain minimum penalty level and maximum 
concavity conditions as described in Section \ref{sec-penalties}. 
The PLSE can be viewed as a statistical choice among local minimizers of the penalized loss.

Among PLSE methods, the Lasso \citep{TibshiraniR96} with the $\ell_1$ penalty 
$\Pen(\bbeta) =\lam\|\bbeta\|_1$ is the most widely used and extensively studied. 
The Lasso is relatively easy to compute as it is a convex minimization problem, 
but it is well known that the Lasso is biased. 
A consequence of this bias is the requirement of a 
neighborhood stability/strong irrepresentable condition on the design matrix $\bX$ 
for the selection consistency of the Lasso \citep{MeinshausenB06,ZhaoY06, Tropp06, Wainwright09}. 
\citet{FanL01} proposed a concave penalty to remove the bias of the 
Lasso and proved an oracle property for one of the local minimizers of the resulting penalized loss. 
\citet{Zhang10} proposed a path finding algorithm PLUS for concave PLSE and proved 
the selection consistency of the PLUS-computed local minimizer under a rate optimal 
signal strength condition on the coefficients 
and the sparse Riesz condition (SRC) \cite{ZhangH08} on the design.  
%including the $p\gg n$ setting. 
The SRC, which requires bounds on both the lower and upper sparse eigenvalues 
of the Gram matrix and is closely related to the restricted isometry property (RIP) \cite{CandesT05}, 
is substantially weaker than the strong irrepresentable condition. 
This advantage of concave PLSE over the Lasso has since become well understood. 

For prediction and coefficient estimation, the existing literature somehow presents an opposite story. 
Consider hard sparse coefficient vectors satisfying $|\supp(\bbeta^*)| \le s$ 
with $\log(p/s)\asymp \log p$ and small $(s/n)\log p$. 
Although rate minimax error bounds were proved under the RIP and SRC respectively 
for the Dantzig selector and Lasso in \citep{CandesT07} and \citep{ZhangH08}, 
\citet{BickelRT09} sharpened their results by weakening the RIP and SRC to the RE condition, 
and \citet{VanB09} proved comparable prediction and $\ell_1$ estimation 
error bounds under an even weaker compatibility or $\ell_1$ RE condition. 
Meanwhile, rate minimax error bounds for concave PLSE still 
require two-sided sparse eigenvalue conditions like the SRC 
\citep{Zhang10, Zhang10-multistage,WangLZ14,fan2015tac} or a proper known upper bound for the $\ell_1$ 
norm of the true coefficient vector \citep{loh2015regularized}. 
It turns out that the difference between the SRC and RE conditions are quite significant as 
\citet{RudelsonZ13} proved that the RE condition is a consequence of 
a lower sparse eigenvalue condition alone.  
This seems to suggest a theoretical advantage of the Lasso, in addition to its relative 
computational simplicity, compared with concave PLSE.

Emerging from the above discussion, 
an interesting question is whether the RE condition alone on the design matrix is also sufficient for the 
above discussed results for concave penalized prediction, coefficient estimation and variable selection, 
provided proper conditions on the true coefficient vector and the noise. 
An affirmative answer of this question, which we provide in this paper, amounts to the removal of the upper 
sparse eigenvalue condition on the design matrix and actually also a relaxation of the lower sparse 
eigenvalue condition or the restricted strong convexity (RSC) condition \cite{negahban2012unified} 
imposed in \cite{loh2015regularized}; 
or equivalently, to the removal of the remaining analytical advantage of the Lasso 
as far as error bounds for the afore mentioned aims are concerned. 
Specifically, we prove that when the true $\bbeta$ is sparse, 
concave PLSE achieves rate minimaxity in prediction and coefficient estimation under 
the RE condition on the design. 
Furthermore, the selection consistency of concave PLSE is also guaranteed under 
the same RE condition and an additional uniform signal strength condition on the nonzero coefficients, 
and these results also cover non-separable multivariate penalties imposed on the vector $\bbeta$ 
as a whole, including sorted and mixed penalties such as the spike-and-slab Lasso 
\citep{RockovaG16}.

In addition to the above conservative prediction and estimation error bounds 
for the concave PLSE that are comparable with those for the Lasso in both rates and regularity 
conditions on the design, we also prove faster rates for concave PLSE when the signal is partially strong. 
For example, instead of the prediction error rate $(s/n)\log p$ in the worst case scenario, 
the prediction rate for concave PLSE is actually $\sigma^2(s + s_1\log p)/n$ where $s_1$ is the number of 
small nonzero signals under the same RE condition on the design. 
Thus, concave PLSE adaptively benefits from signal strength with no harm to the performance in 
the worst case scenario where all signals are just below the radar screen. 
This advantage of concave PLSE is known under the sparse Riesz and 
comparable conditions, but not under the RE condition as presented in this paper. 

The bias of the Lasso can be also reduced by 
taking a smaller penalty level than those required for variable selection consistency, regardless 
of signal strength. 
%A smaller or adaptive penalty level may also reduce the bias of the Lasso. 
In the literature, PLSE is typically studied in a standard setting at penalty level 
$\lam\ge \lam_*=(\sigma/\eta)\sqrt{(2/n)\log p}$. 
This lower bound has been referred to as the universal penalty level.
However, as the bias of the Lasso is proportional to its penalty level, 
rate minimaxity in prediction and coefficient estimation requires 
smaller $\lam \asymp \sigma\sqrt{(2/n)\log(p/s)}$ \cite{SunZ13,BellecLT16}. 
Unfortunately, this smaller penalty level depends on $s=\|\bbeta^*\|_0$, which is typically unknown. 
For the $\ell_1$ penalty, a remedy for this issue is to apply the Slope or a Lepski type procedure \cite{SuC16,BellecLT16}.  
However, it is unclear from the literature whether the same can be done with 
concave penalties.  

We propose a class of sorted concave penalties to combine the advantages of 
concave and sorted penalties. This extends the Slope beyond $\ell_1$. 
Under an RE condition, we prove that the sorted concave PLSE inherits the benefits of both 
concave and sorted PLSE, 
namely bias reduction through signal strength and adaptation to the smaller penalty level.  
This provides prediction and $\ell_2$ estimation error bounds of the order 
$\sigma^2(s+s_1\log(p/s))/n$ and comparable $\ell_1$ estimation error bounds. 
Moreover, our results apply to approximate local solutions which can be viewed as 
output of computational algorithms for sorted concave PLSE. 

To prove the computational feasibility of our theoretical results in polynomial time, 
we develop an LCA algorithm for a large class of multivariate concave PLSE 
to produce approximate local solutions to which our theoretical results apply. 
The LCA is a majorization-minimization (MM) algorithm and is closely related 
to the local quadratic approximation (LQA) \cite{FanL01} 
and the local linear approximation (LLA) \cite{ZouL08} algorithms. 
The development of the LCA is needed 
as the LLA does not majorize sorted concave penalties in general. 
Our analysis of the LCA can be viewed as extension of the results in 
\cite{Zhang10-multistage, HuangZ12, negahban2012unified,agarwal2012fast,WangLZ14,loh2015regularized,fan2015tac}
where separable penalties are considered, typically at larger penalty levels.

The rest of this paper is organized as follows. 
In Section \ref{sec-penalties} we study penalty level and concavity of general multivariate penalties 
in a general optimization setting, 
including separable, multivariate mixed and sorted penalties, 
and also introduce the LCA for sorted penalties. 
In Section \ref{sec:concave plse}, we develop a unified treatment of 
prediction, coefficient estimation and variable selection properties of concave PLSE 
under the RE condition at penalty levels required for variable selection consistency. 
In Section \ref{sec:small penalty} we provide error bounds for approximate solutions 
at smaller and sorted penalty levels and output of LCA algorithms. 
%Section \ref{sec:discussion} contains some discussion.
%Section \ref{sec:discussion} contains discussion.

%\noindent
{\bf Notation:} 
We denote by $\bbeta^*$ the true regression coefficient vector, 
$\bSigmabar=\bX^T\bX/n$ the sample Gram matrix, 
$\calS=\supp(\bbeta^*)$ the support set of the coefficient vector, 
$s=|\calS|$ the size of the support, 
and $\Phi(\cdot)$ the standard Gaussian cumulative distribution function. 
For vectors $\bv=(v_1,...,v_p)$, we denote by $\|\bv\|_q=\sum_j(|v_j|^q)^{1/q}$ the $\ell_q$ norm, 
with $\|\bv\|_{\infty}=\max_j|v_j|$ and $\|\bv\|_0=\#\{j:v_j\neq 0\}$. Moreover, $x_+=\max(x,0)$.
%For functions $f(x)$ and $g(x)$, $f(x)\succsim g(x)$ means that $f(x) \le cg(x)$ for some constants $c \in (0,\infty)$. For matrices $M$ and $N$, we write $\|M\|_{\infty}=\max_{i,j}|M|_{ij}$ to denote the $\ell_\infty$-norm of $M$, we write $M\odiv N$ indicate element-wise division.  

\section{Penalty functions}\label{sec-penalties}
We consider minimization of penalized loss 
\bel{general-penalized-loss}
L(\bbeta) + \Pen(\bbeta),\quad \bbeta\in\R^p, 
\eel
with a general Fr\'echet differentiable loss function $L(\bbeta)$ and 
a general multivariate penalty function $\Pen(\bbeta)$ satisfying 
certain minimum penalty level and maximum concavity conditions. 

Penalty level and concavity of univariate penalty functions are well understood 
as we will briefly describe in our discussion of separable penalties 
in Subsection \ref{subsec-separable-pen} below. 
However, for multivariate penalties, we need to carefully define their penalty level 
and concavity in terms of sub-differential. This is done in Subsection \ref{subsec-sub-diff}. 
We then study in separate subsections three types of penalty functions, 
namely separable penalties, multivariate mixed penalties, and sorted penalties. 
Moreover, we develop the LCA for sorted penalties 
in Subsection \ref{subsec-4-3}. 

\subsection{Sub-differential, penalty level and concavity}\label{subsec-sub-diff} 
The sub-differential of a penalty $\Pen(\cdot)$ at a point $\bb\in\R^p$, 
denoted by $\pa\Pen(\bb)$ as a subset of $\R^p$, 
can be defined as follows. A vector $\bg\in\R^p$ belongs $\pa\Pen(\bb)$ iff  
\bel{sub-diff}
\liminf_{t\to 0+}\,t^{-1}\big\{\Pen(\bb + t\bu) - \Pen(\bb)\big\} \ge \bg^T\bu,\quad \forall\ \bu\in\R^p. 
\eel
As $\bg^T\bu$ is continuous in $\bg$, 
$\pa\Pen(\bb)$ is always a closed convex set. 

Suppose $L(\bb)$ is everywhere Fr\'echet differentiable with derivative ${\dot L}(\bb)$. 
It follows immediately from the definition of the sub-differential in (\ref{sub-diff}) that 
\bes
\liminf_{t\to 0+}\,\frac{1}{t}\Big[\big\{L(\hbbeta+t\bu) + \Pen(\hbbeta+t\bu)\big\} 
- \big\{L(\hbbeta) + \Pen(\hbbeta)\big\}\Big] \ge 0 
\ees
for all $\bu\in\R^p$ iff $- {\dot L}(\hbbeta)\in \pa\Pen(\hbbeta)$. 
This includes all local minimizers. 
Let $\dPen(\bb)$ denote a member of $\pa\Pen(\bb)$. 
We say that $\hbbeta$ is a local solution for minimizing (\ref{general-penalized-loss}) iff 
the following estimating equation is feasible: 
\bel{KKT}
- {\dot L}(\hbbeta) = \dPen(\hbbeta). 
\eel
As (\ref{KKT}) characterizes all minimizers of the penalized loss when the penalized loss is convex, 
it can be viewed as a KarushÐKuhnÐTucker (KKT) condition. 

We define the penalty level of $\Pen(\cdot)$ at a point $\bb\in \R^p$ as 
\bel{penalty-level}
\lam(\bb) = \sup\Big[\lam: \big\{\bg_{\calS_{\bb}^c}: \bg\in \pa\Pen(\bb)\big\} 
\supseteq [-\lam,\lam]^{|\calS_{\bb}^c|}\Big]. 
\eel
This definition is designed to achieve sparsity for solutions of (\ref{KKT}).
Although $\lam(\bb)$ is a function of $\bb$ in general, it depends solely on $\Pen(\cdot)$ 
for many commonly used penalty functions. Thus, we may denote $\lam(\bb)$ by $\lam$ 
for notational simplicity. 
For example, in the case of the $\ell_1$ penalty $\Pen(\bb) = \lam\|\bb\|_1$, 
(\ref{penalty-level}) holds with $\lam(\bb)=\lam$ for all $\bb$ with $|\calS_{\bb}|<p$. 
We consider a somewhat weaker penalty level for the sorted penalty 
in Subsection \ref{subsec-concave-slope}. 

We define the concavity of $\Pen(\cdot)$ at $\bb$, 
relative to an oracle/target coefficient vector $\bbeta^o$, as 
\bel{concavity}\qquad & 
\kappabar(\bb)= \kappabar(\bb,\bbeta^o)
=\sup\big\{(\bbeta^o-\bb)^T\big(\dPen(\bb) - \dPen(\bbeta^o)\big)/\|\bb-\bbeta^o\|_2^2\big\} 
%=\inf\big\{\kappa:  \sup \bh^T\big(\dPen(\bb) - \dPen(\bb+\bh)\big)\le \kappa\|\bh\|_2^2,\ 
%\forall\,\bh\in\R^p\big\},
%\\ \nonumber &&\qquad\qquad \forall \bg\in \pa\Pen(\bb),  \tbg\in \pa\Pen(\bb+\bh)\Big\}
\eel
with the convention $0/0=0$, 
where the supreme is taken over all choices 
$\dPen(\bb)\in \pa\Pen(\bb)$ and $\dPen(\bbeta^o)\in \pa\Pen(\bbeta^o)$. 
We use $\kappabar = \kappabar(\Pen) = \sup_{\bb,\tbb}\kappabar(\bb,\tbb)$ to denote 
the maximum concavity of $\Pen(\cdot)$. 
For convex penalties, 
$- \kappabar(\bb,\bbeta^o)\|\bb-\bbeta^o\|_2^2$ is the symmetric Bregman divergence. 
A penalty function $\Pen(\bb)$ is convex if and only if $\kappabar \le 0$. 
Given %an oracle coefficient vector $\bbeta^o\in\R^p$ and 
$s\ge\|\bbeta^o\|_0$ and $\xi>0$, 
%Instead of (\ref{concavity}), 
we may consider a relaxed concavity of $\Pen(\cdot)$ at $\bb$ as 
\bel{relaxed-concavity}\quad
\kappabar_{1,2}(\bb;\xi) = \inf\big\{\kappabar_2(\bb) + (1+\xi)^2s\,%|\supp(\bb)|
\kappabar_1(\bb)\big\}, 
\eel
where infimum is taken over all nonnegative $\kappabar_1(\bb)$ and $\kappabar_2(\bb)$ satisfying 
\bel{concavity-1-2}
\bh^T\big(\dPen(\bbeta^o)-\dPen(\bb)\big) \le \kappabar_1(\bb)\|\bh\|_1^2+\kappabar_2(\bb)\|\bh\|_2^2
\eel
with $\bh = \bb-\bbeta^o$ 
for all $\dPen(\bb)\in\pa\Pen(\bb)$ and $\dPen(\bb+\bh)\in\pa\Pen(\bb+\bh)$.
This notion of concavity is more relaxed than the $\ell_2$ one in (\ref{concavity}) because 
$\kappabar_{1,2}(\bb;\xi)\le \kappabar(\bb)\wedge \kappabar(\bbeta^o)$
always holds due to the option of picking $\kappabar_1(\bb)=0$. 
The relaxed concavity is 
quite useful in our study of multivariate mixed penalties 
in Subsection \ref{subsec-mixed-pen}. 
To include more solutions for %the KKT condition 
(\ref{KKT}) and also to avoid sometimes tedious task of fully characterizing the sub-differential, 
we allow $\dPen(\bb)$ to be a member of the following ``completion" of the sub-differential,
\bel{sub-diff-2}\qquad &
\overline{\pa\Pen(\bb)} 
= \hbox{\rm convex.hull}\Big\{\lim_{t\to 0+}\hbox{\rm closure}
\Big(\hbox{ $\bigcup_{\|\bv-\bb\|_2\le t}$} \,\pa\Pen(\bv)\Big)\Big\}, 
\eel
in the estimating equation (\ref{KKT}), as long as $\pa\Pen(\bb)$ is replaced by 
the same subset in (\ref{KKT}), (\ref{penalty-level}), (\ref{concavity}), (\ref{relaxed-concavity}) 
and (\ref{concavity-1-2}). 
However, for notational simplicity, we may still use $\pa\Pen(\bb)$ to denote 
$\overline{\pa\Pen(\bb)}$.  
We may also impose an upper bound condition on the penalty level of $\Pen(\cdot)$: 
\bel{max-penalty}\qquad&
\sup\big\{\big\|\dPen(\bb)\big\|_\infty/\lam(\bb): \bb\in\R^p, \dPen(\bb)\in\pa\Pen(\bb)\big\} \le \eta_*. 
\eel
It is common to have $\eta_*=1$ although $\eta_*\ge 1$ by (\ref{penalty-level}). Without loss of 
generality, we impose the condition $\Pen({\bf 0})=0$. 

\subsection{Separable penalties}\label{subsec-separable-pen}
In general, separable penalty functions can be written as a sum of penalties on individual variables, 
$\Pen(\bb) = \sum_{j=1}^p \rho_j(b_j)$.  We shall focus on separable penalties of the form 
\bel{separable-pen}
\rho(\bb;\lam) = \hbox{ $\sum_{j=1}^p$ } \rho(b_j;\lam),
\eel
where $\rho(t;\lam)$ is a parametric family of penalties with the following properties:
\begin{enumerate}[(i)]
	\item $\rho(t;\lam)$ is symmetric, $\rho(t;\lam)=\rho(-t;\lam)$ with $\rho(0;\lam)=0$; 
	\item $\rho(t;\lam)$ is monotone, $\rho(t_1;\lam)\le \rho(t_2;\lam)$ for all $0\le t_1<t_2$; 
	\item $\rho(t;\lam)$ is left- and right-differentiable in $t$ for all $t$;
	\item $\rho(t;\lam)$ has selection property, $\drho(0+;\lam) = \lam \ge 0$;   
	\item $|\drho(t-;\lam)|\vee |\drho(t+;\lam)| \le \lam$ for all real $t$,
\end{enumerate}
where $\drho(t\pm;\lam)$ denote the one-sided derivatives. 
Condition (iv) guarantees that the index $\lam$ equals to the penalty level defined in (\ref{penalty-level}), 
and condition (v) bounds the maximum penalty level with $\eta_*=1$ in (\ref{max-penalty}).
We write $\drho(t;\lam)=x$ when $x$ is between the left- and right-derivative of 
$\rho(t;\lam)$ at $t$, including $t=0$ where $\drho(0;\lam)=x$ means $|x|\le\lam$, 
so that $\drho(t;\lam)$ is defined in the sense of (\ref{sub-diff-2}). 
By (\ref{concavity}), the concavity of $\rho(t;\lam)$ is defined as
\bel{eq-2-1}
%\kappabar(t) = 
\kappabar(t;\rho,\lam) = \hbox{$\sup_{t'>t}$}\big\{\drho(t';\lam)-\drho(t;\lam)\big\}\big/(t-t'),
\eel
where the supreme is taken over all possible choices of $\drho(t;\lam)$ and $\drho(t';\lam)$ 
between the left- and right-derivatives. Further, define the overall maximum concavity of 
$\rho(t;\lam)$ as
\bel{eq-2-2}
\kappabar(\rho) = \hbox{$\max_{t\ge 0,\lam>0}$}\,\kappabar(t;\rho,\lam). 
%\kappabar(\rho) = \hbox{$\max_\lam$}\,\kappabar(\rho,\lam),\quad 
%\kappabar(\rho,\lam)= \max_{t\ge 0}\kappabar(t;\rho,\lam). 
\eel
Because $\rho(\bb;\lam)$ is a sum in (\ref{separable-pen}), 
the closure %$\overline{\pa \rho(\bb;\lam)}$ 
in (\ref{sub-diff-2}) 
is the set of all vectors $\drho(\bb;\lam) =(\drho(b_1;\lam),\ldots,\drho(b_p;\lam))^T$, 
so that $\lam$ gives the penalty level (\ref{penalty-level}), and 
\bel{eq-2-2a}
\kappabar(\bb) = \kappabar(\bb;\rho,\lam)  \le  \max_{j\le p}\kappabar(b_j;\rho,\lam)
\eel
%$\kappabar(\rho,\lam)$ in (\ref{eq-2-2}) 
gives the concavity (\ref{concavity}) of the multivariate penalty $\rho(\bb;\lam)$. 
%The relaxed concavity is smaller. 

Many popular penalty functions satisfy conditions (i)--(v) above, including 
the $\ell_1$ penalty $\rho(t;\lam)=\lambda|t|$ for the Lasso with $\kappabar(\rho)=0$, 
the SCAD (smoothly clipped absolute deviation) penalty \cite{FanL01} with 
\bel{eq-2-3}
\rho(t;\lam)=\int_{0}^{|t|}\big\{\lam - \kappabar(x-\lam)_+\big\}_+ dx
%\kappabar(0;\rho,\lam)=\kappa/(1+\kappa)
\eel
and $\kappabar(\rho)=\kappabar$, 
and the MCP (minimax concave penalty) \cite{Zhang10} with 
\bel{eq-2-4}
\rho(t;\lam)=\int_{0}^{|t|}(\lam- \kappabar x)_+dx %,\ \kappabar(0;\rho,\lam)=\kappa. 
\eel
and $\kappabar(\rho)=\kappabar$. 
An interesting way of constructing penalty functions is to mix 
penalties $\rho(t;\lam)$ with a distribution $G(d\lam)$ and a real $r_n$ as follows, 
\bel{univ-mix-1}
\rho_G(t) = - r_n^{-1} \log\Big[\,\hbox{\Large $\int$}\, \exp\big\{ - r_n \rho(t;\lam)\big\} G(d\lam)\Big]. 
\eel
This class of mixed penalties has a Bayesian interpretation as we discuss in Subsection \ref{subsec-mixed-pen}. 
If we treat $\exp\big\{ - r_n \rho(t;\lam)\big\} G(d\lam)/\int \exp\big\{ - r_n \rho(t;x)\big\} G(dx)$ 
as conditional density of $\lam$ under a joint probability $\P_G$, we have 
\bel{univ-mix-2}\qquad
\drho_G(t) &=& \E_G\big[\drho(t;\lam)\big| t\big],\ 
\cr \lam_G &=& \drho_G(0+) = \int\lam G(d\lam)\ \hbox{ (penalty level),}
\\ \nonumber \kappabar(\rho_G) &\le& \kappabar(\rho) + r_n \sup_t \Var_G\big[\drho(t;\lam)\big| t\big], 
\eel
due to $\drho(0+;\lam)=\lam$ and $\{\drho(t_1;\lam) - \drho(t_2;\lam)\}/(t_2-t_1)\le \kappabar(\rho)$ 
for all $t_1\neq t_2$ and $\lam>0$. 
For example, if $G$ puts the entire mass in a two-point set $\{\lam',\lam''\}$, 
\bes
\kappabar(\rho_G) \le \kappabar(\rho) + r_n(\lam'-\lam'')^2/4. 
%\P_G\big\{\lam=\lam_0\big|t\big\}\P_G\big\{\lam=\lam_1\big|t\big\} 
\ees
In particular, for $\rho(t;\lam)=\lam|t|$, $r_n=n$ and two-point distributions $G$, 
(\ref{univ-mix-1}) gives the spike-and-slab Lasso penalty 
as in \cite{RockovaG16}.
%as in \cite{RockovaG16,Rockova15}}.
%\bes
%\drho_G(t) 
%= \frac{\int \drho(t;\lam) \exp\big\{-r_n\rho(t;\lam)\big\}G(d\lam)}
%{\int \exp\big\{-r_n\rho(t;\lam)\big\}G(d\lam)} = \E_G\big[\drho(t;\lam)\big| t\big]  
%\ees
%\bes
%{\ddot\rho}_G(t) = \E_G\big[{\ddot\rho}\big| t\big] -r_n \Var_G\big[\drho(t;\lam)\big| t\big]  
%\ees

%\subsection{Hierarchical penalties} 
\subsection{Multivariate mixed penalties}\label{subsec-mixed-pen}
Let $\pi(\bb|\theta)$ be a parametric family of prior density functions for $\bbeta$. 
When $\bep\sim N(0,\sigma^2\bI_{n\times n})$ with known $\sigma$ and $\theta$ is given, 
the posterior mode can be written as the minimizer of 
\bes
\|\by - \bX\bb\|_2^2\big/(2n) + \Pen_\theta(\bb) 
\ees
when $\Pen_\theta(\bb) = - (\sigma^2/n)\log(\pi(\bb|\theta)/\pi({\bf 0}|\theta))$. 
In a hierarchical Bayes model where $\theta$ has a prior distribution $\pi(d\theta)$, the 
posterior mode corresponds to 
\bel{mixed-pen}
\Pen(\bb) = - r_n^{-1} \log \hbox{$\int$} \exp\big\{ - r_n \Pen_\theta(\bb)\big\}\nu(d\theta) 
\eel
with $r_n = n/\sigma^2$ and $\nu(d\theta) = \pi({\bf 0}|\theta)\pi(d\theta)/\int \pi({\bf 0}|\theta)\pi(d\theta)$. 
This gives rise to (\ref{mixed-pen}) as a general way of mixing penalties $\Pen_\theta(\cdot)$ with suitable $r_n$. 
When $r_n=n/\sigma^2$, it corresponds to the posterior for a proper hierarchical prior if the integration 
$\int \int \exp\big\{ - r_n \Pen_\theta(\bb)\big\}\nu(d\theta) d\bb$ is finite, and an improper one otherwise. 
When $0<r_n\neq n/\sigma^2$, it still has a Bayesian interpretation 
with respect to mis-specified noise level $\sqrt{n/r_n}$ or sample size $\sigma^2r_n$. 
While $r_n=0$ leads to $\Pen(\bb)=\int \Pen_\theta(\bb)\nu(d\theta)$ as the limit at $r_n=0+$, 
the formulation does not prohibit $r_n<0$. 

For $\blam = (\lam_1,\ldots,\lam_p)^T\in [0,\infty)^p$, let 
$\rho(\bb;\blam) = \sum_{j=1}^p \rho(b_j;\lam_j)$ be a separable penalty function 
with different penalty levels for different coefficients $b_j$, 
where $\rho(t;\lam)$ is a family of penalties indexed by penalty level $\lam$ 
as discussed in Subsection \ref{subsec-separable-pen}. 
As in (\ref{mixed-pen}), 
\bel{mixed-rho}
\rho_\nu(\bb) = - r_n^{-1} \log \hbox{$\int$} \exp\big\{ - r_n \rho(\bb;\blam)\big\}\nu(d\blam), 
\eel
with the convention $\rho_\nu(\bb) = \int \rho(\bb;\blam)\nu(d\blam)$ for $r_n=0$, 
is a mixed penalty for any probability measure $\nu(d\blam)$. 
We study below the sub-differential, penalty level and concavity of 
such mixed penalties. 

By definition, the sub-differential of (\ref{mixed-rho}) can be written as 
\bel{mixed-sub-diff}
\qquad & \pa \rho_\nu(\bb) = \bigg\{\frac{\displaystyle \hbox{$\int$} \drho(\bb;\blam)
\exp\big\{ - r_n \rho(\bb;\blam)\big\}\nu(d\blam)}
{\displaystyle \hbox{$\int$} \exp\big\{ - r_n \rho(\bb;\blam)\big\}\nu(d\blam)}: 
\drho(\bb;\blam)\in \pa \rho(\bb;\blam)\bigg\}
\eel
with $\pa \rho(\bb;\blam)$ being the set of all vectors 
$\drho(\bb;\blam)=(\drho(b_1;\lam_1),\ldots,\drho(b_p;\lam_p))^T$,  
provided that the $\liminf$ can be taken under the integration over $\nu(d\blam)$. 
This is allowed when $\|\drho(t;\lam)\|_\infty<\infty$. 
%We shall omit the details. 
As in (\ref{univ-mix-2}), we may write (\ref{mixed-sub-diff}) as 
\bel{mixed-sub-diff-2}\qquad &
\pa \rho_\nu(\bb) = \big\{ \E_{\nu}\big[ \drho(\bb;\blam)\big|\bb\big]: 
\drho_j(\bb;\blam) = (\drho(b_1;\lam_1),\ldots,\drho(b_p;\lam_p))^T \big\}, 
%\drho_j(\bb;\blam) = \drho(b_j;\lam_j)\ \forall j\le p \big\}, 
\eel
where the conditional $\P_\nu\big[d\blam\big|\bb\big]$ is proportional to 
$\exp\big\{ - r_n \rho(\bb;\blam)\big\}\nu(d\blam)$. We recall that 
$\drho(0;\lam_j)$ may take any value in $[-\lam_j,\lam_j]$. 

\begin{proposition}\label{prop-mixed-pen} 
Let $\rho_\nu(\bb)$ be a mixed penalty in (\ref{mixed-rho}) generated from   
a family of penalties $\rho(t;\lam)$ satisfying conditions (i)--(v) in Subsection~\ref{subsec-separable-pen}. 
Let $\calS_{\bb}=\supp(\bb)$ with $s_{\bb} = |\calS_{\bb}| < p$. Then, 
the concavity of $\rho_\nu(\bb)$ satisfies  
\bel{prop-mix-pen-1}
\kappabar(\bb)\le \kappabar(\rho)
+ \hbox{$\sup_{\bu}$}\,\phi_{\max}\big(r_n\Cov_{\nu}\big(\drho(\bu;\blam),\drho(\bu;\blam)\big|\bu\big)\big)
\eel
with $\phi_{\max}$ being the largest eigenvalue, 
and (\ref{concavity-1-2}) holds with 
\bes\quad 
\kappabar_2(\bb)\le \kappabar(\rho),\  
\kappabar_1(\bb)\le (r_n\vee 0)
\hbox{ $\sup_{\bu}$}\,\hbox{$\max_{1\le j\le p}$}\,
\Var_\nu\big(\drho(u_j;\lam_j)\big|\bu\big). 
\ees
If the components of $\blam$ are independent given $\theta$, then (\ref{concavity-1-2}) holds with 
\bes\quad 
& \kappabar_2(\bb)\le \kappabar(\rho) + (r_n\vee 0)
\hbox{ $\sup_{\bu}$}\,\hbox{$\max_{1\le j\le p}$}\,
\E_\nu\big[\Var(\drho(u_j;\lam_j)|\bu,\theta)\big|\bu\big]. 
\cr & \kappabar_1(\bb)\le (r_n\vee 0)
\hbox{ $\sup_{\bu}$}\,\hbox{$\max_{1\le j\le p}$}\,
\Var_\nu\big(\E_{\nu}\big[\drho(u_j;\lam_j)\big|\bu,\theta)\big]\big|\bu\big). 
\ees 
If in addition $\blam$ is exchangeable under $\nu(d\blam)$, 
the penalty level of (\ref{mixed-rho}) is 
\bel{prop-mix-pen-2}
\lam(\bb) = \E\big[\lam_j \big|\bb\big],\quad \forall\ j\not\in \calS_{\bb}. 
\eel
\end{proposition}

Interestingly, (\ref{prop-mix-pen-1}) indicates that mixing $\rho(\bb;\blam)$ with $r_n<0$ 
makes the penalty more convex. 

For the non-separable spike-and-slab Lasso \citep{RockovaG16},
the prior is hierarchical where $\beta_j|\blam\sim (r_n\lam_j/2)e^{-|t|r_n\lam_j}$ are independent, 
$\lam_j|\theta$ are iid with $\pi(\lam_j = \lam'|\theta) = \theta = 1 - \pi(\lam_j=\lam''|\theta)$ 
for some given constants $\lam'$ and $\lam''$, and $\theta\sim \pi(d\theta)$. 
As $\pi({\bf 0}|\theta)=\{(r_n/2)(\theta\lam'+(1-\theta)\lam'')\}^p$ and 
$\pi({\bf 0})=\int\pi({\bf 0}|\theta)\pi(d\theta)$, the penalty can be written as
\bes %as in (\ref{mixed-rho})
\rho_\nu(\bb) = \frac{- 1}{r_n}\log\frac{\pi(\bb)}{\pi({\bf 0})} 
= \frac{-1}{r_n} \log \int \exp\bigg\{ - r_n \sum_{j=1}^p \lam_j|b_j| \bigg\}\nu(d\blam), 
\ees
%with $\rho(t;\lam) = \lam |t|$, 
where $\lam_j\in\{\lam',\lam''\}$ are iid given $\theta$ 
with $\nu(\lam_j=\lam'|\theta) = \theta\lam'/\{\theta\lam'+(1-\theta)\lam''\}$
and $\nu(d\theta) = \pi({\bf 0}|\theta)\pi(d\theta)/\pi({\bf 0})$. %\int\pi({\bf 0}|\theta)\pi(d\theta)$. 
The penalty level is given by 
\bes
\lam(\bb) 
= \frac{\int \lam_\theta\exp\big\{-r_n\sum_{j\in \calS_{\bb}}\rho_\theta(b_j)\big\}\nu(d\theta)}
{\int \exp\big\{-r_n\sum_{j\in \calS_{\bb}}\rho_\theta(b_j)\big\}\nu(d\theta)},
%\ j\not\in \calS_{\bb},
\ees
with 
$\rho_\theta(t)=\{\theta \lam'e^{-r_n\lam'|t|}+(1-\theta)\lam''e^{-r_n\lam''|t|}\}/\{\theta\lam'+(1-\theta)\lam''\}$ 
and $\lam_\theta=\{\theta(\lam')^2+(1-\theta)(\lam'')^2\}/\{\theta\lam'+(1-\theta)\lam''\}$, 
and the relaxed concavity in (\ref{relaxed-concavity})--(\ref{concavity-1-2}) are bounded by 
$\kappabar_2(\bb) = 0,\ \kappabar_1(\bb) \le r_n(\lam'-\lam'')^2/4$.
%\chz{We demonstrate the dependence of the penalty level and concavity of $\rho_{\nu}(\bb)$ 
%on the magnitude and sparsity of $\bb$ in the following figures, with $\pi(d\theta)\sim$ unif (0,1).} 

\subsection{Sorted concave penalties}\label{subsec-concave-slope} 
Given a sequence of sorted penalty levels $\lam_1\ge \lam_2 \ge\cdots\ge\lam_p\ge 0$, 
the sorted $\ell_1$ penalty \citep{SuC16} is defined as 
\bel{slope}
\Pen(\bb) =\,\hbox{$\sum_{j=1}^p$} \lam_jb_j^\#, 
\eel
where $b_j^\#$ is the j-th largest value among $|b_1|,\ldots,|b_p|$. 

Here we extend the sorted penalty beyond $\ell_1$. 
Given a family of univariate penalty functions $\rho(t;\lam)$ 
and a vector $\blam=(\lam_1,\ldots,\lam_p)^T$ with non-increasing nonnegative elements, 
we define the associated sorted penalty as 
\bel{concave-slope}
\rho_\#(\bb;\blam) = \hbox{$\sum_{j=1}^p$} \rho(b_j^\#;\lam_j). 
\eel
Although (\ref{concave-slope}) seems to be a superficial extension of (\ref{slope}), 
it brings upon potentially significant benefits and its properties are nontrivial. 
We say that the sorted penalty is concave if $\rho(t;\lam)$ is concave in $t$ in $[0,\infty)$. 
In Section~\ref{sec:small penalty}, we prove that under an RE condition, 
the sorted concave penalty inherits the benefits of both the concave and sorted penalties, 
namely bias reduction for strong signal components
and adaptation to the penalty level to the unknown sparsity of $\bbeta$. 

The following proposition gives penalty level and an upper bound for the 
maximum concavity for a broad class sorted concave penalties, including the sorted SCAD penalty 
and MCP. In particular, the construction of the sorted penalty does not increase the maximum concavity 
in the class. 

\begin{proposition}\label{prop-slope}
Let $\rho_\#(\bb;\blam)$ be as in (\ref{concave-slope}) with $\lam_1\ge \cdots\ge\lam_p\ge 0$. 
Suppose $\rho(t;\lam) = \int_0^{|t|}\drho(x;\lam)dx$ with a certain $\drho(x;\lam)$  
non-decreasing in $\lam$ almost everywhere in positive $x$. Let $\calS_{\bb}=\supp(\bb)$ and 
$s_{\bb}=|\calS_{\bb}|$. 
Then, the sub-differential of $\rho_\#(\bb;\blam)$ includes all vectors $\bg=\drho_\#(\bb;\blam)$ satisfying 
\bel{prop-slope-1}
\begin{cases} g_{k_j} = \drho(b_{k_j};\lam_j),\ |b_{k_1}|\ge\cdots\ge|b_{k_s}|>0,\ j\le s_{\bb}, 
\cr |g_{k_j}|\le \lam_j, \{k_{s+1},\ldots,k_p\}=\calS_{\bb}^c,\ j>s_{\bb}. 
\end{cases}
\eel
Moreover, the maximum concavity of $\rho_\#(\bb;\blam)$ is no greater than that of the penalty 
family $\rho(t;\lam)$: 
\bes
\bh^T\Big\{\drho_\#(\bb;\blam) - \drho_\#(\bb+\bh;\blam)\Big\} \le \kappabar(\rho)\|\bh\|_2^2,\quad 
\forall\ \bb,\,\bh. 
\ees
\end{proposition} 

The monotonicity condition on $\drho(x;\lam)$ holds for the $\ell_1$, SCAD and MCP. 
It follows from (\ref{prop-slope-1}) that the maximum penalty level at each index $j\in\calS_{\bb}^c$ is 
$\lam_{s+1}$. 
Although the penalty level does not reach $\lam_{s+1}$ simultaneously 
for all $j\in\calS_{\bb}^c$ as in (\ref{penalty-level}), we
still take $\lam_{s+1}$ as the penalty level for the sorted penalty 
$\rho_\#(\bb;\blam)$. This is especially reasonable when $\lam_j$ decreases slowly in $j$. 
In Subsection \ref{subsec-4-6}, %Section~\ref{sec:small penalty}, 
we show that this weaker version of the penalty level 
is adequate for Gaussian errors provided that for certain $A_0 >1>\alpha$ 
\bel{sorted-lam}
&& \lam_j \ge \lam_{*,j} = A_0\sigma\sqrt{(2/n)\log(p/(\alpha j))},\ \ j=1,\ldots,p.
\eel
More important, (\ref{prop-slope-1}) shows sorted penalties automatically pick penalty level 
$\lam_{s+1}$ from the sequence $\{\lam_j\}$ without requiring the knowledge of $s$.

A key element of the proof of Proposition \ref{prop-slope} is to write (\ref{concave-slope}) as 
\bel{slope-2}
\qquad & \rho_\#(\bb;\blam) 
= \max\Big\{\hbox{$\sum_{j=1}^p$} \rho(b_j;\lam_{k_j}): (k_1,\ldots,k_p)^T\in \hbox{\rm perm}(p)\Big\}, 
\eel
where perm$(p)$ is the set of all vectors generated by permuting %the elements of 
$(1,\ldots,p)^T$

\subsection{Local convex approximation}\label{subsec-4-3}
We develop here LCA for penalized optimization (\ref{general-penalized-loss}), 
especially for sorted penalties. 
As a majorization-minimization (MM) algorithm, 
it is closely related to and in fact very much inspired by 
the LQA \cite{FanL01} 
and LLA \cite{ZouL08, Zhang10-multistage, HuangZ12}. 
%the local quadratic approximation (LQA) algorithm \cite{FanL01} 
%and the local linear approximation (LLA) algorithm \cite{ZouL08, Zhang10-multistage, HuangZ12}. 

Suppose for a certain continuously differentiable convex function $\Pen_-(\bb)$, 
\bel{eq-2-5-12}
\Pen_+(\bb) = \Pen(\bb)+\Pen_-(\bb)
\eel
is convex. The LCA algorithm can be written as 
\bel{eq-2-5-13}
&& \bb^{(new)}=\argmin_{\bb}\bigg\{L(\bb)+\Pen_+(\bb) - \bb^T\dPen_-(\bb^{(old)})\bigg\}.
\eel
This LCA is clearly an MM-algorithm: As  
\bes
\Pen^{(new)}(\bb) = \Pen_+(\bb) - \Pen_-(\bb^{(old)}) - (\bb - \bb^{(old)})^T\dPen_-(\bb^{(old)})
\ees
is a convex majorization of $\Pen(\bb)$ 
with $\Pen^{(new)}(\bb^{(old)}) =\Pen(\bb^{(old)})$, 
\bel{eq-2-5-14}
L(\bb^{(new)}) + \Pen(\bb^{(new)}) 
&\le& L(\bb^{(new)}) + \Pen^{(new)}(\bb^{(new)})
\cr &\le& L(\bb^{(old)}) + \Pen^{(new)}(\bb^{(old)})
\\ \nonumber &=& L(\bb^{(old)}) + \Pen(\bb^{(old)}).
\eel

Let $\rho_\#(\bb;\blam)$ be the sorted concave penalty in (\ref{concave-slope}) 
with a penalty family $\rho(t;\lam)$ and a vector of sorted penalty levels 
$\blam = (\lam_1,\ldots,\lam_p)^T$. 
Suppose $\drho(x;\lam)=(\pa/\pa x)\rho(x;\lam)$ is non-decreasing in $\lam$ almost everywhere 
in positive $x$, so that Proposition \ref{prop-slope} applies. 
Suppose for a certain continuously differentiable convex function $\rho_-(t)$ 
\bel{eq-2-5-15}
&& \hbox{$\rho_+(t;\lam_j) = \rho(t;\lam_j) + \rho_-(t)$ is convex in $t$ for $j=1,\ldots, p$.}
\eel
By (\ref{slope-2}), $\rho_{+,\#}(\bb;\blam)=\rho_\#(\bb;\blam)+\rho_-(\bb)$, 
the sorted penalty with $\rho_+(t;\lam)$, is convex in $\bb$,
so that the LCA algorithm for $\rho_\#(\bb;\blam)$ can be written as 
\bel{eq-2-5-16}
&& \bb^{(new)}=\hbox{$\argmin_{\bb}$}\big\{L(\bb)
%+\rho_\#(\bb;\blam)+\rho_-(\bb)
+\rho_{+,\#}(\bb;\blam) - \bb^T\drho_-(\bb^{(old)})\big\},
\eel
where $\drho_-(\bb)$ is the gradient of $\rho_-(\bb)=\sum_{j=1}^p\rho_-(b_j)$. 
The simplest version of LCA takes $\rho_-(t)=t^2\kappabar(\rho)/2$ with the maximum concavity 
defined in (\ref{eq-2-2}), but this is not necessary as (\ref{eq-2-5-15}) is only required to hold 
for the given $\blam$. 

\begin{figure*}[ht]\label{lca-fig-1}
	\centering
	\vspace{-.42in}
	\includegraphics[width=\textwidth]{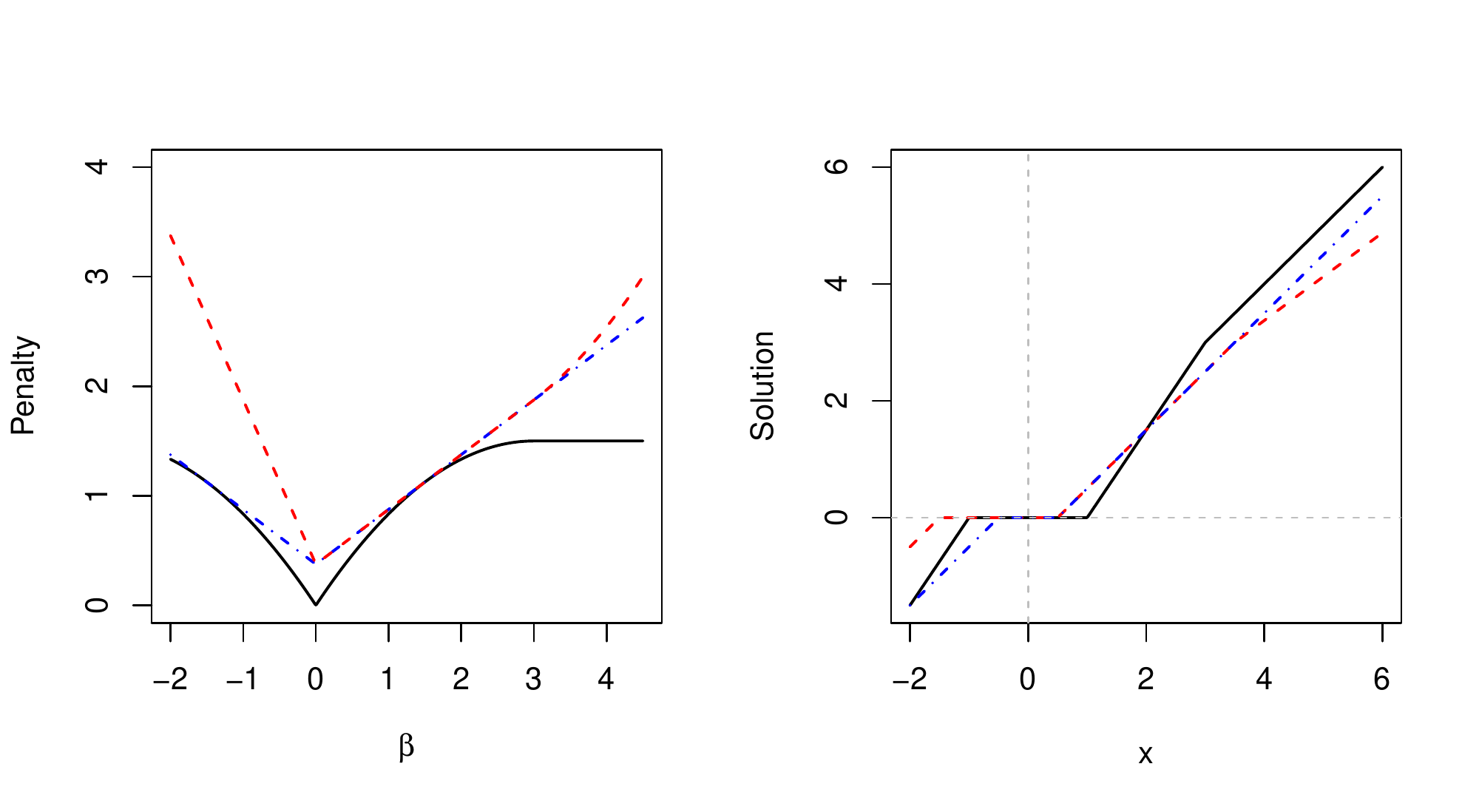}
	\vspace{-.45in}
	\caption{Local convex approximation (red dashed), local linear approximation (blue mixed) 
	and original penalty (black solid) for MCP with $\lambda=1$ and $\kappabar=1/3$ at $b^{(old)}=1.5$. 
	Left: penalty function and its approximations; Right: $\argmin_b\big\{(x-b)^2/2+\Pen(b)\big\}$.}
	\vspace{-.15in}
\end{figure*}

Figure 1 demonstrates that for $p=1$, 
the LCA with $\rho_-(t)=t^2\kappabar(\rho)/2$ also majorizes 
the LLA with $\Pen^{(new)}(b)=\rho(|b^{(old)}|;\lam)+\drho(|b^{(old)}|;\lam)(|b|-|b^{(old)}|)$. 
With $\rho_-(t)= \lam|t| - \rho(t;\lam)$ %and $\rho_+(t)=\lam |t|$ 
in (\ref{eq-2-5-15}), the LCA is identical to an unfolded LLA with 
$\Pen^{(new)}(b)=\lam|b|+\{\drho(b^{(old)};\lam)-\lam\,\sgn(b^{(old)})\big\}(b-b^{(old)})$. 
The situation is the same for separable penalties, i.e. $\lam_1=\lam_p$.  
However, the LLA is not feasible for truly sorted concave penalties with $\lam_1>\lam_p$.  
%as convexity is not preserved under sorting. 
As the LCA also majorized the LLA, it imposes larger penalty on solutions with larger step size 
compared with the LLA, but this has little effect in our theoretical analysis in 
Subsections \ref{subsec-4-5} and \ref{subsec-4-6}. 

%\chz{Long: Would you please plot the LCA versus LLA for the SCAD and MCP and write 
%captions and comments? 
%Please plot both $\rho_+(t;\lam)$ and $\rho_+(t;\lam) - \kappa_{\max}(b_0^2/2+ t-b_0)$ 
%against $\rho(t;\lam)$ and their derivatives, 
%where $b_0$ is the initial value. 
%Find Fan-Li (01) and ZouLi (08) plots as examples. 
%Note that the LCA majorization is not (cannot be) symmetric in $t$. 
%} 

%\bes
%&& L(\bx) + \big(\bb - \bx\big)^T\nabla L(\bx) + \|\bb-\bx\|_2^2/(2t)
%+\rho_{+,\#}(\bb;\blam) - \bb^T\drho_-(\bb^{(old)})
%\cr && \bb^T\{t \nabla L(\bx) - t \drho_-(\bb^{(old)})\} + \|\bb-\bx\|_2^2/2 + t\rho_{+,\#}(\bb;\blam) 
%\cr && \|\bb-\bx +\{t \nabla L(\bx) - t \drho_-(\bb^{(old)})\} \|_2^2/2 +t\rho_{+,\#}(\bb;\blam) 
%\cr && \|\bb-\{\bx - t \nabla L(\bx) + t \drho_-(\bb^{(old)})\} \|_2^2/2 +t\rho_{+,\#}(\bb;\blam) 
%\ees
The LCA (\ref{eq-2-5-16}) can be computed by
proximal gradient algorithms \citep{nesterov2007gradient, beck2009fast, parikh2013proximal}, 
which approximate $L(\bb)$ by $L(\bx) + \big(\bb - \bx\big)^T\nabla L(\bx) + \|\bb-\bx\|_2^2/(2t)$ 
around $\bx$. For example, the ISTA \cite{beck2009fast} for LCA can be written as follows. 
\begin{center}
\vspace{-.12in}
\begin{tabular}{ll}\\ \hline
\multicolumn{2}{l}{{\bf Algorithm 1:} ISTA for LCA}  \\ \hline 
Initialization: & $\bb^0 = \bb^{(old)}$ \\ 
Iteration: & $\bb^{k+1} = {\rm prox}\big(\bb^{k} - t_* \nabla L(\bb^{k}) + t_*\drho_-(\bb^{(old)}); 
t_* \rho_{+,\#}(\cdot;\blam)\big)$ \\ 
%Output: & $\bb^{(new)}$ \\ 
\hline
\end{tabular}
\end{center}
where $\rho_{+,\#}(\bb;\blam)\big) = \sum_{j=1}^p \rho_+(b^\#_j;\lam_j)$, 
$t_*$ is the reciprocal of a Lipschitz constant for $\nabla L$ or determined in the iteration by backtracking, and 
\bel{prox}
{\rm prox}\big(\bx; \Pen) = \hbox{$\argmin_{\bb}$} \big\{\|\bb - \bx\|_2^2/2 + \Pen(\bb)\big\}
\eel
is the so called proximal mapping for convex $\Pen$, e.g. $\Pen(\bb)=t \rho_{+,\#}(\bb;\blam)$. 
We may also apply FISTA \cite{beck2009fast} as an accelerated version of Algorithm 1. 
\begin{center}
\vspace{-.12in}
\begin{tabular}{ll}\\ \hline
\multicolumn{2}{l}{{\bf Algorithm 2:} FISTA for LCA}  \\ \hline 
Initialization: & $\bx^1 = \bb^0 = \bb^{(old)}$, $t_1=1$ \\ 
Iteration: & $\bb^{k} = {\rm prox}\big(\bx^{k} - t_* \nabla L(\bx^{k}) + t_*\drho_-(\bb^{(old)}); 
t_* \rho_{+,\#}(\cdot;\blam)\big)$ \\ 
& $t_{k+1} = \{1+(1+4t_k^2)^{1/2}\}/2$ \\
& $\bx^{k+1} = \bb^k+\{(t_k-1)/t_{k+1}\}(\bb^k - \bb^{k-1})$ \\ 
%Output: & $\bb^{(new)} = \bb^{k_{fin}}$ \\ 
\hline
\end{tabular}
\vspace{.1in}
\end{center}

%For separable penalties of the form $\Pen(\bb) = \sum_{j=1}^p \rho(b_j;\lam_j)$, 
%the proximal mapping (\ref{prox}) is also separable with 
%$b_j = \argmin_b \{ (x_j-b)^2/2+\rho(b;\lam_j)\}$ as elements of the solution. 
For sorted penalties $\rho_\#(\bb;\blam)$, the proximal mapping is not separable but still preserves 
the sign and ordering in absolute value of the input. Thus, after removing the sign and sorting the 
input and output simultaneously, it can be solved with 
the isotonic proximal mapping, 
\bel{iso-prox}
&& {\rm iso.prox}\big(\bx; \Pen) = \hbox{$\argmin_{\bb}$}
\Big\{\|\bb - \bx\|_2^2/2 + \Pen(\bb): 
b_j\downarrow\,\hbox{in}\,j\Big\},
\eel
with $\Pen(\bb) = \sum_{j=1}^p \rho(b_j;\lam_j)$. 
Moreover, similar to the computation of the proximal mapping for the Slope in \cite{bogdan2015slope}, 
this isotonic proximal mapping can be computed by the following algorithm. 

\begin{center}
\vspace{-.12in}
\begin{tabular}{ll}\\ \hline
\multicolumn{2}{l}{{\bf Algorithm 3:} ${\rm iso.prox}\big(\bx; \rho(\cdot;\blam)\big)$}\\ \hline 
{\bf Input:} & $\blam\downarrow$, $\bx \downarrow$  \\ 
{\bf Compute} & $b_j = \argmin_b \{ (x_j-b)^2/2+\rho(b;\lam_j)\}$\\
\multicolumn{2}{l}{{\bf While} $\bb$ is not nonincreasing {\bf do}}\\
&Identify blocks of violators of the monotonicity constraint,\\
&\qquad $b_{j'-1}>b_{j'}\le b_{j'+1}\le \cdots\le b_{j''}>b_{j''+1}$, $b_{j'}<b_{j''}$\\
&Replace $b_j$, $j'\le j\le j''$, with the solution of\\ 
&\qquad $\argmin_b \sum_{j=j'}^{j''}\{(x_j-b)^2/2+\rho(b;\lam_j)\}$ \\
\hline
\end{tabular}
\vspace{.1in}
\end{center}
We formally state the above discussion in the following proposition. 

\begin{proposition}\label{prop-prox} 
For $\bv=(v_1,\ldots,v_p)^T$, let  
$\bv^\# = (v^\#_1,\ldots,v^\#_p)^T$ with $v_j^\#$ being the $j$-th largest among 
$\{|v_1|,\ldots,|v_p|\}$. 
For $\lam_1\ge\ldots\ge\lam_p\ge 0$, 
let $\rho_\#(\bb;\blam) = \sum_{j=1}^p \rho(b^\#_j;\lam_j)$ as in (\ref{concave-slope}) 
and $\rho(\bb;\blam) = \sum_{j=1}^p \rho(b_j;\lam_j)$. Then, 
$\sgn(b_j)=\sgn(x_j)$, $|b_j|\ge |b_k|$ whenever $|x_j| > |x_k|$, and 
\bel{prop-prox-1}
\big\{{\rm prox}\big(\bx; \rho_\#(\cdot;\blam)\big)\big\}^\# 
= {\rm iso.prox}\big(\bx^\#; \rho(\cdot;\blam)\big).
\eel
Moreover, when $x_1\ge\ldots\ge x_p\ge 0$ and $(x_j-b)^2/2+\rho(b;\lam_j)$ are convex in 
$b\ge 0$ for all $j$, ${\rm iso.prox}\big(\bx; \rho(\cdot;\blam)\big)$ is solved by Algorithm 3.  
\end{proposition}

%Proposition \ref{prop-prox} assets that after removing the sign and sorting the absolute 
%value simultaneously for the input and output, 
%the proximal mapping (\ref{prox}) can be solved 
%by the isotonic proximal mapping (\ref{iso-prox}) for the corresponding unsorted penalty. 

\begin{figure*}%[ht]
	\vspace{-.4in}
	\centering
	\includegraphics[width=\textwidth]{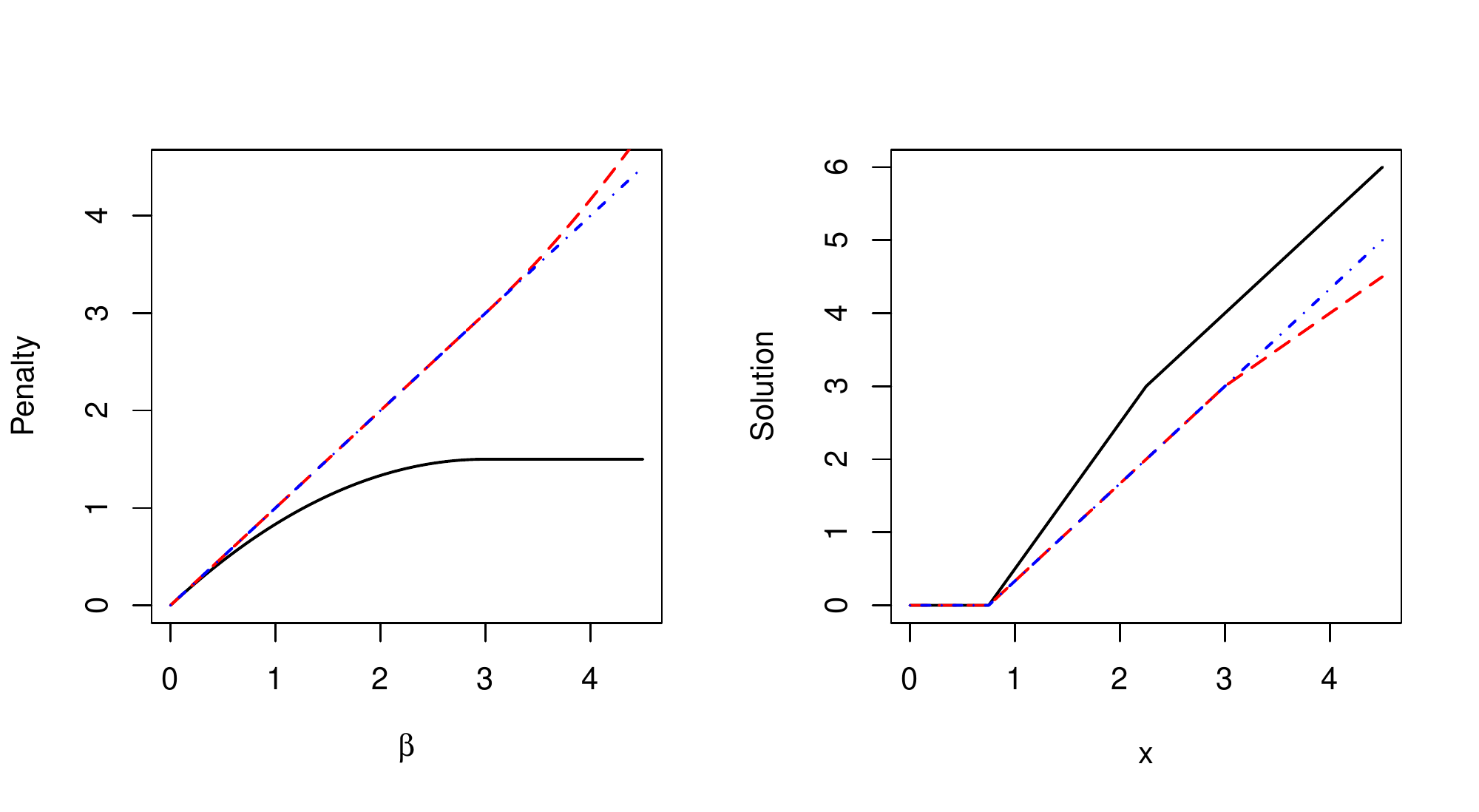}
	\vspace{-.43in}
	\caption{$\rho_+(\cdot;\lambda)=\rho(t;\lam)+\kappabar t^2/2$ for LCA (red dashed), 
	Lasso (blue mixed) and MCP (black solid) 
	with $\lambda=1$ and $\kappabar=1/3$. 
	Left: penalties; Right: proximal mappings}
	\vspace{-.2in}
\end{figure*}

%\bes
%&& (b-x)^2/(2t) + \rho_+(b;\lam), x > 0, 
%\cr && x-b = t \drho_+(b;\lam) = t \max(\lam,\kappabar b)
%\cr && b = \min\big\{ (x - t\lam)_+,x/(1+t\kappabar)\big\}
%\ees
For the MCP with $\rho_+(t;\lam)=\rho(t;\lam)+\kappabar t^2/2$, the univariate solution is 
\bes
{\rm prox}\big(x; t\rho_+(\cdot;\lam)) = \sgn(x)\min\big\{(|x| - t\lam)_+, |x|/(1+t\kappabar)\},
\ees
which is a combination of the soft threshold and shrinkage estimators. 
Figure~2 plots this univariate proximal mapping for a specific $(\lam,t)$. 
Algorithm 3 is more explicitly given as follows. 
\begin{center}
\vspace{-.12in}
\begin{tabular}{ll}\\ \hline
\multicolumn{2}{l}{{\bf Algorithm 4:} $\hbox{\rm iso.prox}(\bx; t\,\rho_+(\cdot;\blam))$ 
for the LCA with sorted MCP}  \\ \hline 
{\bf Input:} & $ \blam\downarrow$, $\bx \downarrow$, $t$, $\kappabar$ \\ 
{\bf Compute} & $b_j=\min\big\{(x_j - t\lam_j)_+, x_j/(1+t\kappabar)\}$\\
\multicolumn{2}{l}{{\bf While} $\bb$ is not nonincreasing {\bf do}}\\
&Identify blocks of violators of the monotonicity constraint,\\
&\qquad $b_{j'-1}>b_{j'}\le b_{j'+1}\le \cdots\le b_{j''}>b_{j''+1}$, $b_{j'}<b_{j''}$\\
&Replace $b_j$, $j'\le j\le j''$, with common value $b$ satisfying \\ 
&\qquad $\displaystyle b = \frac{\sum_{j=j'}^{j''}\big(x_j - t\lam_j I\{ \lam_j>\kappabar b\}\big)}
{\sum_{j=j'}^{j''}\big(1 + t\kappabar I\{ \lam_j\le \kappabar b\}\big)}$ \\
\hline
\end{tabular}
\vspace{.1in}
\end{center}
The computation of the isotonic proximal mapping for SCAD can be 
carried out in a similar but more complicated fashion, as the region 
for shrinkage is broken  
by an interval for soft thresholding in each coordinate in the domain.

\section{Properties of Concave PLSE}\label{sec:concave plse}

In this section, we present our results for concave PLSE at a sufficiently high penalty level 
to allow selection consistency.  
Smaller penalty levels are considered in Section~\ref{sec:small penalty}. 
We divide the section into three subsections to describe conditions on the design matrix, 
the collection of PLSE under consideration, and error bounds for prediction, coefficient estimation 
and variable selection.

\subsection{The restricted eigenvalue condition} 
We now consider conditions on the design matrix. 
The RE condition, proposed in \cite{BickelRT09}, 
is arguably the weakest available on the design to guarantee 
rate minimax performance in prediction and coefficient estimation for the Lasso.
The RE coefficient for the $\ell_2$ estimation loss can be defined as 
follows: For $\calS\subset\{1,\ldots,p\}$ and $\xi>0$,
%$\eta\in [0,1)$ and $\delta_*\in [0,1]$, 
\begin{align}\label{eq-2-5}
\RE_2(\calS;\xi) = \inf\left\{\frac{(\bu^T\bSigmabar\bu)^{1/2}}{\|\bu\|_2}: 
\|\bu_{\calS^c}\|_1< \xi \|\bu_{\calS}\|_1 \right\}
%\RE_2^2(\calS;\xi) = \inf\left\{\frac{\bu^T\bSigmabar\bu}{\|\bu\|_2^2}: 
%(1-{\eta})\|\bu_{\calS^c}\|_1\le (1+\delta_*\eta)\|\bu_{\calS}\|_1 \right\}.
\end{align}
with $\inf\emptyset=0$ for $\calS=\emptyset$. 
The RE condition refers to the property that $\RE_2(\calS;\xi)$ %$\RE_2(\calS;\xi)$ 
is no smaller than a certain positive constant for all design matrices under consideration.
For prediction and $\ell_1$ estimation, 
it suffices to impose a somewhat weaker compatibility condition \cite{VanB09}. 
The compatibility coefficient, also called $\ell_1$-RE \cite{VanB09}, is defined as
\begin{align}\label{eq-2-6}
\RE_1(\calS;\xi) = \inf\left\{\frac{(\bu^T\bSigmabar\bu)^{1/2}}{\|\bu_{\calS}\|_1/|\calS|^{1/2}}: 
\|\bu_{\calS^c}\|_1 < \xi \|\bu_{\calS}\|_1\right\}. 
%\RE_1^2(\calS;\xi) = \inf\left\{\frac{\bu^T\bSigmabar\bu|\calS|}{\|\bu_{\calS}\|_1^2}: 
%(1-{\eta})\|\bu_{\calS^c}\|_1\le (1+\delta_*\eta)\|\bu_{\calS}\|_1\right\}. 
\end{align}
In addition to the RE coefficients above, we define
a relaxed cone invertibility factor (RCIF) for prediction as 
\begin{align}\label{eq-2-7}
\RCIF_{\rm pred}(\calS;{\eta},\bw)
= \inf\left\{\frac{\|\bSigmabar\bu\|_\infty^2|\calS|}{\bu^T\bSigmabar\bu}: 
(1-\eta)\|\bu_{S^c}\|_1 < -\bw_{\calS}^T\bu_{\calS}\right\},
\end{align}
with $\eta\in [0,1)$ and a vector $\bw \in \mathbb{R}^p$, 
and a RCIF for the $\ell_q$ estimation as
\begin{align}\label{eq-2-8}
\RCIF_{{\rm est},q}(\calS;{\eta},\bw)
= \inf\left\{\frac{\|\bSigmabar\bu\|_\infty|\calS|^{1/q}}
{\|\bu\|_q}: 
(1-\eta)\|\bu_{\calS^c}\|_1 < -\bw_{\calS}^T\bu_{\calS}\right\}. 
\end{align}
The RCIF is a relaxation of the cone invertibility coefficient \cite{YeZ10} 
for which the constraint $\|\bu_{\calS^c}\|_1 < \xi \|\bu_{\calS}\|_1$ is imposed. 

The choices of $\xi$, $\eta$ and $\bw$ depend on the problem under consideration in the analysis, 
but typically we have $\|\bw\|_\infty\le (1-\eta)\xi$ so that the minimization in 
(\ref{eq-2-7}) and (\ref{eq-2-8}) is taken over a smaller cone. 
However, $\|\bw\|_2$ can be much smaller than $|\calS|^{1/2}(1-\eta)\xi$ 
with partial signal strength. Moreover, it is feasible to have $\bw_{\calS}={\bf 0}$ 
under a beta-min condition for selection consistency. 
In our analysis, we use an RE condition to prove cone membership of the estimation error of the concave PLSE and 
the RCIF to bound the prediction and coefficient estimation errors. 
%in Theorem \ref{th-1}. 
 The following proposition, 
which follows from the analysis in Section 3.2 of \cite{YeZ10},
 shows that the RCIF may provide sharper bounds than the RE does. 
 
\begin{proposition}\label{prop-1}
	%Let $\RE$, $\RCIF$ be as in (\ref{eq-2-5})-(\ref{eq-2-8}). 
	If $\,\|\bw_S\|_{\infty}\le (1-\eta)\xi$, then
	\bel{prop-1-1}
	\RCIF_{\rm pred}(\calS;{\eta},\bw) &\ge& \RE_1^2(\calS;\xi)/(1+\xi)^2,
	\cr
	\RCIF_{{\rm est},1}(\calS;{\eta},\bw)&\ge & \RE_1^2(\calS;\xi)/(1+\xi)^2,
	\\ \nonumber
	\RCIF_{{\rm est},2}(\calS;{\eta},\bw)&\ge & \RE_1(\calS;\xi)\RE_2(\calS;\xi)/(1+\xi).
	\eel
\end{proposition}

%\noindent{\bf Proof of Proposition \ref{prop-1}.} 
%Since $\|\bw_S\|_{\infty}\le (1-\eta)\xi$, we have 
%$\|\bu_{\calS^c}\|_1\le \|\bw_{\calS}\|_\infty\|\bu_{\calS}\|_1/(1-\eta) \le \xi \|\bu_\calS\|_1$.
%It follows that
%\bes
%{\|\bSigmabar\bu\|_{\infty}^2|\calS|\over \bu^T\bSigmabar\bu} 
%\ge {\bu^T\bSigmabar\bu|\calS|\over \|\bu\|_1^2}\ge {\bu^T\bSigmabar\bu|\calS|\over (1+\xi)^2\|\bu_\calS\|_1^2}.
%\ees
%The first inequality of (\ref{prop-1-1}) is obtained by taking infimum in the cone 
%$\scrC(\calS;\xi)=\{\bu: \|\bu_{\calS^c}\|_1\le \xi\|\bu_{\calS}\|_1\}$.  Similarly,
%\bes
%{\|\bSigmabar\bu\|_{\infty}|\calS|\over \|\bu\|_1} 
%\ge {\bu^T\bSigmabar\bu|\calS|\over \|\bu\|_1^2}\ge {\bu^T\bSigmabar\bu|\calS|\over (1+\xi)^2\|\bu_\calS\|_1^2},
%\ees
%and
%\bes
%{\|\bSigmabar\bu\|_{\infty}|\calS|^{1/2}\over \|\bu\|_2} 
%\ge {\bu^T\bSigmabar\bu|\calS|^{1/2}\over \|\bu\|_1\|\bu\|_2}
%\ge {\bu^T\bSigmabar\bu|\calS|^{1/2}\over (1+\xi)\|\bu_\calS\|_1\|\bu\|_2},
%\ees
%give the second and third inequalities of (\ref{prop-1-1}). $\hfill\square$

\subsection{Concave PLSE}
As discussed in Subsection \ref{subsec-sub-diff}, 
local minimizers of the penalized loss (\ref{eq-1-2}) must satisfy the KKT condition 
\bel{eq-2-10}
\bX^T(\by -\bX\hbbeta)/n = \dPen\big(\hbbeta\big)
\eel
for a certain member $\dPen\big(\hbbeta\big)$ of the sub-differential $\pa\Pen(\hbbeta)$. 
We shall treat (\ref{eq-2-10}) as an estimating equation to allow somewhat more 
general solutions, including solutions with $\dPen\big(\hbbeta\big)$ in 
the completion of the sub-differential $\pa\Pen(\bb)$ at $\bb=\hbbeta$ 
as defined in (\ref{sub-diff-2}). 
Local solutions of form (\ref{eq-2-10}) include all estimators $\hbbeta$ satisfying 
\bel{eq-2-10b}
&& \begin{cases} (\lambda-\kappa_*|\hbeta_j|)_+ 
\le \sgn(\hbeta_j)\bX_j^T(\by-\bX\hbbeta)/n \le \lambda, \ \  & \hbeta_j\neq 0, 
	\cr |\bX_j^T(\by-\bX\hbbeta)/n|\le \lambda,\ \  & \hbeta_j = 0, 
\end{cases}
\eel
as solutions of (\ref{eq-2-10b}) can be constructed with separable penalties 
$\Pen(\bb) = \sum_{j=1}^p \rho_j(b_j;\lam)$ with a common penalty level $\lam$ 
and potentially different concavity satisfying $\kappabar(\rho_j)\le \kappa_*$. 
Unfortunately, the more explicit (\ref{eq-2-10b}) does not cover solutions 
with sorted penalties and some mixed penalties. 

We study solutions of (\ref{eq-2-10}) %\chzR{(\ref{eq-2-10-approx})}
by comparing them with an oracle coefficient vector $\bbeta^o$. 
We assume that for a certain sparse subset $\calS$ of $\{1,\ldots,p\}$, 
$\lam_*>0$ and $\eta\in [0,1)$, the oracle $\bbeta^o$ satisfy the following, 
\bel{eq-2-12}
\supp(\bbeta^o) \subseteq \calS,\quad 
\|\bX^T(\by-\bX^T\bbeta^o)/n\|_\infty < \eta\lam_*. 
%\|\bX_{\calS^c}^T(\by-\bX^T\bbeta^o)/n\|_\infty < \eta\lam_*. 
\eel
We may take $\bbeta^o\in\R^p$ as the true coefficient vector $\bbeta^*$ with $\calS=\supp(\bbeta^*)$, 
or the oracle LSE $\hbbeta^o$ given by 
\bel{oracle-LSE}
\hbbeta^o_{\calS} = (\bX_{\cal S}^T\bX_{\cal S})^{-1}\bX_{\cal S}^T\by,\quad 
\hbbeta^o_{\calS^c} = {\bf 0}. 
\eel
When $\bep=\by-\bX\bbeta^*\sim N({\bf 0},\bV)$ with 
$\phi_{\max}(\bV)\vee \max_{j\le p}\bx_j^T\bV\bx_j/n\le \sigma^2$, 
\bel{normality}\qquad &&
\hbox{$\by-\bX\bbeta^o\sim N({\bf 0},\bV^o)$ with 
$\phi_{\max}(\bV^o)\vee\max_{j\le p}\bx_j^T\bV^o\bx_j/n\le \sigma^2$} 
\eel
for such $\bbeta^o$, so that (\ref{eq-2-12}) holds with at least probability $1 - \sqrt{2/(\pi\log p)}$ for
\bes
\lam_* = (\sigma/\eta)\sqrt{(2/n)\log p}. 
\ees

Our analysis also applies to approximate local solutions satisfying  
\bel{eq-2-10-approx}
\bX^T(\by -\bX\hbbeta)/n = \dPen\big(\hbbeta\big) + \bnu_{approx}, 
\eel
including approximate solutions of (\ref{eq-2-10b}), 
with $\bnu_{approx}$ an approximation error satisfying $\hbeta_j(\bnu_{approx})_j\ge 0$ 
and a proper upper bound for $\|\bnu_{approx}\|_\infty$. 
For computational efficiency, many algorithms only provide approximate 
solutions. In Subsection \ref{subsec-4-1}, we consider solutions with less restrictive 
approximation errors and provide a brief discussion with some references. 

Given an oracle $\bbeta^o$ satisfying (\ref{eq-2-12}), 
we consider solutions of (\ref{eq-2-10-approx}) %(\ref{eq-2-10}) 
with penalties satisfying the following penalty level and concavity conditions: 
\bel{cond-pen}\qquad &
\lam(\bbeta^o)\ge \lam_*,\ 
\kappabar_{1,2}(\hbbeta;\xi) \le \kappa_*,\ 
\|\bw_{\calS}\|_\infty \le (1-\eta)\xi, 
\eel 
where $\lam_*$ is as in (\ref{eq-2-12}), $\kappa_*$ and $\xi$ are positive constants, and 
\bel{eq-2-13}
\bw = \big\{\dPen(\bbeta^o) + \bnu_{approx}  - \bX^T(\by-\bX\bbeta^o)/n\big\}\big/\lam(\bbeta^o). 
\eel
We recall that for any $\bb$, $\lam(\bb)$ and %$\kappabar(\bb)$ and 
$\kappabar_{1,2}(\bb;\xi)$ 
are the penalty level and relaxed concavity of $\Pen(\cdot)$ as defined in 
(\ref{penalty-level})
%, (\ref{concavity})  
and (\ref{relaxed-concavity}) respectively,
%with $\kappabar_{1,2}(\bb;\xi)\le \kappabar(\bb)$,  
and $\dPen(\bbeta^o)$ is a member of the sub-differential $\pa\Pen(\bbeta^o)$ in 
(\ref{sub-diff}) or its completion in (\ref{sub-diff-2}). 
The third inequality in (\ref{cond-pen}), assumed to hold %in (\ref{cond-pen}) 
for all choices of $\dPen(\bbeta^o)\in \pa\Pen(\bbeta^o)$ in (\ref{eq-2-13}), 
follows from (\ref{eq-2-12}) and the penalty level condition (\ref{max-penalty}) when $\xi$ satisfies 
$\eta_*\lam(\bbeta^o) +\|(\bnu_{approx})_{\calS}\|_\infty + \eta\lam_*\le (1-\eta)\xi\lam(\bbeta^o)$,
%$\eta_*\lam(\bbeta^o) \le (1-\eta)\xi\lam(\bbeta^o)-\eta\lam_*$, 
but this may not be sharp. 
%The third inequality in (\ref{cond-pen}) is 
%assumed to hold for all choices of $\dPen(\bbeta^o)\in \pa\Pen(\bbeta^o)$. 
%Under (\ref{max-penalty}), $\|\bw_{\calS}\|_\infty\le\eta_*+\eta$ by (\ref{eq-2-12}), 
%but $\xi<(\eta_*+\eta)/(1-\eta)$ is possible. 

For separable penalties $\rho(\bb;\lam)$ in (\ref{separable-pen}), $\lam(\bb)=\lam$, 
$\kappabar_{1,2}(\bb;\xi) \le\kappabar(\bb) \le \kappabar(\rho)$ and 
$\|\bw_{\calS}\|_\infty\le 1+\eta + \|(\bnu_{approx})_{\calS}\|_\infty/\lam$
under conditions (i)--(v) and (\ref{eq-2-12}). The penalty level and upper bounds for the 
concavity are given in (\ref{univ-mix-2}) for univariate mixed penalties like 
the spike-and-slab Lasso, and in multivariate mixed penalties in Proposition \ref{prop-mixed-pen}. 
These facts leads to more explicit versions of (\ref{cond-pen}).

Let $\scrB(\lam_*,\kappa_*)$ be the set of all 
approximate local solutions satisfying (\ref{eq-2-10-approx})
%local solutions of (\ref{eq-2-10}) 
with different types of penalty functions satisfying (\ref{cond-pen}); Given $\bbeta^o$, 
\bel{scrB}\qquad 
&\scrB(\lambda_*,\kappa_*) = 
\Big\{\text{$\hbbeta$: (\ref{eq-2-10-approx}) %(\ref{eq-2-10}) 
holds with some $\Pen(\cdot)$ satisfying (\ref{cond-pen})}\Big\}. 
\eel
Given a penalty satisfying (\ref{cond-pen}), $\scrB(\lambda_*,\kappa_*)$ 
contains all local minimizer of the penalized loss (\ref{eq-1-2}). 
Our theory is applicable to the subclass 
\begin{align}\label{eq-2-11}
\scrB_0(\lambda_*,\kappa_*) 
= \Big\{\hbbeta: \text{$\hbbeta$ and \textbf{0} are connected in } 
\scrB(\lambda_*,\kappa_*)\Big\}.  
\end{align}
Here $\hbbeta$ and \textbf{0} are not connected iff there exist disjoint closed sets 
$\scrB_0$ and $\scrB_1$ in $\R^p$ such that 
${\bf 0}\in\scrB_0$, $\hbbeta\in \scrB_1$ and 
$\scrB(\lambda_*,\kappa_*) \subseteq \scrB_0\cup \scrB_1$.
However, this condition will be relaxed in Proposition \ref{prop-2} below.

By definition, $\scrB_0(\lambda_*,\kappa_*)$ contains the set of all local solutions (\ref{eq-2-10})
computable by path following algorithms starting from the origin, with constraints
on the penalty and concavity levels respectively. 
This is a large class of statistical solutions as 
%It is clear that $\scrB_0(\lambda_*,\kappa_*)$ is a quite general set in the sense that 
it includes all local solutions connected to the origin regardless of the specific algorithms 
used to compute the solution and different types of penalties can be used in a single solution path. 
For example, the Lasso estimator belongs to the class 
as it is connected to the origin through the LARS algorithm 
\citep{OsbornePT00, OsbornePT00b, EfronHJT04}. 
The SCAD and MCP solutions with $\lam \ge \lambda_*$ and $\kappabar \le \kappa_*$
belong to the class if they are computed 
by the PLUS algorithm \cite{Zhang10} or by a 
continuous path following algorithm from the Lasso solution. 
As $\hbbeta= {\bf 0}$ is the sparsest solution, $\scrB_0(\lambda_*,\kappa_*)$ can be viewed as 
the sparse branch of the solution space $\scrB(\lambda_*,\kappa_*)$. 
In Proposition~\ref{prop-2}, we prove that our theory is also applicable to 
approximate local solutions (\ref{eq-2-10-approx}) 
%all local solutions of the estimating equation (\ref{eq-2-10}) 
computable through a discrete solution path in 
$\scrB(\lambda_*,\kappa_*)$ from ${\bf 0}$ with $\ell_1$ step size of order $\lam$.

\subsection{Error bounds} 
Let $\scrB(\lambda_*,\kappa_*)$ be as in (\ref{scrB}), $\scrB_0(\lambda_*,\kappa_*)$ 
as in (\ref{eq-2-11}), $\calS\supseteq \supp(\bbeta^o)$, $\xi>0$, 
$\scrC(\calS;\xi)  = \big\{\bu:\|\bu_{\calS^c}\|_1\le \xi\|\bu_{\calS}\|_1\big\}$ and 
\bel{eq-2-11a}\qquad &
\scrB_0^*(\lambda_*,\kappa_*)
%=\big\{\hbbeta \in \scrB(\lambda_*,\kappa_*): 
%\big\|(\hbbeta-\bbeta^o)_{\calS^c}\big\|_1 \le \xi\big\|(\hbbeta-\bbeta^o)_{\calS}\big\|_1\big\}
%\cup \scrB_0(\lambda_*,\kappa_*)
=\scrB_0(\lambda_*,\kappa_*)\cup \big\{\hbbeta \in  
\scrB(\lambda_*,\kappa_*): \hbbeta-\bbeta^o \in \scrC(\calS;\xi)\big\}.
\eel
Here under an RE condition on the design matrix, 
we provide prediction and coefficient estimation error bounds 
for solutions of (\ref{eq-2-10-approx})
and conditions for variable selection consistency  
for solutions of (\ref{eq-2-10}) in the set (\ref{eq-2-11a}) above.

\begin{theorem} \label{th-1a}
Suppose (\ref{eq-2-12}) holds for certain $\bbeta^o\in\R^p$ and 
$\RE_2^2(\calS;\xi) \ge \kappa_*$. 
Let $\hbbeta$ be a solution of 
(\ref{eq-2-10-approx})  %(\ref{eq-2-10}) 
in $\scrB_0^*(\lambda_*,\kappa_*)$ 
with $\|\dPen(\hbbeta)\|_\infty\le \lam$. Then, 
\bel{th-1-1}\quad&&\displaystyle
	\|\bX\hbbeta-\bX\bbeta^o\|_2^2/n 
	\le \frac{\lambar^2|\calS|}{\RCIF_{\rm pred}(\calS;{\eta},\bw)}
	\le \frac{(1+\xi)^2\lambar^2|\calS|}{\RE_1^2(\calS;\xi)}, 
\eel
with 
$\lambar = \lam+\eta\lam_*+\|\bnu_{approx}\|_\infty$ and 
the $\bw$ in (\ref{eq-2-13}), and
\begin{align}\label{th-1-2}
& \|\hbbeta - \bbeta^o\|_1 \le \displaystyle
\frac{\lambar|\calS|}{\RCIF_{{\rm est},1}(\calS;{\eta},\bw)} 
\le \frac{(1+\xi)^2\lambar|\calS|}{\RE_1^2(\calS;\xi)},
\cr & \|\hbbeta - \bbeta^o\|_2 \le\displaystyle
		\frac{\lambar|\calS|^{1/2}}{\RCIF_{{\rm est},2}(\calS;{\eta},\bw)} 
		\le \frac{(1+\xi)\lambar|\calS|^{1/2}}{\RE_1(\calS;\xi)\RE_2(\calS;\xi)}, 
\\ \nonumber & \|\hbbeta - \bbeta^o\|_q \le \displaystyle
		\frac{\lambar|\calS|^{1/q}}{\RCIF_{{\rm est},q}(\calS;{\eta},\bw)},\quad\forall q\ge 1. 
\end{align}
Moreover, if $\hbbeta$ is a solution satisfying (\ref{eq-2-10b}), then 
\bes
\big(\|\bX\hbbeta-\bX\bbeta^o\|_2^2/n\big)\vee\big(\lambar\|\hbbeta-\bbeta^o\|_1\big) 
\le \lambar^2|\calS|/\RE_1^2(\calS;\xi) 
\ees
and $\|\hbbeta-\bbeta^o\|_2^2
\le \lambar^2|\calS|/\{\RE_1(\calS;\xi)\RE_2(\calS;\xi)\}$. 
\end{theorem}

Compared with Theorem \ref{th-1a}, existing results on statistical solutions of concave PLSE 
\cite{Zhang10,WangLZ14,fan2015tac} cover only separable penalties of form (\ref{separable-pen}), 
require stronger conditions on the design (such as upper sparse eigenvalue) 
and offer less explicit error bounds. 
While the error bounds for concave penalty in Theorem \ref{th-1a} match existing ones 
for the Lasso \cite{BickelRT09,VanB09}, 
they also hold when $\RE_1(\calS;\xi)$ and $\RE_2(\calS;\xi)$ in (\ref{th-1-1}) 
and (\ref{th-1-2}) are replaced by the larger version with the constraint 
$\|\bu_{\calS^c}\|_1\le \xi\|\bu_{\calS}\|_1$ replaced by the more stringent 
$(1-{\eta})\|\bu_{\calS^c}\|_1\le -\bw_{\calS}^T\bu_{\calS}$ in their definition, 
as $-\bw_{\calS}^T\bu_{\calS}$ could be much smaller than $(1-\eta)\xi\|\bu_{\calS}\|_1$ 
when significant proportion of the signals are strong. 
We describe the benefit of concave PLSE over the Lasso in such scenarios 
in the following two theorems.

\begin{theorem} \label{th-1b}
Suppose (\ref{eq-2-12}) holds for an oracle solution $\bbeta^o$ of (\ref{eq-2-10}) and 
$\RE_2^2(\calS;\xi) \ge \kappa_*$. 
Let $\hbbeta$ be a solution of (\ref{eq-2-10}) in $\scrB_0^*(\lambda_*,\kappa_*)$. 
Then, 
\bel{th-1-4}
\phi_{\min}(\bSigmabar_{\calS,\calS}) > \kappabar_{1,2}(\hbbeta;\xi)
\ \Rightarrow\ \hbbeta = \bbeta^o.
\eel
where $\phi_{min}(\cdot)$ denotes the minimum eigenvalue for symmetric matrices. 
If a separable penalty of the form $\rho(\bb;\lam)$ is taken as in (\ref{separable-pen}), then 
\bel{th-1-3}
\hbbeta_{\calS^c}={\bf 0}\ \hbox{ and }\ \sgn(\hbeta_j)\sgn(\beta^o_j)\ge 0\ \forall\,j\in \calS, 
\eel 
and under the additional condition $\kappabar(0;\rho,\lam)<\phi_{\min}(\bSigmabar_{\calS,\calS})$, 
\bel{th-1-3a}
\sgn(\hbbeta) = \sgn(\bbeta^o). 
\eel
\end{theorem}

Theorem \ref{th-1b} provides selection consistency of concave PLSE under the restricted 
eigenvalue condition, compared with the required irrepresentable condition for the Lasso 
\cite{MeinshausenB06,Tropp06,ZhaoY06}. 
This is also new as existing results require the stronger sparse Riesz condition 
\cite{Zhang10} or other combination of lower and upper sparse eigenvalue conditions on the design 
\cite{ZhangZ12,WangLZ14,fan2015tac} for selection consistency. 
Theorem \ref{th-1b} also extends the selection consistency theory to non-separable multivariate concave penalties, 
e.g. non-separable spike-and-slab Lasso with suitable choice of tuning parameters. 

The proof of Theorem \ref{th-1b} unifies the analysis for prediction, estimation and variable selection, as 
the proof of selection consistency
is simply done by inspecting the case of $\bw_{\calS}=0$ in the prediction and estimation error 
bounds. %in Theorem \ref{th-1c} below. 
The condition $\bw_{\calS}=0$, a consequence of the feasibility of $\bbeta^o$ as an oracle solution 
of (\ref{eq-2-10}) in general, holds for the oracle LSE (\ref{oracle-LSE}) under the beta-min condition 
$\min_{j\in\calS}|\hbeta_j^o| \ge \gamma\lam$ when the separable penalty (\ref{separable-pen}) 
is used with $\drho(t;\lam) = 0$ for all $|t|\ge \gamma\lam$. 
As $\RE_2^2(\calS;\xi) \le \phi_{\min}(\bSigmabar_{\calS,\calS})$ always holds 
and $\kappabar_{1,2}(\hbbeta;\xi)\le \kappa_* \le \RE_2^2(\calS;\xi)$ 
by (\ref{cond-pen}) and the RE condition, 
the condition in (\ref{th-1-4}) holds when $\bSigmabar_{\calS,\calS^c}\neq {\bf 0}$ 
as it implies $\RE_2^2(\calS;\xi) < \phi_{\min}(\bSigmabar_{\calS,\calS})$.

\begin{theorem} \label{th-1c}
Suppose (\ref{eq-2-12}) holds for $\bbeta^o\in\R^p$ and 
$\kappa_*\le \RE_2^2(\calS;\xi)$. 
Let $\hbbeta$ and $\lam$ be as in Theorem \ref{th-1a} 
with $\kappabar_{1,2}(\hbbeta;\xi) \le (1-1/C_0) \RE_2^2(\calS;\xi)$.
Then,
\bel{th-1-2a}\quad&\displaystyle  
	\|\bX\hbbeta-\bX\bbeta^o\|_2^2/n 
	\le (C_0\lam)^2\sup_{\bu\neq 0}
	\frac{\big[\bw_{\calS}^T\bu_{\calS}-(1-\eta)\|\bu_{\calS^c}\|_1\big]_+^2}{\bu^T\bSigmabar\bu}
%\cr &\le& (C_0\lam)^2\sup_{\bu\neq 0}
%	\frac{(1+\eta)^2|S|}{\RE_1^2}
\eel
with the $\bw$ in (\ref{eq-2-13}), and for any seminorm $\|\cdot\|$ as a loss function  
\bel{th-1-2b}\quad &\displaystyle 
	\|\hbbeta - \bbeta^o\| \le C_0\lam\ \sup_{\bu\neq 0}
	\frac{\|\bu\|\big[\bw_{\calS}^T\bu_{\calS}-(1-\eta)\|\bu_{\calS^c}\|_1\big]}{\bu^T\bSigmabar\bu}. 	
\eel
\end{theorem}

\begin{corollary}\label{cor-th-1} 
Suppose conditions of Theorem \ref{th-1c} with 
$\lam = \lam_* = (\sigma/\eta)\sqrt{(2/n)\log p}$, $C_0^2/\RE_2^2(\calS;\xi)=O(1)$, and 
$\bbeta^o =\hbbeta^o$ being the oracle estimator in (\ref{oracle-LSE}).  
Let $s=|\calS|$ and $s_1=\|(\dPen(\bbeta^o)+\bnu_{approx})_{\calS}/\lam\|_2$.
Then, 
%Let $\hbbeta$ be a solution of 
%\chzR{(\ref{eq-2-10-approx})}  %(\ref{eq-2-10}) 
%in $\scrB_0^*(\lambda_*,\kappa_*)$. 
%Let $\lam = \lam_* = (\sigma/\eta)\sqrt{(2/n)\log p}$, 
%$\bbeta^o =\hbbeta^o$ be the oracle estimator in (\ref{oracle-LSE}) in Theorem~\ref{th-1c}, 
%$s=|\calS|$, and \chzR{$s_1=\|(\dPen(\bbeta^o)+\bnu_{approx})_{\calS}/\lam\|_2$.}  
%%with the convention $s_1=0$ for $\bw_{\calS}={\bf 0}$. 
%Suppose $\kappa_*\le \RE_2^2(\calS;\xi)$ and 
%$\kappabar_{1,2}(\hbbeta;\xi) \le (1-1/C_0)\RE^2(\calS;\xi)$ with $C_0^2/\RE_2^2(\calS;\xi)=O_P(1)$. 
%Then, 
\bes\displaystyle
\|\bX\hbbeta-\bX\hbbeta^o\|_2^2/n+\|\hbbeta - \hbbeta^o\|_2^2
+\|\hbbeta - \hbbeta^o\|_1^2/s
=O_P(\sigma^2/n)s_1\log p, 
\ees 
implying $\hbbeta = \hbbeta^o$ when $s_1=0$, and for the true $\bbeta^*$ 
\bel{cor-th-1-1}
&& \|\bX\hbbeta-\bX\bbeta^*\|_2^2/n+\|\hbbeta - \bbeta^*\|_2^2
+\|\hbbeta - \bbeta^*\|_1^2/s
\cr &=& O_P(\sigma^2/n)\big(s_1\log p + s\big). 
\eel
For $\bnu_{approx}={\bf 0}$ and penalties (\ref{separable-pen}) with 
$\supp(\drho(t;\lam))\subseteq[-\gamma\lam,\gamma\lam]$, 
\bel{cor-th-1-2}
s_1\le \#\{j\in\calS: |\hbeta_j^o| < \lam\gamma\}, 
\eel
\end{corollary}

Theorem \ref{th-1c} and Corollary \ref{cor-th-1} demonstrate the benefits of concave PLSE, 
as $s_1=s=|\calS|$ in (\ref{cor-th-1-1}) for the Lasso 
but $s_1$ could be much smaller than $s$ for concave penalties. 
For separable concave penalties, (\ref{cor-th-1-2}) 
holds with $\gamma = 1+1/\kappabar$ for the SCAD (\ref{eq-2-3}) and 
$\gamma = 1/\kappabar$ for the MCP (\ref{eq-2-4}).
 
For the exact Lasso solution with $\kappabar(\bb)=0$, $C_0=1$ and 
$\|\bw_{\calS}\|_\infty\le 1+\eta$, 
Theorem \ref{th-1c} yields the sharpest possible 
prediction and estimation error bounds based on the basic inequality 
$\bu^T\bSigmabar\bu + (1-\eta)\|\bu_{\calS^c}\|_1 \le (1+\eta)\|\bu_{\calS}\|_1$ 
with $\bu=(\hbbeta-\bbeta^o)/\lam$, as stated in the following corollary. 

%Theorem \ref{th-1c} also sharpens existing prediction and estimation error bounds for the Lasso 
%for which $\kappabar(\bb)=0$, $C_0=1$ and $\|\bw_{\calS}\|_\infty\le (1+\eta)=(1-\eta)\xi$. 
%For example, it is strictly sharper than the familiar prediction error bound in \cite{VanB09},
%$\|\bX\hbbeta-\bX\bbeta^o\|_2^2/n \le (1+\eta)^2\lam^2|S|/\RE_1^2(\calS;\xi)$,  
%when $\RE_1^2(\calS;\xi)< \RE_1^2(\calS;1)$.  
%In fact, the sharper prediction and coefficient estimation error bounds in 
%Corollary \ref{cor-th-1-b} below are the sharpest possible based on the basic inequality 
%$\bu^T\bSigmabar\bu + (1-\eta)\|\bu_{\calS^c}\|_1 \le (1+\eta)\|\bu_{\calS}\|_1$ 
%with $\bu=(\hbbeta-\bbeta^o)/\lam$. 

\begin{corollary}\label{cor-th-1-b}
Let $\hbbeta$ be the Lasso estimator with penalty level $\lam \ge \lam_*$. %=(\sigma/\eta)\sqrt{(2/n)\log p}$. 
If (\ref{eq-2-12}) holds for a coefficient vector $\bbeta^o\in\R^p$, then 
\bes
\frac{\|\bX\hbbeta-\bX\bbeta^o\|_2^2}{n(1+\eta)^2\lam^2}
\le \sup_{\bu\neq 0}
\frac{\psi^2(\bu)}{\bu^T\bSigmabar\bu}
= \max_{0<t<1}\frac{|\calS|(1-t)^2}{\RE_1^2(\calS;t\xi)}
\ees 
with $\psi(\bu) =\big[\|\bu_{\calS}\|_1-\|\bu_{\calS^c}\|_1/\xi\big]_+$ and 
$\xi=(1+\eta)/(1-\eta)$, 
	\bes
	\frac{\|\hbbeta - \bbeta^o\|_2}{(1+\eta)\lam} 
	\le \sup_{\bu\neq 0}
	\frac{\|\bu\|_2\psi(\bu)}{\bu^T\bSigmabar\bu}
= \max_{0<t<1}\frac{|\calS|^{1/2}(1-t)}{\RE_{1,2}^2(\calS;t\xi)} 
	\ees
	with $\RE_{1,2}(\calS;\xi) = 
	\inf_{\|\bu_{\calS^c} < \xi\|\bu_{\calS}\|_1}\bu^T\bSigmabar\bu/(\|\bu\|_2\|\bu_{\calS}\|_1/|\calS|^{1/2})$, 
	%\{\RE_2(\calS;\xi)\RE_1(\calS;\xi)\}^{1/2}\ge\RE_2^2(\calS;\xi)$, 
	and 
	\bes
	\frac{\|\hbbeta - \bbeta^o\|_1}{(1+\eta)\lam} 
	\le \sup_{\bu\neq 0}
	\frac{\|\bu\|_1\psi(\bu)}{\bu^T\bSigmabar\bu}
= \max_{0<t<1}\frac{|\calS|^{1/2}(1+t\xi)(1-t)}{\RE_{1}^2(\calS;t\xi)}.
	\ees
	\end{corollary}

As Theorems \ref{th-1a}--\ref{th-1c} deal with the same estimator 
under the same RE conditions on the design, 
they give a unified treatment of the prediction, coefficient estimation and variable selection performance 
of the PLSE, including the $\ell_1$ and concave penalties. 
For prediction and coefficient estimation, (\ref{th-1-1}) and (\ref{th-1-2}) 
match those of state-of-art for the Lasso in both the convergence rate 
and the regularity condition on the design, while (\ref{th-1-2a}), (\ref{th-1-2b}) and 
Corollary \ref{cor-th-1} demonstrate the advantages 
of concave penalization when $s_1$ is much smaller than $s$. 
Meanwhile, for selection consistency, Theorem \ref{th-1b} weakens existing conditions on the design to the same 
RE condition as required for $\ell_2$ estimation with the Lasso.
%Moreover, the prediction and estimation error bounds in Theorem \ref{th-1} (ii) and 
%Corollary~\ref{cor-th-1} directly and naturally provide selection consistency when 
%$\bw_{\calS}={\bf 0}$ or $|\calS_1|=0$. 
%More precisely, for selection consistency, Theorem \ref{th-1} (iii) requires 
%only the RE condition for (\ref{th-1-3}) 
%and mild additional eigenvalue conditions for (\ref{th-1-3a}) and (\ref{th-1-4}), 
%provided the existence of an oracular solution $\bbeta^o$ 
%with $\supp(\bbeta^o)\subseteq\calS$ or equivalently $\bw_{\calS}={\bf 0}$. 
%Note that $\kappabar(0;\rho_j,\lam) \le \kappa_*$ and 
%\chzG{$\RE^2_2(\calS;\eta,\delta^*) \le \phi_{\min}(\bX_\calS^T\bX_\calS/n)$ for all $\{\eta,\delta^*\}\subset [0,1]$} 
%by definition. 
%For concave penalties, this condition $\bw_{\calS}={\bf 0}$ can be fulfilled by the rate-optimal signal 
%strength condition $\min_{j\in\calS}|\hbeta_j^o|>\gamma\lam$ as in Corollary \ref{cor-th-1}. 
%However, the condition $\bw_{\calS}={\bf 0}$ for the Lasso requires more restrictive 
%$\ell_\infty$-type conditions such as the irrepresentable condition on the design. 
These RE-based results are significant as the existing theory for concave penalization, 
which requires substantially stronger conditions on the design,  
%such as the sparse Riesz condition in \cite{Zhang10}, 
leaves a false impression that the Lasso has a technical advantage 
in prediction and parameter estimation by requiring much weaker conditions on the design 
than the concave PLSE. 
%Moreover, compared with existing analysis, the proofs of Theorem \ref{th-1a}--\ref{th-1c} are much simpler.

%The RE condition can be viewed as nearly the mildest condition on the design matrix to guarantee the desired oracle properties of the Lasso. Based only on the RE condition, Theorem \ref{th-1} proves that any path following solution of concave PLSE starting from the origin achieves rate-optimal error bounds in prediction and coefficient estimation. Furthermore, the concave PLSE guarantees variable selection consistency as long as $\bbeta^*$ is a local solution. This completely removes the need for the restrictive $\ell_\infty$-type assumptions, such as the irrepresentable condition, required by the Lasso for selection consistency.

The following lemma, which can be viewed as a basic inequality for analyzing 
concave PLSE, is the beginning point of our analysis. 

\begin{lemma}\label{lm-1} 
Let $\lam = \lam(\bbeta^o)$
%, $\kappabar(\bb)$ and $\kappabar_{1,2}(\bb;\xi)$ be 
%the penalty and concavity levels of $\Pen(\cdot)$ 
be as in (\ref{penalty-level}), 
%, (\ref{concavity}) and (\ref{relaxed-concavity}) respectively. 
$\hbbeta$ a solution of (\ref{eq-2-10-approx}), %(\ref{eq-2-10}), 
$\bh = \hbbeta - \bbeta^o$, 
and $\bw = \big\{\dPen(\bbeta^o) + \bnu_{approx} - \bX^T(\by - \bX\bbeta^o)/n\big\}\big/\lam$. Then, 
\bes\quad&
\bh^T\bSigmabar\bh 
\le -\lam \bh^T\bw + \kappabar_1(\hbbeta)\|\bh\|_1^2+\kappabar_2(\hbbeta)\|\bh\|_2^2
\ees
%with either $\bb=\hbbeta$ or $\bb=\bbeta^o$, 
for all choices of $\{\kappabar_1(\bb),\kappabar_2(\bb)\}$ satisfying (\ref{concavity-1-2}) at $\bb=\hbbeta$.
Moreover, for $\calS\supseteq \supp(\bbeta^o)$ and a proper choice of $\dPen(\bbeta^o)\in \pa\dPen(\bbeta^o)$, 
\bel{lm-1-1}&&
\bh^T\bSigmabar\bh + \big\{\lam- \|\bX_{\calS^c}^T(\by - \bX\bbeta^o)/n\|_\infty\big\}\|\bh_{\calS^c}\|_1
%\cr &\le&\bh^T\bSigmabar\bh + \lam\|\bh_{\calS^c}\|_1 - \bh_{\calS^c}^T\bX_{\calS^c}^T(\by - \bX\bbeta^o)/n
\\ \nonumber &\le& -\lam \bh_{\calS}^T\bw_{\calS} 
+ \kappabar_1(\hbbeta)\|\bh\|_1^2+\kappabar_2(\hbbeta)\|\bh\|_2^2.
\eel
\end{lemma}

Next, we prove that the solutions in  
$\scrB_0^*(\lambda_*,\kappa_*)$ in (\ref{eq-2-11a}) and other 
approximate local solutions (\ref{eq-2-10-approx})
%local solutions of (\ref{eq-2-10}) 
in $\scrB(\lambda_*,\kappa_*)$ %\setminus \scrB_0^*(\lambda_*,\kappa_*)$ 
are separated by a gap of size $a_0\lam_*$ in the $\ell_1$ distance for some $a_0>0$. 
Consequently, $\hbbeta\in \scrB_0(\lambda_*,\kappa_*)$ implies 
$\hbbeta-\bbeta^o\in \scrC(\calS;\xi)$ for 
the $\scrB_0(\lambda_*,\kappa_*)$ in (\ref{eq-2-11}) and $\scrC(\calS;\xi)$ 
in (\ref{eq-2-11a}). 

\begin{proposition}\label{prop-2} 
Let $\kappabar_j(\bb)$ 
%$\kappabar_1(\bb)$ and $\kappabar_2(\bb)$
be as in (\ref{concavity-1-2}), 
$\scrB(\lambda_*,\kappa_*)$ as in (\ref{scrB}) 
%$\calS\supseteq \supp(\bbeta^o)$, $\xi>0$, 
%$\scrC(\calS;\xi)  = \big\{\bu:\|\bu_{\calS^c}\|_1\le \xi\|\bu_{\calS}\|_1\big\}$,  
and $\scrC(\calS;\xi)$ as in (\ref{eq-2-11a}). 
% be the set of all solutions $\hbbeta$ in $\scrB(\lambda_*,\kappa_*)$ 
%with estimation error $\hbbeta - \bbeta^o$ in the cone $\scrC(\calS;\xi)$. 
Suppose $\RE_2^2(\calS;\xi)\ge\kappa_*$. 
Let ${a}_1 =\eta - \|\bz_{\calS^c}\|_\infty/\lam_*$, 
${a}_2 = {a}_1\xi/[2(\xi+1)\{\kappabar_1(\tbbeta)+\kappabar_2(\tbbeta)\}]$ 
and ${a}_3 = {a}_1(1-\eta)\xi/\{(1-\eta)(\xi + 1)+{a}_1\}$. 
Let $\hbbeta\in \scrB(\lambda_*,\kappa_*)$ with penalty level 
$\lam$, and $\tbbeta\in \scrB(\lambda_*,\kappa_*)$ with $\tlam$. 
%Let ${a}_1 =\eta - \|\bz_{\calS^c}\|_\infty/\lam_*$, 
%${a}_2 = {a}_1\xi/[2(\xi+1)\{\kappabar_1(\tbbeta)+\kappabar_2(\tbbeta)\}]$, 
%and ${a}_3 = {a}_1(1-\eta)\xi/\{(1-\eta)(\xi + 1)+{a}_1\}$. 
Then, 
\bes 
\hbox{$\hbbeta-\bbeta^o\in \scrC(\calS;\xi)$ and 
$\big\|\tbbeta - \hbbeta\big\|_1\le \tlam a_0$ imply 
$\tbbeta-\bbeta^o\in \scrC(\calS;\xi)$,}
\ees
with $a_0 = \min\big[{a}_2, {a}_2{a}_3/\{(1-\eta)(\xi\vee 1)\}\big]$. Consequently, 
\bes
\scrB_0^*(\lambda_*,\kappa_*)=\big\{\hbbeta\in \scrB(\lambda_*,\kappa_*): 
\hbbeta-\bbeta^o\in \scrC(\calS;\xi)\big\}. 
\ees
\end{proposition}

It follows from Proposition \ref{prop-2} that our theoretical results are applicable to 
statistical choices of 
approximate local solutions (\ref{eq-2-10-approx})
%local solutions in (\ref{eq-2-10}) 
computable through a discrete 
path of solutions $\hbeta^{(t)}$ satisfying $\|\hbbeta^{(t)} - \hbbeta^{(t-1)}\|_1\le {a}_0\lam^{(t)}$ 
and beginning from ${\bf 0}$ or the Lasso solution, as $\{\bf 0\}$ and the Lasso solution 
both satisfy the condition $\hbbeta-\bbeta^o\in \scrC(\calS;\xi)$. 

%In Section \ref{sec:small penalty}, we prove that $\hbbeta$ are connected 
%to the Lasso solution in no more than $O(1)$ steps under a strong convexity condition. 
%Let $\hbbeta$ be a solution of (\ref{eq-2-10}) in $ \scrB_0^*(\lambda_*,\kappa_*)$. 

\section{Smaller penalty levels, sorted penalties and LCA}
% Concave PLSE with smaller penalty}
\label{sec:small penalty}
We have studied in Section~{\ref{sec:concave plse}} exact solutions (\ref{eq-2-10}) 
and approximate solutions (\ref{eq-2-10-approx}) 
for no smaller penalty level than $\lambda_*$ in the event where $\lam_*$ is a strict upper bound of 
the supreme norm of $\bz = \bX^T(\by-\bX\bbeta^o)/n$ as in (\ref{eq-2-12}). 
Such penalty or threshold levels are commonly used in the literature to study regularized 
methods in high-dimensional regression, but 
this is conservative and may yield poor numerical results. 
Under the normality assumption (\ref{normality}), 
(\ref{eq-2-12}) requires $\lam\ge \lam_*=(\sigma/\eta)\sqrt{(2/n)\log p}$, 
but the rate optimal penalty level is 
the smaller $\lam\ge \sigma\sqrt{(2/n)\log(p/s)}$ for prediction and coefficient estimation 
with $s=|\calS|$. 
It is known that for $\log(p/s)\ll\log p$, 
rate optimal performance in prediction and coefficient estimation  
can be guaranteed by the Lasso with the nonadaptive smaller penalty level 
depending on $s$ \cite{SunZ13,BellecLT16} or the Slope to achieve adaptation 
in the smaller penalty level \cite{SuC16,BellecLT16}.  
However, it is unclear from the literature whether the same can be done with 
concave penalties to also take advantage of noise strength, 
and under what conditions on the design.
%Moreover, for computational efficiency, many algorithms only provide approximate 
%solutions of (\ref{eq-2-10}). 
Moreover, for computational considerations, it is desirable to relax the condition 
imposed in Section~{\ref{sec:concave plse}} 
on the supreme norm of the approximation error $\bnu_{approx}$ in (\ref{eq-2-10-approx}).
In this section, we consider approximate solutions for general penalties with nonadaptive 
penalty levels which are allowed to be smaller, concave slopes, and solutions produced by 
the LCA algorithm. 
Our analysis imposes a somewhat stronger RE condition, but that is the only 
condition required on the design matrix. 

%This is the most commonly used threshold level in high-dimensional regression literature.
%Notably, \citet{SunZ13} found a smaller penalty level that may enjoy better prediction and estimation accuracy. They use a sparse $\ell_2$-norm to control the excess of $\bX^T\bep/n$ over a threshold level $\lambda_*$ when $\lam_* < \|\bX^T\bep/n\|_{\infty}$. Since their theoretical results limited to the $\ell_1$ penalty Lasso, we in this section extend their results to the concave PLSE and prove its consistency property. 

%%Besides introduce the smaller penalty in the concave PLSE, we also extend the general sparsity condition, i.e., $\supp(\bbeta^*)\in \calS$ in section 3.  Now we consider the $\bbeta^*$ satisfying the following complexity bound for certain threshold levels $\lambda_0>0$,
%% \bel{sparse2}
%% |\calS|+\sum_{j\notin \calS}|\beta^*_j|/\lam_0 \le s_*
%% \eel
%%with certain $\calS \in \{1,...,p\}$. This include the general sparsity condition with $\|\bbeta^*\|_0=s^*$, $\supp(\bbeta^*)\in\calS$.

\subsection{Approximate local solutions}\label{subsec-4-1}
As in (\ref{eq-2-10-approx}), we write 
\bel{approx-sol}
\bX^T(\by-\bX\hbbeta)/n = \dPen(\hbbeta)+\bnu_{approx}
\eel
as approximate solution of (\ref{eq-2-10}), including (\ref{eq-2-10b}) for separable penalties. 
However, instead of imposing $\ell_\infty$ bounds on the approximation error as in 
Section \ref{sec:concave plse}, we will impose in this section a more practical 
$\ell_\infty$-$\ell_2$ split bound  
%where $\bnu_{approx}$ is sufficiently small under a certain error measure, 
as in (\ref{eq-4-10}) below. 
Solutions of form (\ref{approx-sol})
are called approximate local solutions in \cite{ZhangZ12}
where their uniqueness, variable selection properties and relationship to the global solution were studied. 
Computational algorithms for approximate solutions and statistical properties of the 
resulting estimators have been considered in 
\citep{negahban2012unified,agarwal2012fast,WangLZ14,loh2015regularized,fan2015tac} among others. 
However, these studies of approximate solutions 
all focus on separable penalties with penalty level (\ref{eq-2-12}) or higher. 

Let $\bbeta^o$ be the oracle coefficient vector as in Section~{\ref{sec:concave plse}} 
and $\calS\subset \supp(\bbeta^o)$. 
Let $\bnu_{div} = \dPen(\hbbeta) - \dPen(\bbeta^o)$. 
It follows from (\ref{approx-sol}) that 
\bel{eq-4-2}
\bX^T(\by-\bX\hbbeta)/n = \dPen(\bbeta^o)+\bnu_{div}+\bnu_{approx}. 
\eel
With $\bz=\bX^T(\by-\bX\bbeta^o)/n$, (\ref{eq-4-2}) leads to the identity
\bes
\bh^T\bSigmabar \bh = \bh^T\big\{\bz - \dPen(\bbeta^o)\big\} - \bh^T\bnu_{div} - \bh^T\bnu_{approx}
\ees 
with $\bh = \hbbeta - \bbeta^o$. While upper bounds for $-\bh^T\bnu_{div}$ and 
$-\bh^T\bnu_{approx}$ can be obtained via the maximum concavity of the penalty and 
the size restriction on $\bnu_{approx}$, a favorable choice of the sub-derivative 
$\dPen(\bbeta^o)$ must be used to control the noise $\bz$. 
When the penalty function is endowed with a constant penalty level $\lam=\lam(\bbeta^o)$ on 
$\calS^c$ as in (\ref{penalty-level}) and concavity level as in (\ref{concavity}), 
\bel{eq-4-3}
&& \bh^T\bSigmabar \bh + (1-\eta)\lam\big\|\bh_{\calS^c}\big\|_1 - \kappabar(\hbbeta)\|\bh\|_2^2
\\ \nonumber &\le& \big(\bh_{\calS^c}^T\bz_{\calS^c}
 - \eta\lam\|\bh_{\calS^c}\|_1\big) - \lam \bh_{\calS}^T\bw_{\calS}- \bh^T\bnu_{approx}. 
\eel
This can be viewed as a basic inequality in our analysis of (\ref{approx-sol}). 
%to be carried out in Subsections \ref{subsec-4-5} and \ref{subsec-4-6}. 
We may use the relaxed concavity (\ref{concavity-1-2}) 
%$\kappabar_1(\hbbeta)\|\bh\|_1^2 + \kappabar_2(\hbbeta)\|\bh\|_2^2$ 
as in Lemma \ref{lm-1}, but (\ref{eq-4-3}) provides notational simplicity 
for a unified treatment with sorted concave penalties.

\subsection{Sorted concave penalized estimation}\label{subsec-4-2}
Approximate solutions for sorted concave penalties are still defined by (\ref{approx-sol}), 
so that (\ref{eq-4-2}) also holds. However, as sorted penalties are defined in 
(\ref{concave-slope}) with a sequence of penalty levels $\lam_1\ge\cdots\ge\lam_p\ge 0$, 
they do not provide $\ell_1$ control of the noise as in (\ref{eq-4-3}). 
Instead, sorted concave penalties control the noise through a sorted $\ell_1$ norm. 
More precisely, they regularize 
the correlation of the noise and design vectors in $\calS^c$ with the norm
\bel{eq-4-4}
\hbox{$\|\bb\|_{\#,s} = \sum_{j=1}^{p-s}(\lam_{s+j}/\lam_{s+1})b_j^\#,\quad \forall\ \bb\in\R^{p-s},\ s=|\calS|,$}
\eel
a standardized dual norm induced by the set of $\bg_{\calS^c}$ given in (\ref{prop-slope-1}). 
Here $b_j^\#$ is the $j$-th largest value among $|b_1|,\ldots,|b_{p-s}|$. 

As (\ref{eq-4-2}) holds, taking the most favorable $\dPen(\bbeta^o)$ and 
the concavity bound in Proposition \ref{prop-slope} in the derivation of (\ref{eq-4-3}), 
we find that 
\bel{eq-4-5}
&& \bh^T\bSigmabar \bh + (1-\eta)\lam\big\|\bh_{\calS^c}\big\|_{\#,s}- \kappabar(\hbbeta)\|\bh\|_2^2
\\ \nonumber &\le & \big(\bh_{\calS^c}^T\bz_{\calS^c}
 - \eta\lam\big\|\bh_{\calS^c}\big\|_{\#,s}\big) - \lam \bh_{\calS}^T\bw_{\calS}- \bh^T\bnu
%\\ \nonumber && + \kappabar_1(\hbbeta)\|\bh/\lam\|_\#^2 + \kappabar_2(\hbbeta)\|\bh\|_2^2
\eel
with $\calS\supseteq\supp(\bbeta^o)$, $\bz =\bX^T(\by-\bX\bbeta^o)/n$, 
$\bw = \{\dPen(\bbeta^o)-\bz\}/\lam$, $\bnu = \bnu_{approx}$, 
$\lam = \lam_{s+1}$, and $\kappabar(\hbbeta)=\kappabar(\rho)$. 
%where $\|\bb\|_\#$ is as in (\ref{slope}). 
This includes (\ref{eq-4-3}) as a special case, because $\|\bb_{\calS^c}\|_{\#,s}=\|\bb_{\calS^c}\|_1$ 
when $\lam_{s+1}=\cdots=\lam_p$. 
More important, (\ref{eq-4-5}) is adaptive to the penalty level $\lam_{s+1}$ 
without the knowledge of $s$. 

%We assume 
%\bel{eq-4-5a}
%&& \|\bh_{\calS^c}\|_{\#,s}\le \xi s^{1/2}\|\bh\|_2\ \Rightarrow\ 
%\kappabar_1(\hbbeta)\|\bh/\lam\|_\#^2 + \kappabar_2(\hbbeta)\|\bh\|_2^2
%\le \kappabar(\hbbeta;\xi)\|\bh\|_2^2
%\eel
%as in (\ref{relaxed-concavity}). This is guaranteed to hold when 
%\bes
%\kappabar(\hbbeta;\xi) \ge \kappabar_1(\hbbeta)\Big(\xi s^{1/2} + \|\lam_{1:s}\|_2/\lam_{s+1}\Big)^2 
%+ \kappabar_2(\hbbeta). 
%\ees

We outline our analysis of (\ref{eq-4-5}) as follows. 
For $r>0$ and $\gamma>0$, define 
\bel{eq-4-6}
\qquad &\displaystyle \Delta(r,\bw,\bnu) =  \sup_{\bu\neq 0} \,
 \frac{\bu_{\calS^c}^T\bz_{\calS^c}/\lam%\bX_{\calS^c}^T(\by-\bX\bbeta^o)/n 
 - \eta\|\bu_{\calS^c}\|_{\#,s} -  \bu_{\calS}^T\bw_{\calS}- \bu^T\bnu/\lam}
{r\max\big\{\|\bu\|_2,\|\bX\bu\|_2\sqrt{\gamma/n}\big\}}. 
\eel
We assume that for a certain constant $\xi>0$, 
\bel{eq-4-7}
\Delta((1-\eta)\xi |\calS|^{1/2},\bw,\bnu) <1
\eel 
for all approximate solutions (\ref{approx-sol}) under consideration. 
This and (\ref{eq-4-5}) yield 
\bel{eq-4-8}
&& \bh^T\bSigmabar \bh + (1-\eta)\lam\big\|\bh_{\calS^c}\big\|_{\#,s}- \kappabar(\hbbeta)\|\bh\|_2^2
\\ \nonumber &\le & (1-\eta)\xi |\calS|^{1/2}\lam\max\big\{\|\bh\|_2,\|\bX\bh\|_2\sqrt{\gamma/n}\big\}.   
%\\ \nonumber && + \kappabar_1(\hbbeta)\|\bh\|_\#^2 + \kappabar_2(\hbbeta)\|\bh\|_2^2
\eel
When $\kappabar(\hbbeta)\|\bh\|_2^2 \le \bh^T\bSigmabar \bh$, (\ref{eq-4-8}) provides 
the membership of $\bh$ in the cone 
\bel{eq-4-9}
&& \scrC_\#(\calS;\xi,\gamma) = \big\{\bu: \|\bu_{\calS^c}\|_{\#,s}
\le \xi|\calS|^{1/2}\max\big(\|\bu\|_2,(\gamma\bu^T\bSigmabar\bu)^{1/2}\big)\big\}. 
\eel
On the other hand, when $\bh\in \scrC_\#(\calS;\xi,\gamma)$, a suitable RE condition 
provides $\kappabar(\hbbeta)\|\bh\|_2^2 \le \bh^T\bSigmabar \bh$. 
A key step in our analysis is to break this vicious circle, which will be done in 
Subsection \ref{subsec-4-4}.  

The seemingly complicated (\ref{eq-4-6}),
which summarizes conditions in our analysis 
on the noise vector $\bep$, the penalty level $\lam$ and approximation error $\bnu$, 
is actually not hard to decipher. 
Consider $\bnu$ satisfying
\bel{eq-4-10}&&
\bu^T\bnu\le \eta_1\lam_{s+1}\|\bu_{\calS^c}\|_{\#,s}+ r_2\lam_{s+1} \|\bu\|_2\ \forall\ \bu\in\R^p
%\|\bw_{\calS}\|_2\le w,
\eel
with $\eta_1\in (0,\eta)$ and $r_2>0$. 
The above condition on $\bnu$ is fulfilled when 
\bes
\hbox{$\bnu=\bnu_1+\bnu_2$ with $\nu_{1,j}^\# \le \eta_1\lam_j$ 
and  $\|\bnu_2\|_2 + \eta_1\|\lam_{1:s}\|_2  \le \lam_{s+1} r_2$,}
\ees
as this implies 
$\bu^T\bnu \le \eta_1\sum_{j=1}^p u_j^\#\lam_j +\|\bu\|_2\|\bnu_2\|_2 
\le \eta_1\lam_{s+1}\|\bu_{\calS^c}\|_{\#,s}+(\eta_1 \|\lam_{1:s}\|_2+\|\bnu_2\|_2)\|\bu\|_2$. 
Suppose $\lam_j$ are bounded from below, $\lam_j\ge\lam_{*,j}, 1\le j\le p$, such as (\ref{sorted-lam}). Let 
\bel{eq-4-11}
&& \Deltabar(r,\eta,w,r_2) = 
\bigg\{\sum_{j=1}^{p-s}\frac{\{z^\#_j- \eta \lam_{*,s+j}\}^2}{r^2\lam_{*,s+1}^2}
+\frac{w^2}{r^2}\bigg\}^{1/2}+\frac{r_2}{r}. 
\eel
It follows from the monotonicity of (\ref{eq-4-11}) in $\lam_{*,s+j}$ and some algebra that 
\bes
\Delta(r,\bw,\bnu)I_{\|\bw_{\calS}\|_2\le w} \le \Deltabar(r,\eta-\eta_1,w,r_2)
\ees
when (\ref{eq-4-10}) holds. 
In Subsection \ref{subsec-4-6}, we derive upper bound for the median of $\Deltabar$ by combing the 
arguments in \citep{BellecLT16,SunZ13}, and then apply concentration inequality to $\Delta$. 
In the simpler case with unsorted smaller penalty, 
$\lam= \sigma L/\{n^{1/2}(\eta-\eta_1)\}$ with $L=\sqrt{2\log(p/s)}$, 
\bes
\E\,\Deltabar(r,\eta-\eta_1,w,r_2) 
&\le& r^{-1}\Big\{\sqrt{4s(\eta-\eta_1)^2/(L^4+2L^2) + w^2}+r_2\Big\}
\cr &=& o(1)+(w+r_2)/r
\ees
by Proposition 10 in \cite{SunZ13} under the normality assumption (\ref{normality}). 
Thus, (\ref{eq-4-7}) holds with high probability when 
$\|\bw_{\calS}\|_\infty\le (1-\eta)\xi'$ 
with $\xi'<\xi$, $s=|\calS|$ and $r_2/s^{1/2}$ is sufficiently small, 
in view of the third condition in (\ref{cond-pen}). 

The idea of including $\|\bX\bu\|_2$, by setting $\gamma>0$ in (\ref{eq-4-6}), 
comes from \cite{BellecLT16}. It provides upper bound $1/\big(r\lam\sqrt{\gamma n}\big)$
for the Lipschitz norm of $\Delta(r,\bw,\bnu)$ in the noise $\by-\bX\bbeta^o$, 
and thus large deviation bounds for (\ref{eq-4-6}). 

\subsection{Approximate solutions of the LCA}

As in (\ref{approx-sol}), we consider approximate solutions of the LCA (\ref{eq-2-5-13}) of the form 
\bel{eq-4-12}
{\bf 0} = {\dot L}(\bb^{(new)}) + \dPen_+(\bb^{(new)}) - \dPen_-(\bb^{(old)}) + \bnu_{approx}. 
\eel
Such approximate solutions can be viewed as output of 
iterative algorithms such as those discussed in Subsection \ref{subsec-4-3}.
The following lemma provides the LCA version of the basic inequality (\ref{eq-4-5}). 

\begin{lemma}\label{lm-LCA}
Let $\bh = \bb^{(new)} - \bbeta^o$ and $\bz = - {\dot L}(\bbeta^o)$. Then, 
\bes
&& D_L(\bb^{(new)},\bbeta^o)+D_+(\bb^{(new)},\bbeta^o) 
\cr &=& \bh^T\big\{\bz - \dPen_+(\bbeta^o) + \dPen_-(\bb^{(old)}) - \bnu_{approx}\big\}
\cr &=& \bh^T\big\{\bz - \dPen(\bbeta^o) + \bnu_{carry} - \bnu_{approx}\big\}, 
\ees 
where $D_L(\bb,\bbeta)=(\bb-\bbeta)^T\{{\dot L}(\bb)-{\dot L}(\bbeta)\}$ and 
$D_+(\bb,\bbeta)=(\bb-\bbeta)^T\{\dPen_+(\bb)-\dPen_+(\bbeta)\}$ are respectively 
the symmetric Bregman divergence for the loss $L(\bb)$ and the majorization penalty 
$\Pen^{(new)}(\bb)$ in (\ref{eq-2-5-14}), and 
$\bnu_{carry}=\dPen_-(\bb^{(old)}) - \dPen_-(\bbeta^o)$ is the carryover error in gradient. 
\end{lemma}

In linear regression, $L(\bb)=\|\by-\bX\bb\|_2^2/(2n)$ and (\ref{eq-4-12}) can be written as 
\bel{eq-4-13}&& 
\bX^T(\by - \bX\bb^{(new)})/n = \dPen_+(\bb^{(new)}) - \dPen_-(\bb^{(old)}) + \bnu_{approx}
\eel
as an approximate solution to the convex minimization problem of LCA penalized LSE. 
As in the derivation of (\ref{eq-4-2}), (\ref{eq-4-3}) and (\ref{eq-4-5}), this leads to 
\bel{eq-4-14}&& 
\bX^T(\by - \bX\bb^{(new)})/n = \dPen(\bbeta^o) + \bnu_{div} - \bnu_{carry} + \bnu_{approx},
\eel
with $\bnu_{div}=\dPen_+(\bb^{(new)}) - \dPen_+(\bbeta^o)$ and $\bnu_{carry}$ as in 
Lemma \ref{lm-LCA}, and then 
\bel{eq-4-15}
&& \bh^T\bSigmabar \bh + (1-\eta)\lam\big\|\bh_{\calS^c}\big\|_{\#,s}
\\ \nonumber &\le & \bh^T\bSigmabar \bh + (1-\eta)\lam\big\|\bh_{\calS^c}\big\|_{\#,s} 
+ D_+(\bb^{(new)},\bbeta^o) % D_+(\bbeta^o+\bh,\bbeta^o)
\\ \nonumber &\le & \big(\bh_{\calS^c}^T\bz_{\calS^c}
 - \eta\lam\big\|\bh_{\calS^c}\big\|_{\#,s}\big) - \lam \bh_{\calS}^T\bw_{\calS}- \bh^T\bnu
%\\ \nonumber && + \kappabar_1(\hbbeta)\|\bh/\lam\|_\#^2 + \kappabar_2(\hbbeta)\|\bh\|_2^2
\eel
when the favorable $\dPen(\bbeta^o)$ is taken as in (\ref{prop-slope-1}), 
where $\bh = \bb^{(new)} - \bbeta^o$,
$D_+(\bb^{(new)},\bbeta^o) = \bh^T\bnu_{div}$ is the symmetric Bregman divergence 
as in Lemma~\ref{lm-LCA} and 
$\bnu =  \bnu_{approx} - \bnu_{carry}$. 
Thus, (\ref{eq-4-5}) holds with $\kappabar(\hbbeta)=0$ and an extra carryover term 
$\bh^T\bnu_{carry}$ on the right-hand side. 

For the sorted penalty $\rho_\#(\bb;\blam)$ in (\ref{concave-slope}), we can always take $\rho_-(t)$ with 
$|(\pa/\pa t)^2\rho_-(t)|\le \kappabar(\rho)$, so that the carryover error is bounded by 
\bel{eq-4-16}
&& \big\|\bnu_{carry}\big\|_2 
\le \big\|\drho_-(\bb^{(old)}) - \drho_-(\bbeta^o)\big\|_2 
\le \kappabar(\rho)\big\|\bh^{(old)}\big\|_2 
\eel
with $\bh^{(old)} = \bb^{(old)} - \bbeta^o$ in the analysis of (\ref{eq-4-13}) through (\ref{eq-4-5}). 

\subsection{Good solutions}\label{subsec-4-4}
We define here a set of ``good solutions'' to which our error bounds 
in the following two subsections apply. 

First, as in (\ref{cond-pen}), we need to impose in addition to (\ref{eq-4-7})
a condition on the concavity of the penalty 
in (\ref{approx-sol}). This leads to the solution set
\bel{eq-4-17}\qquad 
&& \scrB(\lambda_*,\kappa_*,\gamma) 
= \Big\{\text{$\hbbeta$: (\ref{approx-sol}) and (\ref{eq-4-7}) hold},\, 
\kappabar(\hbbeta) \le \kappa_* \Big\}. 
%and (\ref{cond-pen-2}) hold}
\eel
Here $\lam_*$ is a minimum penalty level requirement implicit in (\ref{eq-4-7}), e.g.
\bel{eq-4-18}\quad &
\lam \ge \lam_* = (\eta-\eta_1)^{-1}\sigma L/n^{1/2}
\eel
with $\Phi(-L)\le s/p$, e.g. $L=\sqrt{2\log(p/s)}$ 
for fixed penalty levels or $\lam_*=\lam_{*,s+1}$ for sorted ones satisfying (\ref{sorted-lam}).
We note that $\scrB(\lambda_*,\kappa_*,\gamma)$ contains all solutions in the set 
$\scrB(\lambda_*,\kappa_*)$ in (\ref{scrB}) with $\kappabar_1(\hbbeta)=0$
although (\ref{eq-4-17}) allows approximate solutions 
with smaller minimum penalty level $\lam_*$. 

Let $\|\bu\|_{\#,*} = \|\bu_{\calS^c}\|_{\#,s}+|\calS|^{1/2}\max\big\{\|\bu\|_2, \|\bX\bu\|_2\sqrt{\gamma/n}\big\}$. 
We say that two solutions $\hbbeta$ and $\tbbeta$ in 
$\scrB_\#(\lambda_*,\kappa_*,\gamma)$ are connected by an $a_0$-chain 
%($a_0>0$, scaled to the error in norm $\|\cdot\|_{\#,*}$
if there exist $\hbbeta^{(k)}\in\scrB_\#(\lambda_*,\kappa_*,\gamma)$ 
with sorted penalty levels $\{\lam^{(k)}_1,\ldots,\lam_p^{(k)}\}$ such that 
\bel{eq-4-19}&&
\hbbeta^{(0)}=\tbbeta,\ \hbbeta^{(k^*)}=\hbbeta,\ 
\big\| \hbbeta^{(k)} - \hbbeta^{(k-1)}\big\|_{\#,*} \le a_0|\calS|\lam_{s+1}^{(k)},
\eel
$k=1,\ldots,k^*$ with the $a_0$ specified in Proposition \ref{prop-2a} below.
This condition holds if $\hbbeta$ and $\tbbeta$ are connected through a 
continuous path %$\{\hbbeta^{(t)}, 0\le t\le 1\}$ 
in $\scrB(\lambda_*,\kappa_*,\gamma)$. 
Similar to (\ref{eq-2-11}), we define 
\bel{eq-4-20}%\label{eq-2-11}
\quad &&
\scrB_0(\lambda_*,\kappa_*,\gamma)
= \big\{\hbbeta\in \scrB(\lambda_*,\kappa_*,\gamma): \text{(\ref{eq-4-19}) holds with $\tbbeta={\bf 0}$}\big\}.
%\cap \scrB(\lambda_*,\kappa_*,\gamma). 
\eel
As the solutions are connected to ${\bf 0}$ though a chain, 
$\scrB_0(\lambda_*,\kappa_*,\gamma)$ can be viewed as the sparse 
branch of the solution set $\scrB(\lambda_*,\kappa_*,\gamma)$. 

With the cone (\ref{eq-4-9}), we define a set of ``good solutions'' as follows, 
\bel{eq-4-21}\quad &&
\scrB_0^*(\lambda_*,\kappa_*,\gamma)
\\ \nonumber  &=& 
\big\{\hbbeta \in \scrB(\lambda_*,\kappa_*,\gamma): \hbbeta-\bbeta^o \in \scrC_\#(\calS;\xi,\gamma)\big\}
\cup \scrB_0(\lambda_*,\kappa_*,\gamma)
\\ \nonumber &&
\cup \big\{\text{$\hbbeta = \bb^{(new)}$: (\ref{eq-4-13}) and (\ref{eq-4-7}) hold 
with the $\bnu$ in (\ref{eq-4-15})}\big\}.
\eel
This is the set of approximate solutions (\ref{approx-sol}) 
with estimation error inside the cone or connected to the origin through a chain, 
or approximate LCA solutions 
(\ref{eq-4-13}) with
a reasonably small 
carryover component $\bh^T\bnu_{carry}$ in the $\bh^T\bnu$ in (\ref{eq-4-15}). 
Iterative applications of the LCA do provide a chain of good solutions to the final output 
without having to specify the step size. 
We note that for sorted penalties, (\ref{eq-4-16}) can be used to bound 
$\bh^T\bnu_{carry}$. 

We prove below that good solutions all belong to the cone (\ref{eq-4-9}) when the 
following restricted eigenvalue is no smaller than the $\kappa_*$ in (\ref{eq-4-17}),
\bel{eq-4-22}&& 
\RE_\#(S;\xi,\gamma) = \inf\bigg\{\frac{(\bu^T\bSigmabar\bu)^{1/2}}{\|\bu\|_2}: 
{\bf 0}\neq \bu\in \scrC_\#(S;\xi,\gamma)\bigg\}. 
\eel
As the cone $\scrC_\#(S;\xi,\gamma)$ in (\ref{eq-4-9}) depends on the sorted 
$\blam = (\lam_1,\ldots,\lam_p)^T$, the infimum in (\ref{eq-4-22}) is taken over 
all $\blam$ under consideration.
We note that
\bes&& 
\RE_\#(S;\xi,\gamma) \ge \bigg[\inf\left\{\frac{(\bu^T\bSigmabar\bu)^{1/2}}{\|\bu\|_2}: 
\|\bu_{\calS^c}\|_{\#,s} < \xi |\calS|^{1/2}\|\bu\|_2\right\}\bigg]\wedge\frac{1}{\gamma}. 
\ees
When the cone is confined to $\lam_1=\lam_p$, i.e. $\|\bu_{\calS^c}\|_{\#,s} =\|\bu_{\calS^c}\|_1$,   
the RE condition on $\RE_\#^2(S;\xi,\gamma)$ is equivalent to the 
restricted strong convexity condition \citep{negahban2012unified} as 
$\|\bu_{\calS^c}\|_1 < \xi |\calS|^{1/2}\|\bu\|_2$ implies $\|\bu\|_1< (\xi+1)|\calS|^{1/2}\|\bu\|_2$. 
Compared with (\ref{eq-2-5}), the RE in (\ref{eq-4-22}) is smaller due to the use of a larger cone. 
However, this is hard to avoid 
because 
the smaller penalty does not control the $\ell_\infty$ 
measure of the noise as in (\ref{eq-2-12}) and 
we do not wish to impose uniform bound on the approximation error $\bnu$.

\begin{proposition}\label{prop-2a} 
Let $\kappabar(\bb)$ be as in (\ref{concavity}), 
$\scrB(\lambda_*,\kappa_*,\gamma)$ as in (\ref{eq-4-17}), 
and $\scrC_\#(\calS;\xi,\gamma)$ as in (\ref{eq-4-9}). 
Suppose $\RE_\#^2(\calS;\xi,\gamma)\ge\kappa_*$ and 
\bel{eq-4-23}
\Delta((1-\eta-a_1)\xi|\calS|^{1/2},\bw,\bnu) \le 1,\quad 
0<a_1 < 1-\eta, 
\eel
as in (\ref{eq-4-7}) for all $\{\lam, \bw, \bnu\}$ 
associated with solutions in $\scrB(\lambda_*,\kappa_*,\gamma)$. 
Let ${a}_2 = {a}_1\xi^2/\{2\kappabar(\tbbeta)\}$,
${a}_3= {a}_1(1-\eta)\xi/\{(1-\eta-a_1)(\xi + 1)+{a}_1\}$ and 
$a_0 = \min\big[{a}_2, {a}_2{a}_3/\{(1-\eta)(\xi\vee 1)\}\big]$. 
Then, 
\bes
\scrB_0^*(\lambda_*,\kappa_*,\gamma) = \big\{\hbbeta\in \scrB(\lambda_*,\kappa_*,\gamma): 
\hbbeta-\bbeta^o\in \scrC_\#(\calS;\xi,\gamma)\big\}. 
\ees
%for the norm $\|\cdot\|_{\#,*}$ in (\ref{eq-4-19}), 
%\bes 
%\hbox{$\big\{\hbbeta-\bbeta^o\in \scrC_\#(\calS;\xi,\gamma), 
%\|\tbbeta - \hbbeta\|_{\#,*} \le \tlam |\calS| a_0\big\}
%\ \Rightarrow\ \tbbeta-\bbeta^o\in \scrC_\#(\calS;\xi,\gamma)$.}
%\ees
%Moreover, 
%$\scrB_0^*(\lambda_*,\kappa_*,\gamma) = \big\{\hbbeta\in \scrB(\lambda_*,\kappa_*,\gamma): 
%\hbbeta-\bbeta^o\in \scrC_\#(\calS;\xi,\gamma)\big\}$. 
%%$\scrB_0(\lambda_*,\kappa_*,\gamma)\subseteq \big\{\hbbeta\in \scrB(\lambda_*,\kappa_*,\gamma): 
%%\hbbeta-\bbeta^o\in \scrC_\#(\calS;\xi,\gamma)\big\}$. 
\end{proposition}

\subsection{Analytic error bounds}\label{subsec-4-5}
In this subsection we provide prediction and estimation error bounds for 
the set of approximate solutions in (\ref{eq-4-21}). 
While (\ref{eq-4-7}) is imposed on the entire $\scrB_0^*(\lambda_*,\kappa_*,\gamma)$, 
it can be sharpened to 
\bel{eq-4-26}
\qquad &\displaystyle 
\Delta(r_1,\bw,\bnu) \le 1
\eel
for a given $\hbbeta$ with $r_1\le (1-\eta)\xi |\calS|^{1/2}$, 
in view of our discussion below (\ref{eq-4-10}).

\begin{theorem} \label{th-2c}
Let $\scrB_0^*(\lambda_*,\kappa_*,\gamma)$ be the solution set 
given through (\ref{eq-4-17}), (\ref{eq-4-20}) and (\ref{eq-4-21}), including 
approximate solutions for concave and sorted penalties and the LCA, 
with $\calS\supseteq\supp(\bbeta^o)$. Suppose 
$\kappa_*\le \RE_\#^2(\calS;\xi,\gamma)$. 
Let $\hbbeta\in \scrB_0^*(\lambda_*,\kappa_*,\gamma)$ satisfying 
$\kappabar(\hbbeta) \le (1-1/C_0)\RE_\#^2(\calS;\xi,\gamma)$ and (\ref{eq-4-26}).
Let $F(\bu) = \max\big\{\|\bu\|_2, (\gamma \bu^T\bSigmabar\bu)^{1/2}\big\}$. 
Then, for any semi-norm $\|\cdot\|$, 
\bel{eq-4-27}
\|\bh\| \le C_0\lam\,\sup_{\bu\neq 0}\frac{\|\bu\|\{r_1F(\bu) - (1-\eta)\|\bu_{\calS^c}\|_{\#,s}\}}
{\bu^T\bSigmabar\bu}
\eel
with $\bh = \hbbeta-\bbeta^o$ for (\ref{approx-sol}) and $\bh = \bb^{(new)}-\bbeta^o$ for (\ref{eq-4-13}). 
In particular, 
\bes
\frac{\|\bh_{\calS}\|_1+\|\bh_{\calS^c}\|_{\#,s}}{(1+\xi)|\calS|^{1/2}} \le F(\bh) 
\le \frac{(\bh^T\bSigmabar\bh)^{1/2}}{\RE_\#(\calS;\xi,\gamma)} 
\le \frac{C_0r_1\lam}{\RE_\#^2(\calS;\xi,\gamma)}.
\ees
\end{theorem}

As an extension of Theorem \ref{th-1c}, Theorem \ref{th-2c} provides 
prediction and $\ell_2$ estimation error bounds in the same form along with a 
comparable sorted $\ell_1$ error bound. 
It demonstrate the benefit of sorted concave penalization as 
$r_1^2 \asymp s_1 + s/\log^2(p/s) + r_2^2$ in standard settings as described in Theorem \ref{th-6} 
and Corollary \ref{cor-th-6}, where $s=|\calS|$, $s_1$ can be viewed as the number of small nonzero 
coefficients, and $r_2$ is the $\ell_2$ component of the approximation error as in (\ref{eq-4-10}).

Theorem \ref{th-2c} applies to approximate solutions 
at smaller penalty levels, for penalties with fixed penalty level ($\lam_1=\cdots=\lam_p$), 
for sorted penalties adaptive choice of penalty 
level $\lam_{s+1}$ from assigned sorted sequence $\lam_1\ge \cdots\ge \lam_p$, 
and for general LCA with possibly sorted penalties. 
However, selection consistency of (\ref{approx-sol}) is not guaranteed 
as false positive cannot be ruled out at the smaller penalty level or with general approximation errors. 
On the other hand, as properties of the Lasso at smaller penalty levels and the Slope 
have been studied in \cite{SunZ13, SuC16, BellecLT16} among others, Theorem \ref{th-2c} 
can be viewed as an extension of their results to concave and/or sorted penalties 
discussed in Section \ref{sec-penalties} and to the LCA.  

It is also possible to derive error bounds in the case of $C_0=\infty$ as 
in Theorem \ref{th-1a} if the $\ell_\infty$ norm in the definition of the RCIF is 
replaced by the norm $\|\cdot\|_{\#,*}$ in (\ref{eq-4-19}). 
We omit details for the sake of space. 

%\chz{Long: Please add some more comments if possible, 
%e.g. a description of original results without $C_0$ for smaller penalty levels 
%or possibly in the discussion section}
%
%\lfR{Prof. Zhang: As $\|\bSigmabar\bu\|_\infty$ is not bounded in the smaller penalty case, the original results without $C_0$ introduced some new notations, e.g. a new combination norms $\|\bv\|_c=\max\big\{\|\bv_{\tilde{\calS}}\|_2/m^{1/2}, \|\bv_{\tilde{\calS^c}}\|_\infty \big\}$. This would drive the analysis too long. I wonder if we really need to add the old results? }

Next we apply Theorem \ref{th-2c} to iterative application of the LCA:  
\bel{eq-4-29}
\bb^{(t)} \leftarrow \hbox{\rm LCA}\big(\bb^{(t-1)},\Pen^{(t)},\bnu_{approx}^{(t)}\big)
\eel
where $\bb^{(new)}\leftarrow \hbox{\rm LCA}\big(\bb^{(old)},\Pen,\bnu\big)$ 
is the one-step LCA as in (\ref{eq-4-13}) with a decomposition  
$\Pen^{(t)}(\bb)=\Pen_+^{(t)}(\bb)-\Pen_-^{(t)}(\bb)$ as in (\ref{eq-2-5-12}). 

\begin{theorem}\label{th-5} 
Let $\bb^{(t)}$ generated in (\ref{eq-4-29}) with penalty levels $\lam^{(t)}=\lam^{(t)}_{s+1}$ 
for $\Pen^{(t)}$. Suppose $\|\dPen_-^{(t)}(\bb) - \dPen_-^{(t)}(\bbeta^o)\|_2\le\kappa_0\|\bb-\bbeta^o\|_2$ 
for all $\bb$ and 
\bel{eq-4-30}
\Delta\Big(r_1^{(t)},\bw^{(t)},\bnu_{approx}^{(t)}\Big) \le 1,\ t=1,\ldots,t_{fin}, 
\eel
with the $\Delta(r,\bw,\bnu)$ in (\ref{eq-4-6}) and certain $r_1^{(t)}>0$. 
Let $\bh^{(t)} = \bb^{(t)} - \bbeta^o$ and $\nu_0=\kappa_0\|\bh^{(0)}\|_2$ be the initial 
carryover error. Suppose the RE condition 
\bes
\RE_\#(\calS;\xi,\gamma)\ge \kappa_0\{\lam^{(t)}/\lam^{(t+1)}\}\{r_1^{(t)}\lam^{(1)}/\nu_0+1\}
\ees
with $(1-\eta)\xi s^{1/2} > \nu_0/\lam^{(1)} + r_1^{(t)}$, $t=1,\ldots,t_{fin}$. Then, for $t\le t_{fin}$ 
\bes
F\big(\bh^{(t)}\big) 
\le \frac{r_1^{(t)}\lam^{(t)}}{\RE_\#^2}(\calS;\xi,\gamma) + \theta_0\big\|\bh^{(t-1)}\big\|_2
\le \sum_{k=1}^t \frac{\theta_0^{t-k} r_1^{(k)}\lam^{(k)}}{\RE_\#^2(\calS;\xi,\gamma)} 
+ \theta_0^t \big\|\bh^{(0)}\big\|_2
\ees
with $F(\bu) = \max\big\{\|\bu\|_2, (\gamma \bu^T\bSigmabar\bu)^{1/2}\big\}$ and 
$\theta_0=\kappa_0/\RE_\#^2(\calS;\xi,\gamma)$. 
\end{theorem}

To find an approximate solution of PLSE with sorted penalty $\rho_\#(\bb;\blam_*)$, 
we may implement (\ref{eq-4-29}) with a fixed penalty family $\rho(x;\lam)$ 
and $\blam^{(t)} \to \blam_*$, 
\bel{LCA-1}&&
%{\dot L}(\bb^{(t)}) + \drho_\#(\bb^{(t)};\blam^{(t)}) + \drho_-(\bb^{(t)};\blam^{(t)}) - 
\bX^T(\by-\bX\bb^{(t)})/n = \drho_{+,\#}(\bb^{(t)};\blam^{(t)})  
- \drho_-(\bb^{(t-1)}) + \bnu_{approx}^{(t)}, 
\eel
$t=1,\ldots, t^*$. 
For example, we may change the penalty levels proportionally by taking 
$\blam^{(t)} = A^{(t)}\blam_*$ with sufficiently large $A^{(0)}$ to ensure $\bb^{(0)}=0$ 
and decreasing $A^{(t)}$. 

Alternatively, we may move gradually from the Lasso to $\rho_\#(\bb;\blam_*)$, 
\bel{LCA-2}
&& \lam^{(t)}_j = \max\big(\lam_{*,j}, \theta^t\lam_{*,1}\big),\ j\le p, 
\\ \nonumber 
&& %\rho_\#^{(t)}(\bb;\blam^{(t)}) 
\Pen^{(t)}(\bb) = \theta^t \lam_{*,1}\|\bb\|_1 + (1-\theta^t) \rho_\#(\bb;\blam^{(t)}), 
\eel
with suitable $\theta <1$. 
The prediction and squared $\ell_2$ estimation error bounds are of the order 
$(\sigma^2/n)s\log p$ for the Lasso and $(\sigma^2/n)\{s  + s_1\log(p/s)\}$ for 
sorted concave penalties, where $s_1$ can be understood as the number of 
small nonzero coefficients. Thus, if we implement (\ref{LCA-2}) for a total of $t^*$ steps, 
we need $\theta_0^{t^*} \le \{1 + (s_1/s)\log(p/s)\}/\log p$ to achieve the better rate. 
This and parallel calculation for the weight for the Lasso and 
penalty level in (\ref{LCA-2}) lead to the following requirement on $t^*$: 
\bes
\max\big(\theta_0^{t^*},\theta^{2t^*}\big) \le \{1 + (s_1/s)\log(p/s)\}/\log p,
\ees
as $\lam_{*,1}^2/\lam_{*,s+1}^2 \approx \log(p/s)/\log p$ is of no smaller order. 

%Moreover, as each step of LCA (\ref{eq-4-29}) is a convex minimization problem, 
%it is computationally manageable to obtain an local solution with appropriate 
%approximation error $\bnu^{(t)}_{approx}$ satisfying (\ref{eq-4-30}). 
%More details about the computation of (\ref{eq-4-29}) will be studied in forthcoming papers.
%
%As to our knowledge, Theorem \ref{th-5} is the first results to guarantee the algorithmic 
%performance of approximate solutions with non-separable penalties. 
%It applies to general LCA with truly sorted penalties.
% On the contrary, all the existing studies of approximate solutions focus on separable 
% penalties with sufficient large penalty level. 
% For instance, the LLA is not feasible for truly sorted concave penalties. 
%
%\chz{Long: Please add some comments, e.g. one short paragraph each for 
%(a) the number of steps from zero 
%to sorted MCP, or from the Slope to sorted MCP, (b) each LCA step is computationally 
%manageable as convex minimization (to be studies in forthcoming paper), 
%(c) comparison with existing results on LLA for separable penalties} 
%

\subsection{Probabilistic error bounds}\label{subsec-4-6} 
Here we prove conditions of Theorems \ref{th-2c} and \ref{th-5} holds with high 
probability under the normality assumption (\ref{normality}) and the 
minimum penalty level condition (\ref{eq-4-18}) and (\ref{sorted-lam}) 
respectively for fixed and sorted penalty levels. %, and suitable RE and concavity condition. 

In addition to (\ref{normality}) we assume (\ref{eq-4-10}) holds with $0<\eta_1<\eta$ and $r_2\ge 0$ 
and (\ref{eq-4-18}) holds for fixed penalty level. 
For sorted penalty levels we assume (\ref{sorted-lam}) holds with $\alpha\in (0,1/4)$ and 
$A_0 = A/(\eta - \eta_1)$ for some $A>1$. 
For nonnegative integer $s$ and the above $\{\alpha,A\}$, define 
$p_{\alpha,A} = 2\alpha \sum_{k=0}^\infty \alpha^{(A-1)A^k}$, 
$q_{\alpha,A}=(1-\sqrt{2p_{\alpha,A}})_+$, $x_1=s/q_{\alpha,A}$, 
$L_x = \sqrt{2\log(p/(\alpha x))}$ and 
\bel{def-mu}
\mu_{\#,s} = 
\bigg\{\frac{2s (x_1/p)^{A^2-1}(2/q_{\alpha,A})^{2}}{A^2L_{s+1}^2(A^2L_{x_1}^2+2)}\bigg\}^{1/2}I_{\{s = 0\}}. 
\eel
We assume here that $x_1\le p$. This is reasonable when $p/s$ is large.
%\bes
%\frac{4 p(x_1/p)^{A^2}}{A^2L_{s+1}^2(A^2L_{x_1}^2+2)} = \mu_{\#,s}^2 q_{\alpha,A}/2,\ 
%\mu_{\#,s}^2 = \frac{4 (s/q_{\alpha,A})(x_1/p)^{A^2-1}(2/q_{\alpha,A})}{A^2L_{s+1}^2(A^2L_{x_1}^2+2)} 
%\ees

\begin{theorem}\label{th-6} 
Suppose (\ref{max-penalty}), (\ref{normality}) and (\ref{eq-4-10}) with $\eta_1\in (0,\eta)$ and $r_2>0$. 
Let $\bbeta^o=\hbbeta^o$ be the oracle LSE as in (\ref{oracle-LSE}) with $\calS = \supp(\bbeta^*)$ 
and $s=|\calS|$. Let $\scrB_0^*(\lambda_*,\kappa_*,\gamma)$ be the solution set in (\ref{eq-4-21}), 
$\xi = \{(1-1/a)(1-\eta)\}^{-1}\big[\big\{(\eta-\eta_1)^2\mu_{\#,s}^2+\eta_*^2\big\}^{1/2}+r_2/s^{1/2}\big]$ 
with the $\mu_{\#,s}$ in (\ref{def-mu}) and $a\in (0,1)$. 
Let $F(\bu) = \max\big\{\|\bu\|_2, (\gamma \bu^T\bSigmabar\bu)^{1/2}\big\}$.   
Suppose $\kappa_*\le (1-1/C^*)\RE_\#^2(\calS;\xi,\gamma)$. 
Then, there exists an event $\Omega$ such that $\P\{\Omega\}\ge 1- e^{-\xi_*s\log(p/(\alpha(s+1)))}$ 
with $\xi_* = (1-\eta)^2\xi^2\gamma/\{a(\eta-\eta_1)\}^2$, and that in the event $\Omega$ 
\bel{th-6-1}\quad && 
\big\|\hbbeta - \hbbeta^o\big\| \le (1-\eta)C^*\lam\,
\sup_{\bu\neq 0}\frac{\|\bu\|\{\xi s^{1/2} F(\bu) - \|\bu_{\calS^c}\|_{\#,s}\}}
{\bu^T\bSigmabar\bu}
\eel
for all approximate solutions $\hbbeta\in \scrB_0^*(\lambda_*,\kappa_*,\gamma)$ 
and seminorms $\|\cdot\|$,  
where $\lam=\lam_{s+1}$ for sorted penalties and $\lam=\lam(\bbeta^o)$ as in (\ref{penalty-level}) otherwise. 
In addition, for all sorted or unsorted penalties satisfying $\|\dPen(\bbeta^o)/\lam\|_2 \le s_1$, 
\bel{th-6-2}
\big\|\hbbeta - \hbbeta^o\big\| \le C^*\lam\,\sup_{\bu\neq 0}
\frac{\|\bu\|\{r_1F(\bu) - (1-\eta)\|\bu_{\calS^c}\|_{\#,s}\}}
{\bu^T\bSigmabar\bu}
\eel
with $r_1 = (1-1/a)^{-1}\big[\big\{(\eta-\eta_1)^2\mu_{\#,s}^2+ s_1^2\big\}^{1/2} + r_2\big]$ 
and at least probability 
$\P\{\Omega\}- e^{- \xi_*'r_1\log(p/(\alpha(s+1)))}$ with $\xi_*'= \gamma/\{a(\eta-\eta_1)\}^2$. 
\end{theorem}

\begin{corollary}\label{cor-th-6} 
Suppose $\lam =\lam_{s+1} \asymp \sigma\sqrt{(2/n)\log (p/s)}$ in (\ref{th-6-2}) 
and $\kappabar(\hbbeta) \le (1-1/C_0)\RE_\#^2(\calS;\xi,\gamma)$ 
	with $C_0^2/\RE_2^2(\calS;\xi)=O_P(1)$. 
Then, 
\bes
\displaystyle
\|\bX\hbbeta-\bX\hbbeta^o\|_2^2/n+\|\hbbeta - \hbbeta^o\|_2^2 
+ \|\hbbeta - \hbbeta^o\|_{\#,s}^2/s=O_P(\sigma^2/n)r_1^2\log (p/s) 
\ees
with $r_1^2 = s_1 + s/\log^2(p/s) + r_2^2$, where $s_1 = \|\dPen(\bbeta^o)/\lam\|_2^2$, and 
	\bel{cor-th-2-1}
	&& 
	\|\bX\hbbeta-\bX\bbeta^*\|_2^2/n+\|\hbbeta - \bbeta^*\|_2^2 + \|\hbbeta - \hbbeta^*\|_{\#,s}^2/s
	\\ \nonumber &=& O_P(\sigma^2/n)
	\big\{(s_1+r_2^2)\log (p/s) + s\big\}. 
	\eel
For sorted concave penalty (\ref{concave-slope}) with $\supp(\drho(\cdot;\lam))\subseteq[-\gamma\lam,\gamma\lam]$,
\bes
s_1 \le \#\big\{j\le |\calS|: (\bbeta^o)^\#_j \le \gamma\lam_j\big\}. 
\ees
\end{corollary}

Corollary \ref{cor-th-6} extends Corollary \ref{cor-th-1} to smaller penalty levels 
$\lam = \lam_{s+1} \ge A_0\sigma\sqrt{(2/n)\log(p/\alpha(s+1))}$, sorted penalties,  
one-step application of the LCA, and their approximate solutions. 
In the worst case scenario where $s_1\asymp s$, the error bounds in Corollary \ref{cor-th-6} 
attain the minimax rate \cite{YeZ10,BellecLT16}. 

Theorem \ref{th-6} and Corollary \ref{cor-th-6} provide sufficient conditions to guarantee 
simultaneous adaptation of sorted concave PLSE: (a) picking level $\lam_{s+1}$ automatically 
from $\{\lam_1,\ldots,\lam_p\}$, and (b) partially removing the bias of 
the Slope \cite{SuC16,BellecLT16} when $s_1 \ll s$, 
without requiring the knowledge of $s$ or $s_1$. 

Theorem \ref{th-6} is a direct consequence of Theorem \ref{th-2c} and the following proposition. 

\begin{proposition}\label{prop-slope-a} 
Suppose the normality assumption (\ref{normality}) holds. 
Let $\eta >\eta_1$ and $\{\lam_j\}$ be as in (\ref{sorted-lam}) for all sorted penalties. 
Let $s=|\calS|$ and $\{q_{\alpha,A},\mu_{\#,s}\}$ be as in (\ref{def-mu}) with $q_{\alpha,A}>0$. 
Let $a>0$, $w>0$ and $L=\sqrt{2\log(p/(\alpha(s+1)))}$. 
Suppose $\supp(\bz)\subseteq\calS^c$. Then, 
\bel{prop-5-1}
&& \P\Big\{ \hbox{(\ref{eq-4-10}) and $\|\bw_{\calS}\|_2\le w$ imply (\ref{eq-4-26})}\Big\} 
\ge \Phi\bigg(\frac{r_1L\sqrt{\gamma}}{a(\eta-\eta_1)}\bigg)
\eel
when $(1-1/a)r_1\ge \big\{(\eta-\eta_1)^2\mu_{\#,s}^2+ w^2\big\}^{1/2} + r_2$. 
In particular, when (\ref{approx-sol}) is also confined to penalties satisfying (\ref{max-penalty}), 
\bel{prop-5-1b}\quad
&& \P\Big\{ \hbox{(\ref{eq-4-10}) and $\|\bw_{\calS}\|_2\le \eta_*s^{1/2}$ imply (\ref{eq-4-7})}\Big\} 
\ge \Phi\bigg(\frac{(1-\eta)\xi s^{1/2}L}{a(\eta-\eta_1)\gamma^{-1/2}}\bigg),
\eel
with $\xi \ge \{(1-1/a)(1-\eta)\}^{-1} 
\big[\big\{(\eta-\eta_1)^2\mu_{\#,s}^2/s+\eta_*^2\big\}^{1/2}+r_2/s^{1/2}\big]$, 
where $\mu_{\#,s}$ is as in (\ref{def-mu}).  
Moreover, when (\ref{approx-sol}) is confined to penalties with the minimum fixed penalty 
level $\lam(\bbeta^o)=\lam\ge\lam_*$ as in (\ref{penalty-level}) and (\ref{eq-4-18}), 
(\ref{prop-5-1}) and (\ref{prop-5-1b}) hold with the $L$ in (\ref{eq-4-18}) 
and $\mu_{\#,s} = \sqrt{4s/(L^4+2L^2)}$. 
%(ii) When (\ref{approx-sol}) is confined to penalties satisfying (\ref{max-penalty}), (\ref{prop-5-1b}) holds 
%for $\xi$ satisfying both $\xi \ge \big\{2(\eta-\eta_1)/(L^4+2L^2)^{1/2}\big\}/\{(1- \nu-1/a)(1-\eta)\}$ and 
%$\xi \ge \{(\eta_*+\eta')+  r_2/s^{1/2}\}/\{\nu(1-\eta)\}$ with $s=|\calS|$. \\
%%$r_1 \nu \ge \{(\eta_*+\eta')s^{1/2}+  r_2\}$. 
\end{proposition}

\bibliographystyle{apalike}
\bibliography{Concave-PLSE-8.bib}

%Notes: Let $g_j^\#$ be the ordered $|g_1|,\ldots,|g_p|$ with $g_j\sim N(0,\sigma_j)$ and $\sigma_j\le 1$. 
%	Since $\exp\big(s^{-1}\sum_{j=1}^s x_j\big) \le s^{-1}\sum_{j=1}^s e^{x_j}$, 
%	\bes
%	\E \exp\bigg(\frac{t}{s}\sum_{j=1}^s g_j^\#\bigg)
%	\le \frac{1}{s}\sum_{i=1}^p \E e^{t|g_j|} 
%	= (2p/s)e^{t^2/2}\Phi(t),\ \forall t>0,
%	\ees
%	as $2\int_0^\infty e^{zt}\varphi(z)dz \le 2e^{t^2/2}\int_0^\infty \varphi(z-t)dz=2e^{t^2/2}\Phi(t)$. 
%	It follows that 
%	\bes
%	\P\Big\{ g_s^\# > t\Big\} \le e^{-t^2}
%	\E \exp\bigg(\frac{t}{s}\sum_{j=1}^s g_j^\#\bigg)
%	\le (2p/s)e^{- t^2/2},
%	\ees

%\section{Discussion}\label{sec-discussion}
%\chz{Figure 1 depicts branches of the solution space $\scrB$ in a 2-dimensional example.}
%\chz{Polytope? 
%\bes
%\begin{cases}
%\sgn(\hbeta_j)\bX_j^T(\by -\bX\hbbeta)/n \ge (\lam - \kappabar|\hbeta_j|)_+ \cr
%\|\bX^T(\by -\bX\hbbeta)/n - \bnu_{approx}\|_\infty \le \lam
%\end{cases}
%\ees
%While (\ref{eq-2-10-approx}) naturally define a class of approximate 
%solutions for (\ref{eq-2-10}), the above can be seen as an extension of all 
%}
\newpage
%\begin{supplement}
%%\sname{Supplementary Material to}	
%%\stitle{``Sorted Concave Penalized Regression"}
%	%\slink[url]{Supplementary Material to Sorted Concave Penalized Regression}
%%	\sdescription{We provide all the proofs in the supplementary material. The proofs are given in the following order: 
%%		Proposition \ref{prop-mixed-pen}, Proposition \ref{prop-slope}, Proposition \ref{prop-prox}, Lemma  \ref{lm-1}, Proposition \ref{prop-2}, Theorem \ref{th-1a}--\ref{th-1c}, Proposition  \ref{prop-2a}, Theorem  \ref{th-2c}, Theorem \ref{th-5} and Proposition \ref{prop-slope-a} along with an additional Lemma.}
%	\end{supplement}
%\medskip

%\begin{appendices}
	\appendix
	\section{Appendices}
	\addcontentsline{toc}{section}{Appendices}
	\setcounter{equation}{0}
	\renewcommand{\thesubsection}{A}
%\centerline{\sc\large\bf Appendix}
%\centerline{\sc\large ``Sorted Concave Penalized Regression"}
%\medskip

We provide proofs in the following order: 
		Proposition \ref{prop-mixed-pen}, Proposition~\ref{prop-slope}, Proposition \ref{prop-prox}, Lemma  \ref{lm-1}, Proposition \ref{prop-2}, Theorem \ref{th-1a}--\ref{th-1c}, Proposition  \ref{prop-2a}, Theorem  \ref{th-2c}, Theorem \ref{th-5}, and Proposition \ref{prop-slope-a} along with an additional lemma.
We omit proofs of Proposition \ref{prop-1}, Lemma \ref{lm-LCA} and Theorem \ref{th-6} as explained above or below their statements. 

\medskip
{\sc Proof of Proposition \ref{prop-mixed-pen}.} 
	 For $t_0=0<t_1<\ldots<t_K=1$, (\ref{mixed-sub-diff}) yields  
	 \bes
	 && \bh^T\big\{\E_{\nu}\big[\drho(\bb;\blam)\big|\bb\big] - \E_{\nu}\big[ \drho(\bb+\bh;\blam)\big|\bb+\bh\big]\big\}
	 \cr &=& \hbox{$\sum_{k=1}^K$} \bh^T\big\{\E_{\nu}\big[\drho(\bb_{t_{k-1}};\blam)\big|\bb_{t_{k-1}}\big] 
	 - \E_{\nu}\big[ \drho(\bb_{t_k};\blam)\big|\bb_{t_k}\big]\big\}
	 \cr &\le& \kappabar(\rho)\|\bh\|_2^2+\hbox{$\sum_{k=1}^K$} 
	 \bh^T\big\{\E_{\nu}\big[\drho(\bb_{t_{k-1}};\blam)\big|\bb_{t_{k-1}}\big] 
	 - \E_{\nu}\big[ \drho(\bb_{t_{k-1}};\blam)\big|\bb_{t_k}\big]\big\}
	 \cr &\to& \kappabar(\rho)\|\bh\|_2^2+r_n\hbox{$\int_0^1$} 
	 \Var_{\nu}\big(\bh^T\drho(\bb_t;\blam)\big|\bb_t\big) dt
	 \ees
	 with $\bb_t=\bb+t\bh$, where the limit is taken as $\max_{k\le K}|t_k-t_{k-1}|\to 0$. 
	 Note that $\bh^T\big\{\drho(\bb_{t_{k-1}};\blam) - \drho(\bb_{t_k};\blam)\big\}\le (t_k-t_{k-1})\kappabar(\rho)\|\bh\|_2^2$ 
	 by the definition of $\kappabar(\rho)$ in (\ref{eq-2-2}). 
	 This gives the upper bound for $\kappabar(\bb)$ and $\{\kappabar_1(\bb),\kappabar_2(\bb)\}$ 
	 	as $\bh^T\bM\bh\le \|\bh\|_1^2\|\bM\|_{\max}$ for all matrices $\bM\in \R^{p\times p}$. 
	 	We note that 
	 	$\Cov_{\nu}\big(\drho(\bu;\blam),\drho(\bu;\blam)\big|\bu,\theta\big)$ 
	 	is a diagonal matrix when the components of $\blam$ are independent given $\theta$.
	 
	 Let $\bdelta\in\{-1,1\}^p$ and take $\drho(b_j;\lam_j) = \delta_j\lam_j$ for $j\in \calS_{\bb}^c$. 
	 When $\blam$ is an exchangeable random vector under $\nu(d\blam)$,  
	 \bes
	 \Big(\E_\nu\big[ \drho(\bb;\blam)\big|\bb \big]\Big)_{\calS_{\bb}^c} 
	 = \bdelta_{\calS_{\bb}^c}\E\big[\lam_j \big|\bb\big],\ \forall\ j\in\calS_{\bb}^c. 
	 \ees
	 As $\pa\Pen(\bb)$ is convex, this gives $\lam(\bb) = \E[\lam_j|\bb]$ for $\rho_{\nu}(\bb),\ j\in\calS_{\bb}^c$. 
	 $\hfill\square$ 
	 
	 \medskip
	 {\sc Proof of Proposition \ref{prop-slope}.} 
	 We omit the proof of (\ref{prop-slope-1}) as it follows directly from the definition 
	 in (\ref{sub-diff}). 
	 For $\lam \ge \lam'$ and $t \ge t' \ge 0$, 
	 \bes
	 \rho(t;\lam)+\rho(t';\lam') - \rho(t;\lam') - \rho(t';\lam)
	 = \int_{t'}^t \big\{\drho(x;\lam) - \drho(x;\lam')\big\}dx \ge 0. 
	 \ees
	 Thus, (\ref{slope-2}) follows. By the convexity of $\rho(t;\lam_j)+\kappabar(\rho) t^2/2$ in $t$, 
	 \bes
	 \rho_\#(\bb;\blam) + \kappabar(\rho)\|\bb\|_2^2/2
	 = \max\bigg[\sum_{j=1}^p\big\{\rho(b_j;\lam_{k_j})+\kappabar(\rho)|b_j|^2/2\big\}: 
	 \bk \in \hbox{\rm perm}(p)\bigg]
	 \ees
	 is convex, as the maximum of convex functions is convex. 
	 Hence, the maximum concavity of $\rho_\#(\bb;\blam)$ is no greater than $\kappabar(\rho)$. 
	 $\hfill\square$
	 
	 	\medskip
	 	{\sc Proof of Proposition \ref{prop-prox}}. Because $\rho_\#(\bb;\blam)$ depends only 
	 	on $\bb^\#$, we have $\sgn(b_j)=\sgn(x_j)$, and the problem is to minimize 
	 	\bes
	 	\sum_{\ell=1}^p\frac{(|x_\ell|-|b_\ell|)^2}{2}+\rho_\#(\bb;\blam) 
	 	= \sum_{j=1}^p\frac{(|x_{\ell_j}|-b^\#_j)^2}{2}+\rho_\#(\bb;\blam). 
	 	\ees
	 	By algebra, $(|x_{\ell_j}|-b^\#_j)^2 +(|x_{\ell_k}|-b^\#_k)^2  - (|x_{\ell_k}|-b^\#_j)^2 - (|x_{\ell_j}|-b^\#_k)^2
	 	= (|x_{\ell_j}| - |x_{\ell_k}|)(b^\#_k-b^\#_j)$. 
	 	Thus, when $(|x_{\ell_j}| - |x_{\ell_k}|)(b^\#_k-b^\#_j)>0$, the solution can be improved by 
	 	switching $x_{\ell_j}$ and $x_{\ell_k}$. 
	 	This implies the monotonicity of $|b_j|$ in $|x_j|$ and (\ref{prop-prox-1}). 
	 	
	 	Now consider $x_1\ge\ldots\ge x_p\ge 0$. Let $\btil_j=\argmin\{(x_j-b)^2/2+\rho(b;\lam_j)\}$ 
	 	be the univariate solution and $\bb = \hbox{\rm iso.prox}(\bx;\rho(\cdot;\blam))$ the multivariate 
	 	isotonic solution. 
	 	If $b_j > b_{j+1}$ and $b_j>\btil_j$, then we can improve $\bb$ by decreasing $b_j$. 
	 	If $b_j > b_{j+1}$ and $b_{j+1}<\btil_{j+1}$, then we can improve $\bb$ by increasing $b_{j+1}$. 
	 	Thus, we must have $\btil_j > \btil_{j+1}$. Moreover, when 
	 	\bes
	 	\btil_{j'-1}>\btil_{j'}\le \btil_{j'+1}\le \cdots\le \btil_{j''}>\btil_{j''+1},\quad \btil_{j'}<\btil_{j''},
	 	\ees
	 	we must have $b_{j'} = b_{j'+1} = \cdots = b_{j''}$. This groups the optimization problem 
	 	in the first round of Algorithm 3. The argument also applies to the grouped optimization problem in 
	 	the second round, so on and so forth. $\hfill\square$ 
	 
	 %solutions of \chzR{(\ref{eq-2-10-approx})}
	 %\bel{eq-2-13}
	 %\bw = \big\{\dPen(\bbeta^o) \chz{+ \bnu_{approx}}  - \bX^T(\by-\bX\bbeta^o)/n\big\}\big/\lam(\bbeta^o). 
	 %\eel
	 \medskip
	 {\sc Proof of Lemma \ref{lm-1}.} Let $\bz = \bX^T(\by - \bX\bbeta^o)/n$.  
	 Recall that $\bh = \hbbeta - \bbeta^o$ and 
	 $\bw = \big\{\dPen(\bbeta^o)+ \bnu_{approx} - \bz\big\}\big/\lam$. 
	By (\ref{eq-2-10-approx}), 
	 %As (\ref{eq-2-10}) holds for $\hbbeta$, 
	 $\dPen(\hbbeta)+ \bnu_{approx} = \bz - \bSigmabar\bh$, so that 
	 \bes
	 \bh^T\bSigmabar\bh 
	 &=& \bh^T\big\{-\lam\bw + \dPen(\bbeta^o)- \dPen(\hbbeta)\big\} 
	 \cr &\le& -\lam \bh^T\bw + \kappabar_1(\hbbeta)\|\bh\|_1^2+\kappabar_2(\hbbeta)\|\bh\|_2^2
	 \ees
	 by the definition of $\kappabar_1(\bb)$ and $\kappabar_2(\bb)$ in (\ref{concavity-1-2}). 
	 As $\bh_{\calS^c}^T(\bnu_{approx})_{\calS^c}\ge 0$,
	 this gives (\ref{lm-1-1}) with $\dPen(\bbeta^o)\in \pa\dPen(\bbeta^o)$ 
	 satisfying $\dPen_j(\bbeta^o) = \sgn(h_j)\lam, j\in \calS^c$. 
	 %$-\lam \bh^T\bw\le \sum_{j\in\calS^c}(z_jh_j -\lam|h_j|) - \lam \bw_{\calS}^T\bh_{\calS}$. 
	 $\hfill\square$
	 
	 \medskip
	 {\sc Proof of Proposition \ref{prop-2}.} 
	 %Let ${a}_0=\min\big[{a}_2, {a}_2{a}_3/\{(1-\eta)(\xi\vee 1)\}\big]$, 
	 %$\lam$ be the penalty level $\lam(\bbeta^o)$ associated with the penalty for $\hbbeta$ 
	 %and $\tlam$ with $\tbbeta$. 
	 Let $\bh = \hbbeta-\bbeta^o$ and $\tbh = \tbbeta - \bbeta^o$. 
	 We want to prove that %when $\tbh_{\calS^c}\ge\xi\|\tbh_{\calS}\|_1$, 
	 \bel{pf-th-1-0}\qquad
	 \|\bh - \tbh\|_1\le {a}_0\tlam \hbox{ and } \bh \in\scrC(\calS;\xi) 
	 \hbox{ imply } \tbh \in \scrC(\calS;\xi). 
	 \eel
	 As $\|\bz_{\calS^c}\|_\infty \le (\eta- {a}_1)\lam_* \le (\eta- {a}_1)\tlam$, 
	 Lemma \ref{lm-1} and (\ref{cond-pen}) imply that 
	 \bel{pf-th-1-1a}
	 && \tbh^T\bSigmabar\tbh  + (1  - \eta + {a}_1)\tlam\|\tbh_{\calS^c}\|_1\
	 \\ \nonumber &\le& (1-\eta)\xi \tlam \|\tbh_{\calS}\|_1
	 + \kappabar_1(\tbbeta)\|\bh\|_1^2+\kappabar_2(\tbbeta)\|\bh\|_2^2. 
	 \eel
	 %with either $\bb=\bbeta^o$ or $\bb = \tbbeta$. 
	 Recall ${a}_2 = {a}_1\xi/[2(\xi+1)\big\{\kappabar_1(\bb)+ \kappabar_2(\bb)\big\}]$. 
	 When $\|\bh\|_1\vee\|\tbh-\bh\|_1\le{a}_2\tlam$, 
	 \bes
	 \kappabar_1(\tbbeta)\|\tbh\|_1^2+\kappabar_2(\tbbeta)\|\tbh\|_2^2
	 \le 2{a}_2\tlam\|\tbh\|_1\big\{\kappabar_1(\tbbeta)+ \kappabar_2(\tbbeta)\big\}
	 \le {a}_1\xi\tlam\|\tbh\|_1/(\xi+1),
	 \ees 
	 so that (\ref{pf-th-1-1a}) implies 
	 \bes\quad&
	 \big(1 - \eta + {a}_1 \big)\|\tbh_{\calS^c}\|_1
	 \le (1-\eta)\xi \|\tbh_{\calS}\|_1+ \{{a}_1\xi/(\xi+1)\}\big(\|\tbh_{\calS}\|_1+\|\tbh_{\calS^c}\|_1\big), 
	 \ees
	 %$(1-\eta+{a}_1-x)\xi = (1-\eta)\xi + x$, $({a}_1-x)\xi = x$, $x = {a}_1\xi/(\xi+1)$. 
	 which is equivalent to $\|\tbh_{\calS^c}\|_1\le \xi\|\tbh_{\calS}\|_1$ by algebra.
	 
	 Because $\kappabar_{1,2}(\hbbeta;\xi)\le \kappa_*$ and 
	 $\|\bh\|_1 \le (1+\xi)\|\bh_{\calS}\|_1\le (1+\xi)|\calS|^{1/2}\|\bh\|_2$,
	 %by $\|\bh_{\calS^c}\|_1\le \xi\|\bh_{\calS}\|_1$},  
	 the $\{\bh,\lam\}$ version of (\ref{pf-th-1-1a}) implies 
	 \bes\qquad
	 & \bh^T\bSigmabar\bh  + (1  - \eta + {a}_1)\lam\|\bh_{\calS^c}\|_1\
	 \le (1-\eta)\xi \lam \|\bh_{\calS}\|_1+ \kappa_*\|\bh\|_2^2. 
	 \ees
	 By the RE condition, we have $\kappa_* \|\bh\|_2^2\le \bh^T\bSigmabar\bh$, so that 
	 \bel{pf-th-1-1b}
	 \big(1- \eta+{a}_1\big)\|\bh_{\calS^c}\|_1 \le (1-\eta)\xi\|\bh_{\calS}\|_1. 
	 \eel
	 Recall that ${a}_3 = {a}_1(1-\eta)\xi/\{(1-\eta)(\xi + 1)+{a}_1\}$.  
	 If $\|\bh\|_1 > {a}_2\tlam$ and $\|\tbh-\bh\|_1\le \tlam{a}_2{a}_3/\{(1-\eta)(\xi\vee 1)\}$, 
	 (\ref{pf-th-1-1b}) implies 
	 \bes
	 && \big(1- \eta\big)\|\tbh_{\calS^c}\|_1 - (1-\eta)\xi\|\tbh_{\calS}\|_1
	 \cr &\le& \big(1- \eta\big)\|\bh_{\calS^c}\|_1 - (1-\eta)\xi\|\bh_{\calS}\|_1 
	 + \{(1-\eta)(\xi\vee 1)\}\|\tbh-\bh\|_1
	 \cr &\le& \big(1- \eta\big)\|\bh_{\calS^c}\|_1 - (1-\eta)\xi\|\bh_{\calS}\|_1 
	 + {a}_3\|\bh\|_1
	 \cr &=& \big(1- \eta + {a}_3\big)\|\bh_{\calS^c}\|_1 - \{(1-\eta)\xi - {a}_3\}\|\bh_{\calS}\|_1 
	 \cr &=& \frac{\big(1- \eta+{a}_1\big)\|\bh_{\calS^c}\|_1 - (1-\eta)\xi\|\bh_{\calS}\|_1}
	 {\{(1-\eta)(\xi + 1)+{a}_1\}/\{(1-\eta)(\xi+1)\}}
	 \le 0
	 \ees
	 by algebra.  
	 %
	 %$\{(1-\eta)(\xi + 1)+{a}_1\}(1- \eta) + {a}_1(1-\eta)\xi
	 %= (1-\eta)(\xi + 1)(1-\eta+{a}_1)$
	 %
	 %$\{(1-\eta)(\xi + 1)+{a}_1\}(1-\eta)\xi - {a}_1(1-\eta)\xi
	 %= (1-\eta)(\xi + 1)(1-\eta)\xi$
	 %
	 Thus, (\ref{pf-th-1-0}) holds in either cases. $\hfill\square$
	 
	 \medskip
	 {\sc Proof of Theorem \ref{th-1a}--\ref{th-1c}.}
	 Let $\bh = \hbbeta-\bbeta^o$. 
	 By Proposition \ref{prop-2}, $\bh \in \scrC(\calS;\xi)$.
	 It follows that $\|\bh\|_1^2\le (1+\xi)^2|\calS|\|\bh_{\calS}\|_2^2$, so that by Lemma \ref{lm-1} and (\ref{eq-2-12})
	 \bes
	 \bh^T\bSigmabar\bh+(1-\eta)\lambda\|\bh_{\calS^c}\|_1
	 + \lambda\bw_{\calS}^T\bh_{\calS}\le \kappabar_{1,2}(\hbbeta;\xi)\|\bh\|_2^2 
	 \le \kappa_*\|\bh\|_2^2.
	 \ees
	 %as $\hbbeta$ is connected to ${\bf 0}$ through a continuous path and the origin has the cone property. 
	 Let $C_0$ be as in Theorem \ref{th-1c}. 
	 As $\bh^T\bSigmabar\bh\ge \RE^2(\calS;\xi)\|\bh\|_2^2$, 
	 Lemma \ref{lm-1} implies 
	 \bel{pf-th-1-2}\qquad 
	 & C_0^{-1}\bh^T\bSigmabar\bh +(1-{\eta})\lambda\|\bh_{\calS^c}\|_1 
	 \le - \lambda \bw_S^T\bh_{\calS}
	 \le \lambda\|\bw_S\|_2\|\bh\|_2. 
	 \eel
	 This immediately implies (\ref{th-1-2a}) and (\ref{th-1-2b}). 
	 For (\ref{th-1-1}) and (\ref{th-1-2}), we set $C_0=\infty$. 
	 However, by (\ref{pf-th-1-2}) and the definition of $\RCIF$, 
	 \bes
	 \RCIF_{\rm pred}(\calS;{\eta},\bw) \bh^T\bSigmabar\bh 
	 \le \|\bSigmabar\bh\|_\infty^2|\calS|.
	 \ees
	 Consequently, the first inequality in (\ref{th-1-1}) follows from the fact that
	 \bes
	 \|\bSigmabar\bh\|_\infty 
	 &\le& \|\bX^T(\by-\bX\hbbeta)/n\|_{\infty}
	 +\|\bX^T(\by-\bX\bbeta^o)/n\|_{\infty}
	 \cr &\le& \lam+\|\bnu_{approx}\|_\infty + \eta\lam_*=\lambar,
	 \ees
	 and the second inequality follows from the first inequality in Proposition~\ref{prop-1}.
	 Similarly, the first inequality in (\ref{th-1-2}) follows from 
	 \bes
	 \RCIF_{\rm est, q}(\calS;{\eta},\bw) \|\bh\|_q
	 \le \|\bSigmabar\bh\|_\infty|\calS|^{1/q}\le \lambar|S|^{1/q}, 
	 \ees 
	 and the second from the second and third inequalities in Proposition~\ref{prop-1}.
	 
	 Finally we consider selection consistency 
	 for exact local solutions $\hbbeta$ of (\ref{eq-2-10}) in $\scrB_0^*(\lam_*,\kappa_*)$
	 under the assumption $\bw_{\calS}={\bf 0}$. In this case, 
	$\bnu_{approx}={\bf 0}$,
	 $\bbeta^o$ is a solution of (\ref{eq-2-10}), and $\|\bh_{\calS^c}\|_1=0$ by (\ref{pf-th-1-2}). 
	 Moreover, because both $\bbeta^o$ and $\hbbeta$ are solutions of (\ref{eq-2-10}) 
	 with support in $\calS$, 
	 \bes
	 \bh_\calS^T\bSigmabar_\calS\bh_\calS 
	 = \bh^T\big\{\dPen(\bbeta^o)-\dPen(\hbbeta)\big\}
	 \le \kappabar_{1,2}(\hbbeta;\xi)\|\bh_{\calS}\|_2^2. 
	 \ees
	 Thus, $\phi_{\min}(\bSigma_{\calS,\calS}) > \kappabar_{1,2}(\hbbeta;\xi)$ implies $\bh={\bf 0}$. 
	 
	 For separable penalties of the form $\rho(\bbeta;\lam) = \sum_{j=1}^p \rho(b_j;\lam)$ in (\ref{separable-pen}), 
	 \bes\quad
	 \kappa_*\|\bh_{\calS}\|_2^2 
	 \le \bh^T\big\{\dPen(\bbeta^o)-\dPen(\hbbeta)\big\}
	 = \hbox{$\sum_{j\in\calS}$}(\hbeta_j-\beta_j^o)\big\{\drho(\beta^o_j;\lam)-\drho(\hbeta_j;\lam)\big\}. 
	 \ees
	 As $(\hbeta_j-\beta_j^o)\big\{\drho(\beta^o_j;\lam)-\drho(\hbeta_j;\lam)\big\}\le 
	 \min\{\kappabar(\beta_j^o;\rho,\lam),\kappabar(\hbeta_j;\rho,\lam)\}h_j^2\le \kappa_*h_j^2$, 
	 equality must attain for every $j$. 
	 This is possible only when $\sgn(\hbeta_j)\sgn(\beta^o_j)\ge 0$ for all $j\in \calS$. 
	 Furthermore, $\sgn(\hbeta_j)\sgn(\beta^o_j) > 0$ for all $j\in \calS$ when 
	 $\kappabar(0;\rho,\lam)<\phi_{\min}(\bSigma_{\calS,\calS})$. 
	 $\hfill\square$
	 
	 \medskip
	 {\sc Proof of Proposition \ref{prop-2a}.} 
	 Let $\hbbeta\in \scrB(\lambda_*,\kappa_*,\gamma)$ with penalty level 
	 $\lam=\lam_{s+1}$ and $\hbbeta-\bbeta^o\in \scrC_\#(\calS;\xi,\gamma)$, 
	 and $\tbbeta\in \scrB(\lambda_*,\kappa_*,\gamma)$ with $\tlam=\tlam_{s+1}$. 
	 We first prove that $\tbbeta\in \scrB_0(\lambda_*,\kappa_*,\gamma)$ implies 
	 $\tbbeta-\bbeta^o\in \scrC_\#(\calS;\xi,\gamma)$, or equivalently  
	 \bel{eq-4-24}
	 \hbox{$\|\tbbeta - \hbbeta\|_{\#,*} \le \tlam |\calS| a_0
	 	\ \Rightarrow\ \tbbeta-\bbeta^o\in \scrC_\#(\calS;\xi,\gamma)$.}
	 \eel
	 Let $F(\bu) = \max\big\{\|\bu\|_2, \|\bX\bu\|_2\sqrt{\gamma/n}\big\}$, 
	 $\bh = \hbbeta-\bbeta^o$ and $\tbh = \tbbeta - \bbeta^o$. 
	 Assume without loss of generality 
	 $\|\tbh_{\calS^c}\|_{\#,s} \ge \xi|S|^{1/2}\|\tbh\|_2$. 
	 By (\ref{eq-4-5}), (\ref{eq-4-6}) and (\ref{eq-4-23}) %(\ref{eq-5-prop-5-1}), 
	 \bel{eq-5-pf-prop-6-1a}
	 \quad && \tbh^T\bSigmabar\tbh + (1-\eta)\tlam\|\tbh_{\calS^c}\|_{\#,s} - \kappabar(\tbbeta)\|\tbh\|_2^2 
	 \le (1-\eta- {a}_1)\tlam\xi|\calS|^{1/2}F(\tbh). 
	 \eel
	 Recall that ${a}_2 = {a}_1\xi^2/\{2\kappabar(\tbbeta)\}$. 
	 When $\|\bh\|_{\#,*}\vee\|\tbh-\bh\|_{\#,*} \le{a}_2|\calS|\tlam$ and 
	 $\|\tbh_{\calS^c}\|_{\#,s} \ge \xi|S|^{1/2}\|\tbh\|_2$, 
	 $\kappabar(\tbbeta)\|\tbh\|_2^2\le \kappabar(\tbbeta)\|\tbh_{\calS^c}\|_{\#,s}^2/(\xi^2|S|)
	 \le {a}_1\|\tbh_{\calS^c}\|_{\#,s}$, so that 
	 $\|\tbh_{\calS^c}\|_{\#,s}\le \xi|\calS|^{1/2}F(\tbh)$ by (\ref{eq-5-pf-prop-6-1a}).
	 
	 Because $\kappabar(\hbbeta)\le \kappa_*$,   
	 the $\{\bh,\lam\}$ version of (\ref{eq-5-pf-prop-6-1a}) implies 
	 \bes\qquad
	 & \bh^T\bSigmabar\bh  + (1  - \eta)\lam\|\bh_{\calS^c}\|_{\#,s}\
	 \le (1-\eta - {a}_1)\xi \lam |\calS|^{1/2}F(\bh) + \kappa_*F^2(\bh). 
	 \ees
	 By the RE condition in $\scrC_\#(\calS;\xi,\gamma)$, 
	 we have $\kappa_* F^2(\bh) \le \bh^T\bSigmabar\bh$, so that 
	 \bes
	 \big(1- \eta\big)\|\bh_{\calS^c}\|_{\#,s} \le (1-\eta - {a}_1)\xi |\calS|^{1/2}F(\bh). 
	 \ees
	 Recall that ${a}_3 = {a}_1(1-\eta)\xi/\{(1-\eta - {a}_1)(\xi + 1)+{a}_1\}$.  
	 If $\|\bh\|_{\#,*} > {a}_2\tlam$ and $\|\tbh-\bh\|_{\#,*} \le \tlam{a}_2{a}_3/\{(1-\eta)(\xi\vee 1)\}$, we have 
	 \bes
	 && \big(1- \eta\big)\|\tbh_{\calS^c}\|_{\#,s} - (1-\eta)\xi|\calS|^{1/2}F(\tbh)
	 \cr &\le& \big(1- \eta\big)\|\bh_{\calS^c}\|_{\#,s} - (1-\eta)\xi|\calS|^{1/2}F(\bh) 
	 + \{(1-\eta)(\xi\vee 1)\}\|\tbh-\bh\|_{\#,*}
	 \cr &\le& \big(1- \eta\big)\|\bh_{\calS^c}\|_{\#,s} - (1-\eta)\xi|\calS|^{1/2}F(\bh) 
	 + {a}_3\|\bh\|_{\#,*}
	 \cr &=& \big(1- \eta + {a}_3\big)\|\bh_{\calS^c}\|_{\#,s} - (1-\eta - {a}_3/\xi)\xi|\calS|^{1/2}F(\bh) 
	 \cr &=& \frac{\big(1- \eta\big)\|\bh_{\calS^c}\|_{\#,s} - (1-\eta - {a}_1)\xi|\calS|^{1/2}F(\bh)}
	 {\{(1-\eta)(\xi + 1)+{a}_1\}/\{(1-\eta)(\xi+1)\}}
	 \le 0
	 \ees
	 due to $(1-\eta+a_3)/(1-\eta-a_3/\xi) = (1-\eta)/(1-\eta-a_1)$. 
	 Hence, (\ref{eq-4-24}) holds in either cases. 
	 The proof for the approximate LCA solution is straightforward as 
	 (\ref{eq-4-13}) corresponds to (\ref{eq-4-5}) with $\kappabar(\hbbeta)=0$. 
	 $\hfill\square$
	 
	 \medskip
	 {\sc Proof of Theorem \ref{th-2c}.} 
	 By (\ref{eq-4-5}), (\ref{eq-4-6}) and (\ref{eq-4-26}) that 
	 \bel{eq-4-28}
	 \quad && \bh^T\bSigmabar\bh + (1-\eta)\lam\|\bh_{\calS^c}\|_{\#,s} - \kappabar(\hbbeta)\|\bh\|_2^2 
	 \le r_1\lam F(\bh). 
	 \eel
	 with $F(\bu) = \max\big\{\|\bu\|_2, \|\bX\bu\|_2\sqrt{\gamma/n}\big\}$. 
	 By Proposition \ref{prop-2a}, $\bh\in\scrC_\#(\calS;\xi,\gamma)$, so that 
	 $\RE_\#^2(\calS;\xi,\gamma)F^2(\bh) \le \bh^T\bSigmabar\bh$. As 
	 $\kappabar(\hbbeta)\le (1-1/C_0)\RE_\#^2(\calS;\xi,\gamma)$, 
	 \bes
	 C_0^{-1}\bh^T\bSigmabar\bh + (1-\eta)\lam\|\bh_{\calS^c}\|_{\#,s}
	 \le r_1 \lam F(\bh),
	 \ees
	 which implies (\ref{eq-4-27}). $\hfill\square$ 
	 
	 \medskip
	 {\sc Proof of Theorem \ref{th-5}.} We begin induction by assuming 
	 \bes
	 %\|\bnu_{carry}^{(t-1)}\|_2/\lam^{(t)}\le 
	 \kappa_0\|\bh^{(t-1)}\|_2/\lam^{(t)}\le \nu_0/\lam^{(1)}, 
	 \ees
	 which holds for $t=1$. 
	 By (\ref{eq-4-6}), (\ref{eq-4-30}) and the condition on $\xi$, 
	 \bes
	 1 & \ge & \Delta\Big(r_1^{(t)}+\|\bnu_{carry}^{(t-1)}\|_2/\lam^{(t)},\bw^{(t)},  \bnu_{approx}^{(t)} - \bnu_{carry}^{(t-1)}\Big) 
	 \cr & \ge & \Delta\Big(r_1^{(t)}+\kappa_0\|\bh^{(t-1)}\|_2/\lam^{(t)},\bw^{(t)}, \bnu_{approx}^{(t)} - \bnu_{carry}^{(t-1)}\Big) 
	 \cr & > & \Delta\Big((1-\eta)\xi s^{1/2},\bw^{(t)},\bnu_{approx}^{(t)} - \bnu_{carry}^{(t-1)}\Big),
	 \ees
	 so that Theorem \ref{th-2c} applies. As LCA is convex minimization, $\kappabar(\hbbeta)=0$ 
	 in this application of Theorem \ref{th-2c}. Thus, 
	 \bes
	 F\big(\bh^{(t)}\big) %\max\big\{\|\bu\|_2, (\gamma \bu^T\bSigmabar\bu)^{1/2}\big\} 
	 \le \frac{r_1^{(t)}\lam^{(t)}+\kappa_0\|\bh^{(t-1)}\|_2}{\RE_\#^2(\calS;\xi,\gamma)}
	 \le \lam^{(t)}\frac{r_1^{(t)} + \nu_0/\lam^{(1)}}{\RE_\#^2(\calS;\xi,\gamma)}. 
	 \ees
	 As $\RE_\#^2(\calS;\xi,\gamma)\ge \kappa_0\{\lam^{(t)}/\lam^{(t+1)}\}\{r_1^{(t)}\lam^{(1)}/\nu_0+1\}$, 
	 \bes
	 \frac{\kappa_0F\big(\bh^{(t)}\big)}{\lam^{(t+1)}}
	 \le \kappa_0\frac{\lam^{(t)}}{\lam^{(t+1)}}\frac{r_1^{(t)} + \nu_0/\lam^{(1)}}{\RE_\#^2(\calS;\xi,\gamma)} \le \nu_0/\lam^{(1)}. 
	 \ees
	 This completes the induction as $\|\bnu_{carry}^{(t)}\|_2\le \kappa_0\|\bh^{(t)}\|_2\le \kappa_0F\big(\bh^{(t)}\big)$. 
	 $\hfill\square$ 
	 
	 \medskip
	 We need the following lemma in the proof of Proposition \ref{prop-slope-a}.
	 
	 \begin{lemma}\label{lm-slope-noise} 
	 	Let $z_j$ be normal variables with $\E\, z_j=0$ and $\Var(z_j)\le \sigma_n^2$. 
	 	Let $s, A, L_x, q_{\alpha,A}, x_1$ and $\mu_{\#,s}$ be as in (\ref{def-mu}). 
	 	Then, for $s>0$ and $q_{\alpha,A}>0$ 
	 	\bes&&
	 \P\bigg\{\max_{x_1 \le j \le p} \frac{z_{j-s}^\#}{A\sigma_n L_j} <1, 
	 		\sum_{ s < j < x_1}\frac{\big(z_{j-s}^\#-A\sigma_n L_j\big)_+^2}
	 		{(A\sigma_n L_{s+1})^2} < \mu_{\#,s}^2\bigg\} > 1/2, 
	 	\ees
	 	and the right-hand side is greater than 
	 	$(1+q_{\alpha,A})/2$ for $s=0$. 
	 \end{lemma}
	 
	 \medskip
	 {\sc Proof of Lemma \ref{lm-slope-noise}.} Set $\sigma_n=1$ without loss of generality. 
	 As in \cite{BellecLT16}, 
	 \bes
	 \P\{ z_j^\# > t\big\} \le (2p/j)e^{-t^2/2}. 
	 \ees
	 Recall that $p_{\alpha,A} = 2\alpha \sum_{k=0}^\infty \alpha^{(A-1)A^k}$, 
	 $q_{\alpha,A} =1-\sqrt{2p_{\alpha,A}}$ and $x_1=s/q_{\alpha,A}\le p$. 
	 Define $x_k$ by $L_{x_k}^2 = A^{1-k}L^2_{x_1}$,
	 $2\le k\le k^*$, with $A L^2_{p} > L^2_{j_{k^*}} \ge L^2_{p}$. 
	 Let $j_k$ be the solution of $j_k-1 < x_k \le j_k$. 
	 Let $y = 1-s/x_1=\sqrt{2p_{\alpha,A}}$. 
	 Because $z_j^\#\ge z_{j+1}^\#$ and $x_k/(x_k-s)\le 1/(1-s/x_1) = 1/y = y/(2p_{\alpha,A})$, 
	 \bel{pf-lm-slope-1-1}
	 && \P\Big\{\max_{x_1\le j \le p} z^\#_{j-s}/L_{j} \ge A\big\}
	 \cr &\le& \sum_{j_k < j_{k+1}}\P\big\{ z^\#_{j_k-s} \ge AL_{x_{k+1}} = A^{1/2}L_{x_k} \big\}
	 \cr &\le& \sum_{j_k < j_{k+1}} \frac{2p}{j_k-s}\exp\Big(-A L^2_{x_k}/2\Big)
	 \\ \nonumber &\le& \sum_{j_k < j_{k+1}} \frac{2p}{x_k-s}\bigg(\frac{\alpha x_k}{p}\bigg)
	 \exp\Big(-(A-1)A^{k^*-k}L^2_p/2\Big)
	 \\ \nonumber &\le& \frac{1}{y}\sum_{k\le k^*} 2\alpha^{1+(A-1)A^{k^*-k}} = p_{\alpha,A}/y = y/2. 
	 \eel
	 This completes the proof for $s=0$, which gives $x_1=0$. 
	 For $s>0$, Proposition 10 in \cite{SunZ13} gives 
	 \bes
	 \E \sum_{s < j<x_1}\frac{\big(z_{j-s}^\# - A L_{j}\big)_+^2}{A^2L_{s+1}^2}
	 \le \frac{\E\big\|(|\bz| - A L_{x_1})_+\big\|_2^2}{A^2L_{s+1}^2} 
	 \le \frac{4 p(x_1/p)^{A^2}}{A^2L_{s+1}^2(A^2L_{x_1}^2+2)}. 
	 \ees
	 As the right-hand side above equals to 
	 $\mu_{\#,s}^2 q_{\alpha,A}/2$, by Markov's inequality 
	 \bes
	 \P\bigg\{\sum_{s < j<x_1}\big(z_{j-s}^\# - A L_j\big)_+^2>A^2L_{s+1}^2\mu_{\#,s}^2\bigg\}
	 \le q_{\alpha,A}/2 = 1/2 - y/2.
	 \ees  
	 This and (\ref{pf-lm-slope-1-1}) completes the proof with the union bound. $\hfill\square$
	 
	 \medskip
	 {\sc Proof of Proposition \ref{prop-slope-a}.} 
	 Let $\eta'=\eta-\eta_1$, $\lam_{*,j}$ be as in (\ref{sorted-lam}) and  
	 \bes
	 \qquad && 
	 \Deltatil(r,\eta,w,r_2) =  \sup_{\|\bu\|_2=1} \,
	 \frac{\bu_{\calS^c}^T\bz_{\calS^c}/\lam_{*,s+1}
	 	- \eta \|\bu_{\calS^c}\|_{\#,s} + w\|\bu_{\calS}\|_2 + r_2}
	 {r\max\big\{\|\bu\|_2/a,\|\bX\bu\|_2\sqrt{\gamma/n}\big\}}
	 \ees
	 where $\|\cdot\|_{\#,s}$ is defined as in (\ref{eq-4-4}) with $\lam_j=\lam_{*,j}$. 
	 By Lemma \ref{lm-slope-noise}, 
	 \bes
	 \hbox{\rm median}\Big(\Deltatil(r,\eta,w,r_2)\Big) 
	 \le \big\{\big(\eta^2\mu_{\#,s}^2+w^2\big)^{1/2} + r_2\big\}(a/r). 
	 \ees
	 As $\supp(\bz)\subseteq\calS^c$, 
	 $t\sigma \Deltatil(r,\eta,w,r_2)$ is convex and unit Lipschitz in 
	 $\by-\bX\bbeta^o$, so that by the Gaussian concentration inequality 
	 \bes
	 \P\Big\{ \Deltatil(r,\eta,w,r_2) \ge \hbox{\rm median}\Big(\Deltatil(r,\eta,w,r_2)\Big) + 1\Big\} \le \Phi(-t). 
	 \ees
	 
	 We first prove (\ref{prop-5-1}) and (\ref{prop-5-1b}). Let $r=r_1$. 
	 Although $\Deltabar(r,\eta,w,r_2)$ in (\ref{eq-4-11}) and $\Deltatil(r,\eta,w,r_2)$ are defined with 
	 $\lam_{*,j}$ and $\Delta(r_1,\bw,\bnu)$ is defined in (\ref{eq-4-6}) with $\lam_j\ge \lam_{*,j}$, 
	 (\ref{eq-4-10}) and the monotonicity of $(z/\lam-1)_+$ in $\lam$ still provide 
	 %By (\ref{eq-4-6}), (\ref{eq-4-10}) and (\ref{eq-4-11}), 
	 \bes
	 \Delta(r_1,\bw,\bnu)I_{\{\|\bw_{\calS}\|_2\le w\}} \le \Deltabar(r_1,\eta',w,r_2)
	 \big/\hbox{\rm median}\big(\Deltatil(r,\eta',w,r_2)+1\big)
	 \ees
	 when 
	 \bes
	 \frac{1}{\hbox{\rm median}\big(\Deltatil(r,\eta',w,r_2)\big)+1} 
	 \ge \max\bigg(\frac{1}{a}, \frac{t \sigma}{r_1\lam\sqrt{\gamma n}}\bigg)
	 = \max\bigg(\frac{1}{a}, \frac{\eta't}{r_1L\sqrt{\gamma}}\bigg), 
	 \ees
	 Setting $t = r_1L\sqrt{\gamma}/(a\eta')$ 
	 and $\big\{\big((\eta')^2\mu_{\#,s}^2+w^2\big)^{1/2} + r_2\big\}/r_1+1/a = 1$, 
	 we have (\ref{prop-5-1}). Moreover, (\ref{prop-5-1b}) follows with $r_1=(1-\eta)\xi s^{1/2}$. 
	 For constant $\lam$ with $\|\bu_{\calS^c}\|_{\#,s}=\|\bu_{\calS^c}\|_1$, we replace the median with 
	 \bes
	 \E\, \Deltatil(r,\eta,w,r_2) \le \big\{\big(\eta^24s/(L^4+2L^2)+w^2\big)^{1/2} + r_2\big\}(a/r),
	 \ees
	 so that $\mu_{\#,s}^2$ can be replaced by $4s/(L^4+2L^2)$. $\hfill\square$
	 %For part (ii), we take $r=(1-\eta)\xi s^{1/2}$ and $T=\{1,\ldots,p\}$. 
	 %We have 
	 %\bes
	 %&& \bu_{\calS^c}^T\bz_{\calS^c}/\lam - \eta'\|\bu_{\calS^c}\|_1 - \bu_{\calS}^T\bw_{\calS} +  r_2\|\bu\|_2
	 %\cr &=& \bu^T\bz/\lam - \eta'\|\bu\|_1 - \bu_{\calS}^T\dPen_{\calS}(\bbeta^*)/\lam
	 %+\eta'\|\bu_{\calS}\|_1 +  r_2\|\bu\|_2 
	 %\cr &\le & \bu^T\bz/\lam - \eta'\|\bu\|_1 + \{(\eta_*+\eta')s^{1/2}+  r_2\}\|\bu\|_2 
	 %\ees
	 %due to $\bw = \{\dPen(\bbeta^o)-\bz\}/\lam$. 
	 %The argument in the proof of part (i) still applies with $w=0$ when 
	 %$r_1 \nu \ge \{(\eta_*+\eta')s^{1/2}+  r_2\}$. $\hfill\square$
%\end{appendices}
\end{document}

%% file: 04-pre-am-short.tex
\newcommand{\bel}{\begin{eqnarray}\label}
\newcommand{\eel}{\end{eqnarray}}
\newcommand{\bes}{\begin{eqnarray*}}
\newcommand{\ees}{\end{eqnarray*}}
\newcommand{\bei}{\begin{itemize}}
\newcommand{\beiftnt}{\begin{itemize}\footnotesize}
\newcommand{\eei}{\end{itemize}}

\def\benu{\begin{enumerate}}
\def\eenu{\end{enumerate}}

\def\argmin{\mathop{\rm arg\, min}}
\def\real{{\mathbb{R}}}
\def\R{{\real}}

\def\E{{\mathbb{E}}}
\def\P{{\mathbb{P}}}

\def\complex{\mathop{{\rm I}\kern-.58em\hbox{\rm C}}\nolimits}
\def\pa{\partial}

\def\sgn{\hbox{\rm sgn}}

\def\Cov{\hbox{\rm Cov}}

\def\Var{\hbox{\rm Var}}

\def\supp{\hbox{\rm supp}}

\def\mathbold{\boldsymbol} %\def\mathbold{\mathbf}
\def\bb{\mathbold{b}}\def\btil{\widetilde{b}}
\def\tbb{{\widetilde{\bb}}}
\def\scrB{{\mathscr B}}

\def\scrC{{\mathscr C}}

\def\scrC{{\mathscr C}}

\def\bg{\mathbold{g}}

\def\bh{\mathbold{h}}
\def\tbh{{\widetilde{\bh}}}

\def\bI{\mathbold{I}}

\def\bk{\mathbold{k}}

\def\bM{\mathbold{M}}

\def\calS{{\cal S}}

\def\bu{\mathbold{u}}

\def\bv{\mathbold{v}}

\def\bV{\mathbold{V}}

\def\bw{\mathbold{w}}

\def\bx{\mathbold{x}}

\def\bX{\mathbold{X}}

\def\by{\mathbold{y}}

\def\bz{\mathbold{z}}

\def\bbeta{\mathbold{\beta}}\def\hbeta{\widehat{\beta}}

\def\hbbeta{{\widehat{\bbeta}}}\def\tbbeta{{\widetilde{\bbeta}}}

\def\bdelta{\mathbold{\delta}}

\def\Deltatil{{\widetilde \Delta}}\def\Deltabar{{\overline \Delta}}

\def\ep{\varepsilon}
\def\bep{\mathbold{\ep}}

\def\kappabar{{\overline{\kappa}}}

\def\lam{\lambda}
\def\blam{\mathbold{\lam}}
\def\tlam{\widetilde{\lam}}\def\lambar{{\overline{\lam}}}

\def\bnu{\mathbold{\nu}}

\def\bSigma{\mathbold{\Sigma}}

\def\bSigmabar{{\overline\bSigma}}

%% file: ms.bbl
\begin{thebibliography}{}

\bibitem[Agarwal et~al., 2012]{agarwal2012fast}
Agarwal, A., Negahban, S., and Wainwright, M.~J. (2012).
\newblock Fast global convergence of gradient methods for high-dimensional
  statistical recovery.
\newblock {\em The Annals of Statistics}, pages 2452--2482.

\bibitem[Beck and Teboulle, 2009]{beck2009fast}
Beck, A. and Teboulle, M. (2009).
\newblock A fast iterative shrinkage-thresholding algorithm for linear inverse
  problems.
\newblock {\em SIAM journal on imaging sciences}, 2(1):183--202.

\bibitem[Bellec et~al., 2016]{BellecLT16}
Bellec, P.~C., Lecu{\'e}, G., and Tsybakov, A.~B. (2016).
\newblock Slope meets lasso: improved oracle bounds and optimality.
\newblock {\em arXiv preprint arXiv:1605.08651}.

\bibitem[Bickel et~al., 2009]{BickelRT09}
Bickel, P.~J., Ritov, Y., and Tsybakov, A.~B. (2009).
\newblock Simultaneous analysis of lasso and dantzig selector.
\newblock {\em The Annals of Statistics}, pages 1705--1732.

\bibitem[Bogdan et~al., 2015]{bogdan2015slope}
Bogdan, M., van~den Berg, E., Sabatti, C., Su, W., and Cand{\`e}s, E.~J.
  (2015).
\newblock SlopeÑadaptive variable selection via convex optimization.
\newblock {\em The annals of applied statistics}, 9(3):1103.

\bibitem[Cand{\'e}s and Tao, 2007]{CandesT07}
Cand{\'e}s, E. and Tao, T. (2007).
\newblock The dantzig selector: statistical estimation when p is much larger
  than n.
\newblock {\em The Annals of Statistics}, pages 2313--2351.

\bibitem[Cand{\'e}s and Tao, 2005]{CandesT05}
Cand{\'e}s, E.~J. and Tao, T. (2005).
\newblock Decoding by linear programming.
\newblock {\em Information Theory, IEEE Transactions on}, 51(12):4203--4215.

\bibitem[Efron et~al., 2004]{EfronHJT04}
Efron, B., Hastie, T., Johnstone, I., Tibshirani, R., et~al. (2004).
\newblock Least angle regression.
\newblock {\em The Annals of statistics}, 32(2):407--499.

\bibitem[Fan and Li, 2001]{FanL01}
Fan, J. and Li, R. (2001).
\newblock Variable selection via nonconcave penalized likelihood and its oracle
  properties.
\newblock {\em Journal of the American statistical Association},
  96(456):1348--1360.

\bibitem[Fan et~al., 2015]{fan2015tac}
Fan, J., Liu, H., Sun, Q., and Zhang, T. (2015).
\newblock Tac for sparse learning: Simultaneous control of algorithmic
  complexity and statistical error.
\newblock {\em arXiv preprint arXiv:1507.01037}.

\bibitem[Huang and Zhang, 2012]{HuangZ12}
Huang, J. and Zhang, C.-H. (2012).
\newblock Estimation and selection via absolute penalized convex minimization
  and its multistage adaptive applications.
\newblock {\em Journal of Machine Learning Research}, 13:1809--1834.

\bibitem[Loh and Wainwright, 2015]{loh2015regularized}
Loh, P.-L. and Wainwright, M.~J. (2015).
\newblock Regularized m-estimators with nonconvexity: Statistical and
  algorithmic theory for local optima.
\newblock {\em Journal of Machine Learning Research}, 16:559--616.

\bibitem[Meinshausen and B{\"u}hlmann, 2006]{MeinshausenB06}
Meinshausen, N. and B{\"u}hlmann, P. (2006).
\newblock High-dimensional graphs and variable selection with the lasso.
\newblock {\em The Annals of Statistics}, pages 1436--1462.

\bibitem[Negahban et~al., 2012]{negahban2012unified}
Negahban, S.~N., Ravikumar, P., Wainwright, M.~J., and Yu, B. (2012).
\newblock A unified framework for high-dimensional analysis of m-estimators
  with decomposable regularizers.
\newblock {\em Statistical Science}, pages 538--557.

\bibitem[Nesterov et~al., 2007]{nesterov2007gradient}
Nesterov, Y. et~al. (2007).
\newblock Gradient methods for minimizing composite objective function.

\bibitem[Osborne et~al., 2000a]{OsbornePT00}
Osborne, M.~R., Presnell, B., and Turlach, B.~A. (2000a).
\newblock A new approach to variable selection in least squares problems.
\newblock {\em IMA Journal of Numerical Analysis-Institute of Mathematics and
  its Applications}, 20(3):389--404.

\bibitem[Osborne et~al., 2000b]{OsbornePT00b}
Osborne, M.~R., Presnell, B., and Turlach, B.~A. (2000b).
\newblock On the lasso and its dual.
\newblock {\em Journal of Computational and Graphical statistics},
  9(2):319--337.

\bibitem[Parikh and Boyd, 2013]{parikh2013proximal}
Parikh, N. and Boyd, S. (2013).
\newblock Proximal algorithms, in foundations and trends in optimization.

\bibitem[Ro{\v{c}}kov{\'a} and George, 2016]{RockovaG16}
Ro{\v{c}}kov{\'a}, V. and George, E.~I. (2016).
\newblock The spike-and-slab lasso.
\newblock {\em Journal of the American Statistical Association},
  (just-accepted).

\bibitem[Rudelson and Zhou, 2013]{RudelsonZ13}
Rudelson, M. and Zhou, S. (2013).
\newblock Reconstruction from anisotropic random measurements.
\newblock {\em Information Theory, IEEE Transactions on}, 59(6):3434--3447.

\bibitem[Su and Candes, 2016]{SuC16}
Su, W. and Candes, E. (2016).
\newblock Slope is adaptive to unknown sparsity and asymptotically minimax.
\newblock {\em The Annals of Statistics}, 44(3):1038--1068.

\bibitem[Sun and Zhang, 2013]{SunZ13}
Sun, T. and Zhang, C.-H. (2013).
\newblock Sparse matrix inversion with scaled lasso.
\newblock {\em The Journal of Machine Learning Research}, 14(1):3385--3418.

\bibitem[Tibshirani, 1996]{TibshiraniR96}
Tibshirani, R. (1996).
\newblock Regression shrinkage and selection via the lasso.
\newblock {\em Journal of the Royal Statistical Society. Series B
  (Methodological)}, pages 267--288.

\bibitem[Tropp et~al., 2006]{Tropp06}
Tropp, J. et~al. (2006).
\newblock Just relax: Convex programming methods for identifying sparse signals
  in noise.
\newblock {\em Information Theory, IEEE Transactions on}, 52(3):1030--1051.

\bibitem[van~de Geer and B{\"u}hlmann, 2009]{VanB09}
van~de Geer, S.~A. and B{\"u}hlmann, P. (2009).
\newblock On the conditions used to prove oracle results for the lasso.
\newblock {\em Electronic Journal of Statistics}, 3:1360--1392.

\bibitem[Wainwright, 2009]{Wainwright09}
Wainwright, M.~J. (2009).
\newblock Sharp thresholds for high-dimensional and noisy sparsity recovery
  using-constrained quadratic programming (lasso).
\newblock {\em Information Theory, IEEE Transactions on}, 55(5):2183--2202.

\bibitem[Wang et~al., 2014]{WangLZ14}
Wang, Z., Liu, H., and Zhang, T. (2014).
\newblock Optimal computational and statistical rates of convergence for sparse
  nonconvex learning problems.
\newblock {\em Annals of statistics}, 42(6):2164.

\bibitem[Ye and Zhang, 2010]{YeZ10}
Ye, F. and Zhang, C.-H. (2010).
\newblock Rate minimaxity of the{L}asso and {D}antzig selector for the $\ell_q$
  loss in $\ell_r$ balls.
\newblock {\em Journal of Machine Learning Research}, 11:3481--3502.

\bibitem[Zhang, 2010a]{Zhang10}
Zhang, C.-H. (2010a).
\newblock Nearly unbiased variable selection under minimax concave penalty.
\newblock {\em The Annals of Statistics}, pages 894--942.

\bibitem[Zhang and Huang, 2008]{ZhangH08}
Zhang, C.-H. and Huang, J. (2008).
\newblock The sparsity and bias of the lasso selection in high-dimensional
  linear regression.
\newblock {\em The Annals of Statistics}, pages 1567--1594.

\bibitem[Zhang and Zhang, 2012]{ZhangZ12}
Zhang, C.-H. and Zhang, T. (2012).
\newblock A general theory of concave regularization for high-dimensional
  sparse estimation problems.
\newblock {\em Statistical Science}, pages 576--593.

\bibitem[Zhang, 2010b]{Zhang10-multistage}
Zhang, T. (2010b).
\newblock Analysis of multi-stage convex relaxation for sparse regularization.
\newblock {\em Journal of Machine Learning Research}, 11:1087--1107.

\bibitem[Zhao and Yu, 2006]{ZhaoY06}
Zhao, P. and Yu, B. (2006).
\newblock On model selection consistency of lasso.
\newblock {\em The Journal of Machine Learning Research}, 7:2541--2563.

\bibitem[Zou and Li, 2008]{ZouL08}
Zou, H. and Li, R. (2008).
\newblock One-step sparse estimates in nonconcave penalized likelihood models.
\newblock {\em Annals of statistics}, 36(4):1509.

\end{thebibliography}
